\documentclass{ieeeaccess}
\usepackage{cite}
\usepackage{amsmath,amssymb,amsfonts}
\usepackage{algorithmic}
\usepackage{graphicx}
\usepackage{textcomp}
\usepackage{algorithm}
\usepackage{epstopdf}
\usepackage{mathtools}
\usepackage{mleftright}

\newtheorem{theorem}{Theorem}

\newtheorem{remark}{Remark}

\newtheorem{corollary}{Corollary}
\newcommand*{\QEDB}{\hfill\ensuremath{\square}}%
\usepackage{footnote}
\makesavenoteenv{tabular}
\usepackage{hyperref}
\usepackage{tikz}
\usetikzlibrary{shapes,arrows,fit,positioning,shadows,calc}
\usetikzlibrary{plotmarks}
\usetikzlibrary{decorations.pathreplacing}
\usetikzlibrary{patterns}
\usetikzlibrary{automata}
\usepackage{pgfplots}
\pgfplotsset{compat=newest}
\usepackage{adjustbox}
\usepackage{caption}
  \NewSpotColorSpace{PANTONE}
  \AddSpotColor{PANTONE} {PANTONE3015C} {PANTONE\SpotSpace 3015\SpotSpace C} {1 0.3 0 0.2}
  \SetPageColorSpace{PANTONE}

\def\BibTeX{{\rm B\kern-.05em{\sc i\kern-.025em b}\kern-.08em
    T\kern-.1667em\lower.7ex\hbox{E}\kern-.125emX}}
\begin{document}
\history{Date of publication xxxx 00, 0000, date of current version xxxx 00, 0000.}
\doi{10.1109/ACCESS.2023.3288889}

\title{Randomized Rank-Revealing QLP for Low-Rank Matrix Decomposition}
\author{\uppercase{M. F. Kaloorazi}\authorrefmark{1}, 
\uppercase{K. Liu\authorrefmark{2}, J. Chen\authorrefmark{2}, \IEEEmembership{Senior, IEEE}
R.C. de Lamare\authorrefmark{3}, \IEEEmembership{Senior, IEEE}, S. Rahardja\authorrefmark{2},}
\IEEEmembership{Fellow, IEEE}}
\address[1]{School of Electronic Engineering, Xi’an Shiyou University, Xi'an 710065,  China (e-mail: kaloorazi@xsyu.edu.cn)}
\address[2]{School of Marine Science and Technology, Northwestern Polytechnical University, Xi'an 710072, China}
\address[3]{Center for Telecommunications Studies, Pontifical Catholic University of Rio de Janeiro, Rio de Janeiro 22451-900, Brazil}
\tfootnote{This work was supported in part by the School of Electronic Engineering,
Xi’an Shiyou University Grant No. 103/134010028, by NSFC Grant No. 62171380, and by CNPq and FAPERJ.}

\markboth
{Kaloorazi \headeretal: Randomized Rank-Revealing QLP for Low-Rank Matrix Decomposition}
{Kaloorazi \headeretal: Randomized Rank-Revealing QLP for Low-Rank Matrix Decomposition}

\corresp{Corresponding author: M. F. Kaloorazi (e-mail: kaloorazi@xsyu.edu.cn).}

\begin{abstract}
The pivoted QLP decomposition is computed through two consecutive pivoted QR decompositions. It is an approximation to the computationally prohibitive singular value decomposition (SVD). 
This work is concerned with a \emph{partial} QLP decomposition of matrices through the exploitation of random sampling.
The method presented is tailored for low-rank matrices and called Randomized Unpivoted QLP (RU-QLP). Like pivoted QLP, RU-QLP is rank-revealing and yet it utilizes randomized column sampling and the unpivoted QR decomposition. The latter modifications allow RU-QLP to be highly scalable and parallelizable on advanced computational platforms. 
We provide an analysis for RU-QLP, thereby bringing insights into its characteristics and performance behavior. In particular, we derive bounds in terms of both spectral and Frobenius norms on: i) the rank-revealing property; ii) principal angles between approximate subspaces and exact singular subspaces and vectors; and iii) the errors of low-rank approximations. 
Effectiveness of the bounds is illustrated through numerical tests.  
We further use a modern, multicore machine equipped with a GPU to demonstrate the efficiency of RU-QLP. \textcolor{black}{Our results show that compared to the randomized SVD, RU-QLP achieves a speedup of up to 7.1 and 8.5 times using the CPU and up to 2.3  and 5.8 times using the GPU for the decomposition of dense and sparse matrices, respectively.}
\end{abstract}


\titlepgskip=-15pt

\maketitle

\section{Introduction}
\label{sec:introduction}
\PARstart{L}{ow}-rank matrix approximation and factorization have traditionally been performed using a truncated version of deterministic matrix decomposition methods. These methods include: i) The singular value decomposition (SVD)  \cite{Demmel97}, which constructs the best factorization of a matrix. ii) The pivoted QR \cite{Stewart98}, which, compared to the SVD, is computationally more efficient but less accurate. If computed exactly, due to the column exchange, this decomposition discloses information on the numerical rank (i.e., the gap in the spectrum) of the matrix. iii) The pivoted QLP  \cite[Chapter 5, Section 2.3]{Stewart98}, which is computed by applying two consecutive pivoted QR. This decomposition, in comparison to pivoted QR, constructs a better approximations to both the singular values and singular subspaces of the matrix. 

With modern applications and with the development of advanced computational architectures, these traditional methods, however, are faced with two daunting challenges, which stymie their applicability and hence practicality:
\begin{enumerate}
\item They may need a large number of arithmetic operations.
\item More importantly, they impose high communication costs upon the system, the cost as a results of moving data between the slow and fast memory or, when the processors work together, between processors \cite{Dongarraetal18}. 
\end{enumerate}

In order to address the foregoing bottlenecks, randomized methods \cite{HMT2011}, \cite{MFKC21}, \cite{Gu2015}, \cite{Saibaba19}, \cite{MFKJC19}, \cite{MFKDeJSTSP18}, \cite{MARTINSSON201147}, \cite{MFKC20} have been proposed. They construct approximations to the deterministic decompositions through the utilization of randomization. In comparison to their classical and more established counterparts, randomized methods are arithmetically more efficient, can leverage the parallel structure in modern machines, but less accurate. The latter is justified as optimality in accuracy is not required in many applications.
Randomized methods make use of deterministic decompositions in their computational procedures. Their computations involve three steps: 
\begin{enumerate}
\item Reduce the large dimension of the original input matrix via a random matrix
\item Perform the SVD, the pivoted QR or QLP
\item Construct the factorization through piecing back together the foregoing constituents.
\end{enumerate}
Due to the second step, depending on the dimension of the data and application in hand, randomized techniques may still suffer from the same bottleneck associated with the communication cost. Recently, however, Randomized Unpivoted QLP (RU-QLP), proposed in \cite{MFKC21}, copes with this issue; it makes use of only the unpivoted QR decomposition, which requires less communication among all deterministic decompositions. 

\subsection{Applications}
Randomized methods for low-rank matrix approximation and factorization are applied in signal and image processing, (supervised and unsupervised) machine learning, and modern data analysis. This is in particular motivated by ubiquitousness of large-scale matrices with low-rank structure in diverse range of applications such as weather forecast \cite{Caoetal20}, principal components regression \cite{MorYA19}, latent variable models \cite{FanWZZ21}, image deblurring \cite{JiangCd21}, background modeling \cite{MFKDeTSP18}, deep learning classification and regression \cite{AdelmanLHS21}, large-scale multiple-antenna systems \cite{ShaoLd21}, community detection \cite{Abbe18}, Stokes flow equations \cite{ButtariHLMRW21}, subspace estimation and tracking \cite{CaiKC21}, covariance estimation \cite{CaiZZ16}, Gaussian processes \cite{GardnerPWBW18}, learning mixture models \cite{ChenMVC21}, and latent variable models \cite{HongGBF21}. 

\subsection{Contributions}
\label{subsec_I_Contrib}
Our paper focuses on the RU-QLP decomposition, an approximation of truncated pivoted QLP computed via random sampling and the highly scalable and parallelizable unpivoted QR factorization.
It provides an in-depth analysis of the method, furnishing bounds for i) the rank-revealing property, ii) principal angles between approximate subspaces and exact singular subspaces and vectors, and iii) the errors low-rank approximation. There are, however, differences between the analysis of this paper and the one presented in the original paper. First, the proof technique of our results presented here is different, yet simpler and easier to follow. Our analysis can be viewed as a systematic treatment of randomized methods for low-rank matrix factorization and hence can be adopted by other methods of this class in order to derive error bounds. However, an SVD-like factorization enables deriving more bounds; see Remarks \ref{remSVDlike} and \ref{ReProof}.

Second, all bounds presented here are different, some of which are missing in the original paper.
To be specific, we derive bounds in terms of spectral and Frobenius norm for the following:
\begin{itemize}
\item All principal angles between the approximate left and right subspaces and i) the singular subspaces and ii) the subspaces spanned by the individual singular vectors. 
\item The errors on low-rank approximations formed through the approximate left and right singular vectors. 
\end{itemize}

Third, an empirical evaluation of the derived bounds is given, which provides insight on RU-QLP's characteristics and behavior as well as the tightness of the bounds.  

Finally, we implement RU-QLP and several competing randomized methods on an advanced, multicore machine equipped with a GPU and discuss their performance behaviors.

\subsection{Overview}
In Section \ref{section_PriorWork}, we describe the notations used in our paper, and briefly review the related deterministic and randomized methods. In Section \ref{secRUQLP}, we describe  RU-QLP in detail, its intuition and relation to the classical method of orthogonal iteration, as well as its computational cost. Section \ref{secThAnalysis} establishes a canonical theoretical analysis of RU-QLP, shedding light on its characteristics. Section \ref{secNuSim} presents the runtime results of RU-QLP and several randomized methods implemented on a hybrid GPU-based architecture, as well as those of empirical evaluation of the developed bounds using several matrices. Concluding remarks are given in Section \ref{secConclusion}.

\section{Background}
\label{section_PriorWork}

\subsection{Notation and Conventions}

We consider a real, dense $m\times n$ matrix $\bf A$ with $m\ge n$. The $i$th column of $\bf A$ is denoted by ${\bf a}_i$, and the $i$th largest and the minimum singular value of $\bf A$  by $\sigma_i(\bf A)$ and $\sigma_\text{min}({\bf A})$, respectively. 
The singular values of $\bf A$ are assumed to be in a decreasing order, and $\sigma_k({\bf A})$ and $\sigma_{k+1}({\bf A})$ to be \emph{well-separated}, implying that $\bf A$ has rank $k$. \textcolor{black}{(The rank of a matrix is defined as the number of linearly independent columns or rows in the matrix.)} For a symmetric matrix $\bf A$, the $i$th largest eigenvalue is denoted $\lambda_i({\bf A})$.

$\|{\bf A}\|_2 = \underset{\|{\bf x}\|_2=1}{\text{max}}
\|{\bf Ax}\|_2$ and $\|{\bf A}\|_F = \sqrt{{\sum_{i=1}^{m}\sum_{j=1}^{n}a_{ij}^2}}$ present respectively the spectral norm (or 2-norm) and the Frobenius norm of $\bf A$. The notation ${\|\cdot\|_{2,F}}$ holds for both the spectral and Frobenius norms, and $\|\cdot\|$ for any unitary invariant norm. The dagger $\dagger$ denotes the Moore-Penrose inverse, and ${\bf I}$ refers to the identity matrix whose size is determined by the context. The range and null space of $\bf A$ are denoted by $\mathcal{R}({\bf A})$ and $\mathcal{N}({\bf A})$, respectively. The notation $\texttt{randn}()$ is used to generate random matrices whose entries have a standard Gaussian distribution. $\texttt{orth}({\bf A})$ returns a set of orthonormal basis for the columns of $\bf A$, and $\texttt{qr}({\bf A})$ gives an unpivoted QR factorization. The notation $[{\bf A}]_k$ denotes the optimal rank-$k$ approximation of $\bf A$, formed by the SVD, with respect to 2- and Frobenius norm, and $\mathbb{E}$ denotes the expected value. 

\subsection{Deterministic Methods}
\label{subsec_II_DeterMeth}
\subsubsection{The SVD}
The (reduced) SVD \cite[Section 5.4]{Demmel97} decomposes $\bf A$ into two orthonormal matrices ${\bf U} = [{\bf U}_k \quad {\bf U}_\perp] \in \mathbb R^{m \times n}$ and ${\bf V}=[{\bf V}_k \quad {\bf V}_\perp] \in \mathbb R^{n \times n}$ and one diagonal matrix ${\bf \Sigma} = \text{diag}({\bf \Sigma}_k, {\bf \Sigma}_\perp) \in \mathbb R^{n \times n}$:
\begin{equation}\label{eqSVD}
\begin{aligned}
{\bf A} = {\bf U}{\bf \Sigma}{\bf V}^T,
\end{aligned}
\end{equation}
where the columns of $\bf U$ and $\bf V$ are left and right singular vectors, respectively, and the entries of ${\bf \Sigma} = \text{diag}(\sigma_1, \sigma_2, ..., \sigma_n)$ are the singular values of $\bf A$.

The SVD furnishes information on the four fundamental subspaces $\mathcal{R}({\bf A})= \text{span}\{{\bf u}_1, ..., {\bf u}_k\}$, $\mathcal{N}({\bf A}^T) = \text{span}\{{\bf u}_{k+1}, ..., {\bf u}_m\}$, $\mathcal{N}({\bf A}) = \text{span}\{{\bf v}_{k+1}, ..., {\bf v}_n\}$, and $\mathcal{R}({\bf A}^T) = \text{span}\{{\bf v}_1, ..., {\bf v}_k\}$. The SVD constitutes an optimal choice in constructing a rank-$r$ approximation to ${\bf A}$. Let ${\bf A}_r = \sum_{i=1}^{r}{\sigma_i{\bf u}_i{\bf v}_i^T}$ be the rank-$r$ approximation. Then \cite[Chapter 1, Theorem 4.32]{Stewart98}:
\begin{equation}\notag
\begin{aligned}
\|{\bf A} - {\bf A}_r\|_2 = & \underset{\text{rank}({\bf M})\le r}{\text{min}}
&& \|{\bf A} - {\bf M}\|_2 = \sigma_{r+1},
\end{aligned}
\end{equation}
\begin{equation}\notag
\begin{aligned}
\|{\bf A} - {\bf A}_r\|_F = & \underset{\text{rank}({\bf M})\le r}{\text{min}}
&& \|{\bf A} - {\bf M}\|_F = \sqrt{\sum_{i=r+1}^{n}{\sigma_i^2}}.
\end{aligned}
\end{equation}

The SVD is commonly computed as follows:
\begin{enumerate}
\item Reduce $\bf A$ to bidiagonal form ${\bf B}={\bf U}_1^T{\bf A}{\bf V}_1$, through Householder reflections ${\bf U}_1$ and ${\bf V}_1$. 
\item Compute the SVD of ${\bf B}={\bf U}_2{\bf \Sigma}{\bf V}_2^T$, primarily through the QR algorithm. 
\item Combine the forgoing decompositions to obtain 
${\bf A}=({\bf U}_1{\bf U}_2){\bf \Sigma}({\bf V}_1{\bf V}_2)^T$, where ${\bf U}={\bf U}_1{\bf U}_2$ and ${\bf V}={\bf V}_1{\bf V}_2$.
\end{enumerate}

\subsubsection{Unpivoted and pivoted QR decompositions} 
The (reduced) unpivoted QR decomposition of ${\bf A} = {\bf Q}{\bf R}$ gives ${\bf Q}\in \mathbb R^{m \times n}$ whose columns are orthonormal and an upper triangular factor ${\bf R} \in \mathbb R^{n \times n}$ \cite[Section 2.3]{Bjorck15}.

The pivoted QR \cite[Chapter 5, Section 2.1]{Stewart98}  factorizes $\bf A$ into an exchange matrix ${\bf \Pi}_\text{p}$, an orthonormal matrix ${\bf Q}_\text{p}$, and an upper triangular matrix ${\bf R}_\text{p}$:
\begin{equation}\notag
{\bf A}{\bf \Pi}_\text{p} = {\bf Q}_\text{p}{\bf R}_\text{p} = [{\bf Q}_1 \quad {\bf Q}_2]
\begin{bmatrix}
{\bf R}_{11} & {\bf R}_{12}  \\
{\bf 0} & {\bf R}_{22}
\end{bmatrix}.
\end{equation} 

The difference between unpivoted and pivoted QR is that in the computation of the latter the columns with largest 2-norm are exchanged with other columns before the reduction (through Householder reflectors) proceeds.
The columns of ${\bf Q}_1 \in \mathbb R^{m \times k}$ and ${\bf Q}_2 \in \mathbb R^{m \times n-k}$ span respectively $\mathcal{R}({\bf A})$ and $\mathcal{N}({\bf A}^T)$, the diagonals of ${\bf R}$ are approximations to $\sigma_i({\bf A})$, and ${\bf R}_{11} \in \mathbb R^{k \times k}$ is well-conditioned. If   ${\bf \Pi}_\text{p}$ is chosen carefully, the decomposition is called ``rank-revealing" QR \cite{Chan87, GUEis96}, and the blocks ${\bf R}_{11}$ and ${\bf R}_{22}$ satisfy:
\begin{equation}\label{eqDeterRankR}
\begin{aligned}
\sigma_\text{min}({\bf R}_{11}) & \ge \dfrac{\sigma_k({\bf A})}{x(n, k)},\\
\sigma_1({\bf R}_{22}) & \le \sigma_{k+1}({\bf A})t(n,k),
\end{aligned}
\end{equation}
where $x(n,k)$ and $t(n,k)$ are low degree polynomials in $n$ and $k$. Pivoted QR, in spite of its computational efficiency in comparison to the SVD, gives fuzzy singular-value estimates, and further does not furnish orthogonal bases for $\mathcal{R}({\bf A}^T)$ and $\mathcal{N}({\bf A})$ explicitly. 

\subsubsection{UTV decompositions}
These rank-revealing decompositions \cite[Chapter 5, Section 4]{Stewart98}  factorize the matrix $\bf A$ as:
\begin{equation}\notag
{\bf A}={\bf U}_\text{t}{\bf T}_\text{t}{\bf V}_\text{t}^T,
\end{equation} 
where ${\bf U}_\text{t}$ and ${\bf V}_\text{t}$ are orthogonal, and ${\bf T}_\text{t}$ is upper or lower triangular. There exist two primary stages in the computation of UTVs: an initial unpivoted QR factorization followed by a rank-revealing step, or deflation steps, in which the largest, or the smallest singular values are extracted one at a time. The UTVs provide information on $\mathcal{R}({\bf A}^T)$ and $\mathcal{N}({\bf A})$ as well.

\subsubsection{Pivoted QLP decomposition} The pivoted QLP \cite[Chapter 5, Section 2.3]{Stewart98} is viewed as an approximate SVD. It is formed by applying two consecutive pivoted QR factorizations, first on $\bf A$, then on the transpose of the $R$ factor. Specifically
\begin{equation}\label{eqQLP}
 {\bf A}{\bf \Pi}_\text{p} = {\bf Q}_\text{p}{\bf R}_\text{p}, \quad {\bf R}_\text{p}^T\dot{\bf \Pi} = \dot{\bf P}\dot{\bf L}^T,
\end{equation}
which gives ${\bf A} = {\bf Q}_\text{p}\dot{\bf \Pi}\dot{\bf L}\dot{\bf P}^T{\bf \Pi}_\text{p}^T$. The matrices ${\bf Q}_\text{p}\dot{\bf \Pi}$ and ${\bf \Pi}_\text{p}\dot{\bf P}$ are orthogonal, providing bases for the spaces spanned by the columns and rows of $\bf A$, respectively. $\dot{\bf L}$ is lower triangular and, as demonstrated by Stewart,  its diagonals (the L-values) give better approximations to $\sigma_i({\bf A})$ than the diagonals of ${\bf R}_\text{p}$ (the R-values). In computing pivoted QLP, the first exchange matrix ${\bf \Pi}_\text{p}$ is crucial, whereas the second one $\dot{\bf \Pi}$ is not always necessary; see also \cite{HuckabyChan03}. This principle has been leveraged for developing RU-QLP, as will be discussed later.

\subsection{Shortcomings of Deterministic Methods}
\subsubsection{Arithmetic cost}
Arithmetic operations required by deterministic methods to factorize $\bf A$ are of order $mn^2$. This is the cost of a \emph{full} factorization, giving all relevant information. Considering large-scale matrices in modern applications, this is obviously prohibitively expensive. However, if a rank $k$ factorization is desired, deterministic methods need $\mathcal{O}(mnk)$ operations to give the truncated version. These methods, however, need to repeatedly \emph{access} the data, which brings us to a more important cost associated with decomposition methods on modern computers, namely, the communication cost.  

\subsubsection{Communication cost}
\label{subsubComCost}
The cost of transferring data between different processors or between different levels of the memory hierarchy is defined as the communication cost. On advanced parallel machines, it is far more expensive than the arithmetic cost in terms of time as well as energy consumption \cite{DoGH20}. The communication cost is associated with the use of level-1, 2 and 3 BLAS (the Basic Linear Algebra Subprograms) routines \cite{ButtariLKD08, Dongarraetal18}. Memory-bound level-1 and level-2 BLAS routines cannot attain high performance on modern computers. However, level-3 BLAS routines are CPU-bound, which enable harnessing the parallel architecture of modern machines.

To approximate the SVD or/and UTV decompositions, Krylov subspace methods, such as the Lanczos algorithm, are used \cite{FierroHanHan99}. A large portion of operations of such methods are performed in level-1 and level-2 BLAS. Pivoted QR needs to compute the 2-norm of the matrix columns and swap them, and hence nearly half of its operations are in level-3 BLAS. While, most operations of the unpivoted QR decomposition are in level-3 BLAS, meaning that it is \emph{communication-friendly}, and can be efficiently implemented on parallel machines. Recently, randomized techniques have been used to develop block Krylov subspace methods \cite{MuscoM15,YuanGL18,Tropp22} and ``communication-avoiding" pivoted QR decompositions \cite{MartinssonV16, DuerschGu20}.

\subsection{Randomized Methods}
\label{subsec_II_RandMeth}

Methods based on randomized sampling offer efficient and compatible approximations to traditional decompositions; they are efficient in arithmetic operations, and can harness parallel architectures of advanced computing devices. Efficiency is attained by reducing the dimension of the input matrix, thus keeping only \emph{important} features. This in turn leads to less accurate approximations. However, there are some techniques, such as the power method, that substantially improve the quality of approximations. The above factors make the randomized methods very attractive for processing low-rank matrices. 

The line of research on randomized methods that culminated in RU-QLP \cite{MFKC21} started with \cite{MARTINSSON201147}, and hence we discuss such methods. We refer to \cite[Section 2]{HMT2011} and \cite[Section 2.4]{GitM16}, \cite{TrYUC17}, \cite{ChepurkoCKW22} for other classes of randomized methods. The general strategy underlying the randomized methods involves the following steps:
\begin{itemize}
\item [] Step 1: The input matrix is transformed to a lower-dimensional space by utilizing a random matrix. This renders a smaller matrix formed by the linear combinations of rows or columns. Then orthonormal bases are obtained through the Gram–Schmidt process \cite[Section 2.3.4]{Bjorck15}, \cite{Barlow19} or Householder reflections \cite{GunterVan05, JoffrainTQVV06}.  
\item [] Step 2: Orthonormal bases are (left or right) multiplied to the original matrix, and then the SVD or pivoted QR or QLP of the reduced-size matrix is computed.
\item [] Step 3: The orthogonal and diagonal or (upper/lower) triangular components are combined forming the final approximation.
\end{itemize}

The random matrix used in Step 1 is usually standard Gaussian; see \cite[Section 3.9]{TrYUC17} for a discussion about other choices. The major difference in computational procedure of the randomized methods, however, appears in Step 2. The works in \cite{MARTINSSON201147, HMT2011} make use of the SVD, where the resulting method of the latter is called the ``randomized SVD" (R-SVD). The error analyses presented though are different. Gu \cite{Gu2015} supplanted the full SVD with a truncated version in Step 2. His error bounds explicitly contain the oversampling parameter. The work in \cite{MFKDeTSP18} presents a two-sided R-SVD. It uses the truncated SVD, and the proof techniques are based on \cite{Gu2015}. Saibaba \cite{Saibaba19} presents a new analysis for R-SVD. In particular, he develops a set of bounds for canonical angles between approximate subspaces and exact singular subspaces and individual singular vectors. The works in  \cite{MFKDeJSTSP18}, and \cite{MFKJC19} are based on the two-sided R-SVD, but the difference is that the former employs pivoted QR, and the latter pivoted QLP. The work in \cite{WuX20} presents a randomized QLP method by supplanting the SVD in R-SVD \cite{HMT2011} with pivoted QLP. The method in \cite{MFKC20} makes use of pivoted QLP, and bounds for estimated singular values, canonical angles between the corresponding subspaces, and the errors of low-rank approximations are provided. 

Through utilization of the SVD, pivoted QR or QLP, the foregoing methods each gives three factors: two matrices with orthonormal columns and one matrix with diagonal or triangular structure. The execution of these methods, as expounded earlier, may bring substantial communication cost when processing large-scale matrices, due to the use of level-1 and level-2 BLAS routines. Whereas, RU-QLP \cite{MFKC21} uses the unpivoted QR factorization in Step 2, whose operations are almost entirely in level-3 BLAS. This makes it possible for RU-QLP to leverage modern architectures and therefore to be implemented more efficiently compared to any other randomized method. Below we describe RU-QLP, and the intuition behind its development. We also provide some theory explaining the relationship between RU-QLP and the classical Orthogonal Iteration \cite[Section 8.2.4]{GolubVanLoan96}.     

We add that there is a class of randomized block Krylov subspace methods, e.g., \cite{WangZZ15, DrineasIKM18}, that provide low-rank approximations; they approximate the left invariant subspace of a matrix and use its orthogonal projection to construct the approximation. They do not provide information on the right invariant subspace nor the singular values. If only a low-rank approximation is desired, these methods can also be used.

\section{Randomized Unpivoted QLP}
\label{secRUQLP}
\subsection{Description} Given the matrix ${\bf A}$ with rank $k \ge 1$, and $d=k+p<n$, where $p$ is an oversampling parameter, RU-QLP is computed as follows:
\begin{enumerate}
\item [1.] Form a standard Gaussian matrix ${\bf \Phi} \in \mathbb R^{m \times d}$.
\item [2.] Compress the matrix $\bf A$ via ${\bf \Phi}$ to obtain ${\bf A}^T{\bf \Phi}$.
\item [3.] Form a matrix of orthonormal bases $\bar{\bf P} = \texttt{orth}({\bf A}^T{\bf \Phi})$. 
\item [4.] Compute the matrix product ${\bf A}\bar{\bf P}$, and carry out two unpivoted QR factorizations: \\
$[{\bf Q}, {\bf R}] = \texttt{qr}({\bf A}\bar{\bf P})$, and $[\widetilde{\bf P}, \widetilde{\bf R}] = \texttt{qr}({\bf R}^T)$.
\item [5.] Construct the rank-$d$ approximation $\hat{\bf A}={\bf Q}{\bf L}{\bf P}^T$. Here ${\bf Q}\in \mathbb R^{m \times d}$ approximates $\mathcal{R}({\bf A})$, ${\bf P}\triangleq \bar{\bf P}\widetilde{\bf P}\in \mathbb R^{n \times d}$ approximates $\mathcal{R}({\bf A}^T)$, and ${\bf L} \triangleq \widetilde{\bf R}^T\in \mathbb R^{d \times d}$ is lower triangular. Its diagonal elements (L-values) approximate the first $d$ singular values of $\bf A$. Its leading block, which is $k\times k$, reveals the rank $k$ of $\bf A$. RU-QLP is presented in Algorithm \ref{AlgRUQLP}.
\end{enumerate}

\begin{algorithm}
	\caption{Randomized Unpivoted QLP (RU-QLP)}
	\begin{algorithmic}[1]
		\renewcommand{\algorithmicrequire}{\textbf{Input:}}
		\REQUIRE
		~~ ${\bf A}\in \mathbb R^{m \times n}$, $k$, and $d=k+p$.
		\renewcommand{\algorithmicrequire}{\textbf{Output:}}
		\REQUIRE
		 ${\bf Q} \in \mathbb R^{m \times d}$ and ${\bf P} \in \mathbb R^{n \times d}$ with orthonormal columns, and lower triangular ${\bf L}\in \mathbb R^{d \times d}$ that form an approximation $\hat{\bf A} = {\bf Q}{\bf L}{\bf P}^T$.
		\STATE \textbf{Function} \text{{RUQLP}}$({{\bf A}, d})$
		\STATE \quad ${\bf \Phi}=\texttt{randn}(m, d)$
		\STATE \quad $\bar{\bf P} = \texttt{orth}( {\bf A}^T{\bf \Phi})$
		\STATE \quad $[{\bf Q}, {\bf R}] = \texttt{qr}( {\bf A}\bar{\bf P})$
		\STATE \quad $[\widetilde{\bf P}, \widetilde{\bf R}]=\texttt{qr}({\bf R}^T)$ 
         \STATE \quad \textbf{return} ${\bf Q}$; ${\bf P} \triangleq \bar{\bf P}\widetilde{\bf P}; {\bf L} \triangleq \widetilde{\bf R}^T$
		\STATE \textbf{End function}
	\end{algorithmic}
\label{AlgRUQLP}
\end{algorithm}

\subsection{Improvement and numerical stability} 
\label{subsubImOr}
To ameliorate the approximation quality, particularly when the singular values of the input matrix do not decay relatively fast, we use the power iteration technique: it replaces ${\bf A}^T$ in Step 3 of Algorithm \ref{AlgRUQLP} with ${\bf A}_q=({\bf A}^T{\bf A})^q{\bf A}^T$. Here $q\ge 1$ is the power method factor. The power iteration also affects the convergence of the approximate left and right invariant subspaces: the approximate subspaces converge to invariant subspaces at a rate proportional to $(\sigma_{k+1}/\sigma_k)^q$. 

There is, however, a concern regarding the computation of ${\bf A}_q$ in floating point arithmetic: it is prone to rounding errors, which will lead to loss of information associated with some small singular values. To be specific, considering the machine precision $\epsilon_\text{machine}$, the singular components less than $\sigma_1 \epsilon_\text{machine}^{1/(2q+1)}$ will be lost; see \cite[Section V-D]{MFKC21} for an example. Thus, to have better numerical accuracy, Algorithm \ref{AlgPIOrth} needs to be utilized to compute $\bar{\bf P}$ in Algorithm \ref{AlgRUQLP}, where orthonormalization of the columns of the sample matrices are carried out \cite{HMT2011, Gu2015, MFKC21}. Since orthonormalization incurs additional cost, it is advised to be used once in every few iterations. This helps balance numerical stability and efficiency. 

\begin{algorithm}
	\caption{Power Iteration and Orthorgonalization}
	\begin{algorithmic}[1]
		\renewcommand{\algorithmicrequire}{\textbf{Input:}}
		\REQUIRE
		~ ${\bf A}\in \mathbb R^{m \times n}$,  ${\bf \Phi} \in \mathbb R^{m \times d}$, and $q\ge 1$.
		\renewcommand{\algorithmicrequire}{\textbf{Output:}}
		\REQUIRE  
		$\bar{\bf P} \in \mathbb R^{n \times d}$ with orthonormal columns.
			\STATE \textbf{Function} \text{PI\_ORTH}{(${\bf A}, {\bf \Phi}, q$)}
		\STATE \quad $\bar{\bf P} = \texttt{orth}( {\bf A}^T{\bf \Phi})$
		\STATE \quad \textbf{for} $i=1,..., q$
		\STATE \quad \quad \quad ${\bf B} = {\bf A}\bar{\bf P}$;     $\bar{\bf P} = \texttt{orth}({\bf B})$
		\STATE \quad \quad \quad ${\bf B} = {\bf A}^T\bar{\bf P}$; $\bar{\bf P} = \texttt{orth}({\bf B})$
		\STATE \quad \textbf{end for}
\STATE \quad \textbf{return} $\bar{\bf P}$
		\STATE \textbf{End function}
	\end{algorithmic}
\label{AlgPIOrth}
\end{algorithm}

\subsection{Intuition and relation to Orthogonal Iteration}
The pivoted QLP is computed by applying two consecutive pivoted QR factorizations, first on $\bf A$, then on ${\bf R}^T$. The first pivoting (i.e., column exchanges) is crucial, but the main purpose of the second pivoting is to ensure that the L-values are arranged in a decreasing order. In another words, if the space spanned by the first $k$ columns of $\bf Q$ \eqref{eqQLP} furnishes a good approximation to the left invariant subspace of $\bf A$ and thus giving a well-conditioned leading block of order $k$ in the $R$ factor that reveals the rank of $\bf A$, then the second pivoting can be circumvented. To connect it to RU-QLP, we show that the first $k$ columns of $\bf Q$ obtained by the unpivoted QR factorization on ${\bf A}\bar{\bf P}$ in fact gives a good approximation to $\mathcal{R}({\bf A})$. We further show that the leading block of the $R$ factor reveals the rank of $\bf A$. This is in essence due to the fact that $\bar{\bf P}$ approximates  $\mathcal{R}({\bf A}^T)$, and a high-quality approximation is given when the power iteration is used.  

Orthogonal Iteration \cite[Section 8.2.4]{GolubVanLoan96} is a generalization of the power method. It is used to compute the dominant invariant subspaces of a matrix. Let ${\bf A} \in \mathbb R^{n \times n}$ and ${\bf Q}^{(0)}$ be $m \times k$ with orthonormal columns. The following orthogonal iteration procedure produces a sequence of  ${\bf Q}^{(t)}$:
\begin{itemize}
\item [] \textbf{for} $t=1, 2, ...$
\item [] \quad  ${\bf K}^{(t)}={\bf A}{\bf Q}^{(t-1)}$
\item [] \quad  ${\bf K}^{(t)}={\bf Q}^{(t)}{\bf R}^{(t)}$
\item [] \textbf{end}
\end{itemize}

Under the assumption that the eigenvalues of $\bf A$ are arranged in decreasing order, and that the largest $k$ eigenvalues are separated from the remainder of the spectrum, $\mathcal{R}({\bf Q}^{(t)})$ converges to the dominant left invariant subspace of $\bf A$ as $t \rightarrow \infty$. If $\bf A$ is real, then the diagonal entries of $\bf R$ converge to the dominant eigenvalues. 

Unlike the orthogonal iteration method, we construct $\bar{\bf P}$ through random sampling of ${\bf A}$'s rows. However, we observe that the $Q$ and $R$ factors of the QR factorization of ${\bf A}\bar{\bf P}$ reveal similar information as those of orthogonal iteration: an approximate basis for the invariant subspace and approximate singular values, respectively. These results are enhanced when the power iteration is used to form $\bar{\bf P}$. It is this principle that is leveraged in developing RU-QLP. 

\begin{remark}
\label{remSVDlike}
It is possible to turn RU-QLP to an SVD-like decomposition by computing an SVD of ${\bf R}=\bar{\bf U}\bar{\bf \Sigma}\bar{\bf V}^T$, hence $\hat{\bf A}=({\bf Q}\bar{\bf U})\bar{\bf \Sigma}(\bar{\bf P}\bar{\bf V})^T$. This, however, can be done entirely using unpivoted QR, due to the method of QR-based Dynamically Weighted Halley (QDWH)-SVD \cite{NakatsukasaH13}. On $\bf R$, this method needs up to $52d^3$ operations, instead of $4d^3/3$ or $2d^3$ operations for an unpivoted QR, which clearly demands more resources.  
\end{remark}

\subsection{Computational cost} The number of arithmetic operations required to compute RU-QLP is as follows: Forming the matrix $\bf \Phi$ $\mathcal{O}(md)$;
	 forming ${\bf A}^T\bf \Phi$ $\mathcal{O}(mnd)$;
	 computing $\bar{\bf P}$ $\mathcal{O}(nd^2)$;
	 forming ${\bf A}\bar{\bf P}$ $\mathcal{O}(mnd)$;
	 computing ${\bf Q}$ and ${\bf R}$ $\mathcal{O}(md^2)$;
	 computing $\widetilde{\bf P}$ and $\widetilde{\bf R}$ $\mathcal{O}(d^3)$;
	 forming $\bf P$ $\mathcal{O}(nd^2)$.

The dominant cost is $\mathcal{C}= \mathcal{O}(mnd)$. If the input matrix $\bf A$ is stored out-of-core, RU-QLP is computed by two passes over $\bf A$. If the power iteration is used, as to ameliorate the factorization quality, RU-QLP needs $2q+2$ passes over $\bf A$, and $(q+1)\mathcal{C}$ operations. In addition, if $\bf A$ is sparse with $s$ non-zero entries, computing RU-QLP costs $\mathcal{O}(sd)$.

As discussed earlier, in modern computer architectures, the cost of communication \cite{Dongarraetal18,DoGH20} in performing decompositional methods dominates the arithmetic cost. This is associated with level-1, 2 and 3 BLAS routines. \textcolor{black}{The methods that are rich in level-3 BLAS, that is, most of their operations can be cast as matrix-matrix multiplications, attain higher performance, due to the fact that the movement of their data can be minimized.} In contrast to the SVD, and pivoted QR, unpivoted QR can be computed almost entirely using level-BLAS 3 routines, thereby lending itself much easier to parallel implementation. This in turn makes RU-QLP the fastest randomized method, as will be shown in Section \ref{subRuntime}.

\section{Theoretical Analysis}
\label{secThAnalysis}

We first define a few terms that are used in our results. Let $\bf U$ be the matrix of left singular vector as in \eqref{eqSVD}. Let $\widehat{\bf \Phi}_1 \in \mathbb R^{k \times d}$ and $\widehat{\bf \Phi}_2 \in \mathbb R^{(m-k) \times d}$ be defined as follows:
\begin{equation}
\begin{aligned}
{\bf U}^T{\bf \Phi} =
\begin{bmatrix}
{\bf U}_k^T{\bf \Phi} \\
{\bf U}_\perp^T{\bf \Phi}
\end{bmatrix} \triangleq
\begin{bmatrix}
\widehat{\bf \Phi}_1\\
\widehat{\bf \Phi}_2
\end{bmatrix},
\label{}
\end{aligned}
\end{equation}
and $\delta_i= \frac{\sigma_{k+1}}{\sigma_i}$, for $i=1,..., k,$ and $\gamma=\sigma_n/\sigma_1$. Let further ${\bf A}\bar{\bf P}$ and its QR factorization be written as:
\begin{equation}\label{eqAPbarPart}
{\bf A}\bar{\bf P} = [{\bf A}\bar{\bf P}_1\quad {\bf A}\bar{\bf P}_2]=[{\bf Q}_1\quad {\bf Q}_2]
\begin{bmatrix}
{\bf R}_{11} & {\bf R}_{12}  \\
{\bf 0} & {\bf R}_{22}
\end{bmatrix},
\end{equation}
where $\bar{\bf P}_1 \in \mathbb R^{n \times k}$ and ${\bf Q}_1 \in \mathbb R^{m \times k}$, giving
\begin{equation}\label{eqAPbarPartSep}
\begin{aligned}
& {\bf A}\bar{\bf P}_1 = {\bf Q}_1{\bf R}_{11},\\
& {\bf A}\bar{\bf P}_2 = {\bf Q}_1{\bf R}_{12} + {\bf Q}_2{\bf R}_{22}.
\end{aligned}
\end{equation}

\subsection{${\bf R}$ Reveals the Numerical Rank of {\bf A}}

We first show that the $R$ factor \eqref{eqAPbarPart} generated by an unpivoted QR factorization in the RU-QLP computational procedure reveals the rank of $\bf A$. Our results parallel those of deterministic rank-revealing decompositions \eqref{eqDeterRankR}.    

\begin{theorem}\label{ThSigmakR22}
Let ${\bf A} \in \mathbb R^{m \times n}$ be a rank-$k$ matrix, with $m\ge n$ and an SVD as in \eqref{eqSVD}. Let ${\bf R}$ be generated by RU-QLP and partitioned as in \eqref{eqAPbarPart}, Then, for $i=1,..., k,$
\begin{equation}\label{eqThSigmak}
\sigma_i \ge {\sigma}_i ({\bf R}_{11}) \ge \frac{\sigma_i}{\sqrt{1 + \delta_i^{4q+4}\|\widehat{\bf \Phi}_2\widehat{\bf \Phi}_1^\dagger\|_2^2}}.
\end{equation}
\begin{equation}\label{eqThR22}
\begin{aligned}
\|{\bf R}_{22}\|_{2,F} \le \Bigg(1 + \dfrac{\delta_k^{2q+1}\|\widehat{\bf \Phi}_{2}\widehat{\bf \Phi}_{1}^\dagger\|_2}{1 + \gamma^{4q+4}\|\widehat{\bf \Phi}_{2}\widehat{\bf \Phi}_{1}^\dagger \|_2^2}\Bigg)\|{\bf \Sigma}_\perp\|_{2,F}.
\end{aligned}
\end{equation}
\end{theorem}

\textit{Proof.} The proof is given in Appendix \ref{Proo_R22bound}.

Theorem \ref{ThSigmakR22} makes it explicit that the convergence of ${\sigma}_k ({\bf R}_{11})$ and $\|{\bf R}_{22}\|_{2,F}$ is governed by the ratio $\delta_k= \frac{\sigma_{k+1}}{\sigma_k}$, or simply the gap in the spectrum. Provided that  ${\sigma_k} \gg \sigma_{k+1}$ fast convergence is expected.

\begin{corollary} \label{CoroSigmaD}
Under the hypotheses of Theorem \ref{ThSigmakR22}, let further ${\bf L}$ be the middle factor in the RU-QLP decomposition with ${\bf L}_{11}$ as its top left $k\times k$ block. Then, for $i=1,..., k,$
\begin{equation}\label{eqCorL11}
\sigma_i \ge {\sigma}_i ({\bf L}) \ge {\sigma}_i ({\bf L}_{11}) \ge \frac{\sigma_i}{\sqrt{1 + \delta_i^{4q+4}\|\widehat{\bf \Phi}_2\widehat{\bf \Phi}_1^\dagger\|_2^2}}.
\end{equation}
\begin{equation}\label{eqCorL22}
\begin{aligned}
\|{\bf L}_{22}\|_{2,F} \le \Bigg(1 + \dfrac{\delta_k^{2q+1}\|\widehat{\bf \Phi}_{2}\widehat{\bf \Phi}_{1}^\dagger\|_2}{1 + \gamma^{4q+4}\|\widehat{\bf \Phi}_{2}\widehat{\bf \Phi}_{1}^\dagger \|_2^2}\Bigg)\|{\bf \Sigma}_\perp\|_{2,F}.
\end{aligned}
\end{equation}
\end{corollary}

The second relation in \eqref{eqCorL11} follows because ${\bf L}_{11}$ is a submatrix of ${\bf L}$, and third relation due to \cite[equation (2.7)]{MathiasSte93}. Tighter and more complicated bounds for ${\sigma}_i ({\bf L}_{11})$ can be obtained using the techniques in  \cite{HuckabyChan03}. The relation in \eqref{eqCorL22} follows due to \cite[equation (2.9)]{MathiasSte93}.

\subsection{Bounds for Principal Angles}
\label{subBoPABS}
Closeness of any two subspaces are measured by means of canonical or principal angles \cite{BjorckGolub73, ZhuK13}. The ranges of $\bf Q$ and $\bf P$ built by RU-QLP approximate the ranges of ${\bf U}_k$ and ${\bf V}_k$, respectively. Let $\theta_i = \angle(\mathcal{R}({\bf Q}), \mathcal{R}({\bf U}_k))$ and $\phi_i = \angle(\mathcal{R}({\bf P}), \mathcal{R}({\bf V}_k))$ be the angles between the approximate and exact subspaces. The following results show how accurate the approximations are.

\begin{theorem}\label{ThCanABk}
Let ${\bf A} \in \mathbb R^{m \times n}$ be a rank-$k$ matrix, with $m\ge n$ and an SVD as in \eqref{eqSVD}. Let $\bf A$ have a RU-QLP decomposition as expounded in Section \ref{secRUQLP}. Then, for $i= 1, ..., k$,
\begin{equation}\label{eqThCanABkSinT}
\begin{aligned}
\text{sin}\theta_i \le 
\frac{\delta_i^{2q+2}{\|{\widehat{\bf \Phi}_2}\widehat{\bf \Phi}_1}^\dagger\|_2} {\sqrt{
{1 + \delta_i^{4q+4}{\|{\widehat{\bf \Phi}_2}\widehat{\bf \Phi}_1}^\dagger\|_2^2}}},
\end{aligned} 
\end{equation}
\begin{equation}\label{eqThCanABkSinP}
	\begin{aligned}
	\text{sin}\varphi_i \le 
	\frac{\delta_i^{2q+1}{\|{\widehat{\bf \Phi}_2}\widehat{\bf \Phi}_1}^\dagger\|_2} {\sqrt{{1 + \delta_i^{4q+2}{\|{\widehat{\bf \Phi}_2}\widehat{\bf \Phi}_1}^\dagger\|_2^2}}}.
	\end{aligned} 
	\end{equation}
\end{theorem}

\textit{Proof.} The proof is presented in Appendix \ref{Proo_CanABk}.

Theorem \ref{ThCanABk} makes two points clear. First, $\theta_i$ and $\phi_i$ approach zero at a rate proportional to $(\sigma_{k+1}/\sigma_i)^q$.  Given $\sigma_i> \sigma_{k+1}$, as the power factor $q$ increases, the approximate subspaces become more accurate. Through numerical tests we show that this holds true. Second, $\theta_i$ is smaller than $\phi_i$. This is because the computation of $\bf Q$ requires $\bar{\bf P}$, which forms ${\bf P}$, and hence one more step of iteration.

\begin{corollary}
Under the hypotheses of Theorem \ref{ThCanABk}, we have
\begin{equation}\notag
\begin{aligned}
&  \text{cos}\theta_i \ge 
   \frac{1}{\sqrt{{1 + \delta_i^{4q+4}{\|{\widehat{\bf \Phi}_2}\widehat{\bf \Phi}_1}^\dagger\|_2^2}}}, \hspace{2mm}
\text{tan}\theta_i \le 
   \delta_i^{2q+2}{\|{\widehat{\bf \Phi}_2}\widehat{\bf \Phi}_1}^\dagger\|_2,\\
& 	\text{cos}\varphi_i \ge 
	\frac{1}{\sqrt{{1 + \delta_i^{4q+2}{\|{\widehat{\bf \Phi}_2}\widehat{\bf \Phi}_1}^\dagger\|_2^2}}}, \hspace{2mm}
  \text{tan}\varphi_i \le 
	\delta_i^{2q+1}{\|{\widehat{\bf \Phi}_2}\widehat{\bf \Phi}_1}^\dagger\|_2.
\end{aligned} 
\end{equation}
\end{corollary}

\begin{remark}
The lower bounds for the cosines can be obtained directly using  \cite[equation 12]{BjorckGolub73}. The proofs follow along the same lines as those of the sines given in Appendix \ref{Proo_CanABk}.  
\end{remark}


The following theorem bounds from above in the 2- and Frobenius norms the largest principal angles between the approximate and exact subspaces. This is called the \emph{distance} between the two subspaces \cite[Section 2.6.3]{GolubVanLoan96}.   
\begin{theorem}\label{ThCanASMax}
Under the hypotheses of Theorem \ref{ThCanABk}, we have
\begin{equation}\notag 
\begin{aligned}
& \|\text{sin}\angle(\mathcal{R}({\bf Q}), \mathcal{R}({\bf U}_k))\|_{2,F} \le  \dfrac{\delta_k^{2q+1}\|\widehat{\bf \Phi}_{2}\widehat{\bf \Phi}_{1}^\dagger \|_2\|{\bf \Sigma}_\perp \|_{2,F}}{\sigma_k\sqrt{1 + \gamma^{4q+4}\|\widehat{\bf \Phi}_{2}\widehat{\bf \Phi}_{1}^\dagger \|_2^2}},\\
& 
\|\text{sin}\angle(\mathcal{R}({\bf P}), \mathcal{R}({\bf V}_k))\|_{2,F} \le  \dfrac{\delta_k^{2q}\|\widehat{\bf \Phi}_{2}\widehat{\bf \Phi}_{1}^\dagger \|_2\|{\bf \Sigma}_\perp \|_{2,F}}{\sigma_k\sqrt{1 + \gamma^{4q+2}\|\widehat{\bf \Phi}_{2}\widehat{\bf \Phi}_{1}^\dagger \|_2^2}}.
\end{aligned} 
\end{equation}
\end{theorem}

\textit{Proof.} The proof is presented in Appendix \ref{ProoThCanASMax}.


The theorem that follows furnishes upper bounds for the sine of principal angles between subspaces spanned by individual singular vectors and the approximate subspaces. 

\begin{theorem}\label{ThCanASIn}
Under the hypotheses of Theorem \ref{ThCanABk}, we have  for $i=1, ..., k$,
\begin{equation}
\begin{aligned}\notag
& \text{sin}\angle(\mathcal{R}({\bf Q}), \mathcal{R}({\bf u}_i)) \le \dfrac{\delta_i^{2q+2}\|\widehat{\bf \Phi}_{2}\widehat{\bf \Phi}_{1}^\dagger \|_2}{\sqrt{1 + \gamma^{4q+4}\|\widehat{\bf \Phi}_{2}\widehat{\bf \Phi}_{1}^\dagger \|_2^2}},\\
& \text{sin}\angle(\mathcal{R}({\bf P}), \mathcal{R}({\bf v}_i)) \le \dfrac{\delta_i^{2q+1}\|\widehat{\bf \Phi}_{2}\widehat{\bf \Phi}_{1}^\dagger \|_2}{\sqrt{1 + \gamma^{4q+2}\|\widehat{\bf \Phi}_{2}\widehat{\bf \Phi}_{1}^\dagger \|_2^2}}.
\end{aligned} 
\end{equation}
\end{theorem}

\textit{Proof.} The proof is presented in Appendix \ref{Proo_ThCanASIn}.

Theorem \ref{ThCanASIn} states that $\theta_i \le \angle(\mathcal{R}({\bf Q}), \mathcal{R}({\bf u}_i))$ and $\phi_i \le \angle(\mathcal{R}({\bf P}), \mathcal{R}({\bf v}_i))$. This is because ${\bf u}_i \subseteq {\bf U}_k$ and ${\bf v}_i \subseteq{\bf V}_k$.

\subsection{Bounds for Low-Rank Approximation Errors}

When orthogonal projections associated with left or right singular vectors are utilized to form low-rank approximations, the error incurred remains the same. This is due to optimality of the SVD. In contrast, the quality of approximation by RU-QLP differs, as the error incurred depends upon which bases of $\bf Q$ and $\bf P$ are used. The theorem that follows demonstrates this in 2- and Frobenius norm.

\begin{theorem}\label{ThLRAEB}
Let ${\bf A} \in \mathbb R^{m \times n}$ be a rank-$k$ matrix, with $m\ge n$ and an SVD as in \eqref{eqSVD}. Let $\bf A$ have a RU-QLP decomposition as expounded in Section \ref{secRUQLP}, and $[{\bf Q}^T{\bf A}]_k$ and $[{\bf A}{\bf P}]_k$ be the rank-$k$ approximations provided by the SVD. We then have
\begin{equation}\notag
\begin{aligned}
 \|{\bf A}-{\bf Q}{\bf Q}^T{\bf A}\|_{2,F} & \le \|{\bf A} - {\bf Q}[{\bf Q}^T{\bf A}]_k\|_{2,F} \\
& \le  \Bigg(1 + \dfrac{\delta_k^{2q+1}\|\widehat{\bf \Phi}_{2}\widehat{\bf \Phi}_{1}^\dagger\|_2}{1 + \gamma^{4q+4}\|\widehat{\bf \Phi}_{2}\widehat{\bf \Phi}_{1}^\dagger \|_2^2}\Bigg)\|{\bf \Sigma}_\perp\|_{2,F}.\\
\end{aligned}
\end{equation}

\begin{equation}\notag
\begin{aligned}
\|{\bf A}-{\bf A}{\bf P}{\bf P}^T\|_{2,F} &\le \|{\bf A}-[{\bf A}{\bf P}]_k{\bf P}^T\|_{2,F} \\
& \le \Bigg(1 + \dfrac{\delta_k^{2q}\|\widehat{\bf \Phi}_{2}\widehat{\bf \Phi}_{1}^\dagger\|_2}{1 + \gamma^{4q+2}\|\widehat{\bf \Phi}_{2}\widehat{\bf \Phi}_{1}^\dagger \|_2^2}\Bigg)\|{\bf \Sigma}_\perp\|_{2,F}.
\end{aligned}
\end{equation}
\end{theorem}

\textit{Proof.} The proof is presented in Appendix \ref{Pro_ThLRAEB}.

The first error bound of Theorem \ref{ThLRAEB} is tighter than the second. This is expected and the interpretation is similar to that of principal angles $\theta_i$ and $\phi_i$: since the computation of $\bf Q$ requires $\bar{\bf P}$, one more step of iteration is therefore needed to form $\bf Q$, which in turn enhances its approximation quality. 

\subsection{Probabilistic Bounds}

This section provides the average case analysis for RU-QLP. The random matrix ${\bf \Phi}$ is standard Gaussian. $\widehat{\bf \Phi}_1$ and $\widehat{\bf \Phi}_2$ have the standard normal distribution and are statistically independent, due to rotational invariance.  

\begin{theorem}
Under the hypotheses of Theorem \ref{ThSigmakR22}, we have 
\begin{equation} \label{eqThSigR11Prob}
\mathbb{E}(\sigma_i({\bf R}_{11})) \ge \frac{\sigma_i}{\sqrt{1 + \delta_i^{4q+4}\omega^2}}, \quad i=1,...k,
\end{equation}
\begin{equation}  \label{eqThSigR22Prob}
\mathbb{E} \|{\bf R}_{22}\|_{2,F} \le (1 + C\delta_k^{2q+1})\|{\bf \Sigma}_\perp\|_{2,F},
\end{equation}
where $\omega = \omega_1\omega_2$, with $\omega_1 = \sqrt{m-k} + \sqrt{k+p}+7$, $\omega_2 = \frac{4\text{e}\sqrt{k+p}}{p+1}$, and $C= \sqrt{\dfrac{k}{p-1}}+{\dfrac{e\sqrt{(m-k)(p+k)}}{p}}$.
\label{ThSigmR22Probab}
\end{theorem}   

\textit{Proof.} The proof is given in Appendix \ref{Proo_ThSigmR22Probab}.

\begin{theorem}\label{ThSigmCanAngiProbab}
Under the hypotheses of Theorem \ref{ThCanABk}, we have 
\begin{equation}\notag
\begin{aligned}
& \mathbb{E}\text{sin}\theta_i \le \frac{\delta_i^{2q+2}\omega}{\sqrt{1 + \delta_i^{4q+4}\omega^{2}}}, \\
&
\mathbb{E}\text{sin}\phi_i \le \frac{\delta_i^{2q+1}\omega}{\sqrt{1 + \delta_i^{4q+2}\omega^{2}}},
\label{}
\end{aligned}
\end{equation}
where $\omega$ is defined in Theorem \ref{ThSigmR22Probab}.
\end{theorem}

\textit{Proof.} The proof is presented in Appendix \ref{Proo_ThSigmCanAngiProbab}.

\begin{theorem}\label{ThSigmCanAngFNProbab}
Under the hypotheses of Theorem \ref{ThCanABk}, we have 
\begin{equation}\notag
\begin{aligned}
&\mathbb{E} \|\text{sin}\angle(\mathcal{R}({\bf Q}), \mathcal{R}({\bf U}_k))\|_{2,F} \le \frac{\delta_k^{2q+1}\omega \|{\bf \Sigma}_\perp \|_{2,F}}{\sigma_k\sqrt{1 + \delta_k^{4q+4}\omega^{2}}},\\
& \mathbb{E}\|\text{sin}\angle(\mathcal{R}({\bf P}), \mathcal{R}({\bf V}_k))\|_{2,F} \le \frac{\delta_k^{2q}\omega \|{\bf \Sigma}_\perp \|_{2,F}}{\sigma_k\sqrt{1 + \delta_k^{4q+2}\omega^{2}}},
\label{}
\end{aligned}
\end{equation}
where $\omega$ is defined in Theorem \ref{ThSigmR22Probab}.
\end{theorem}

\begin{theorem}\label{ThSigmCanAngIndProbab}
Under the hypotheses of Theorem \ref{ThCanABk}, we have 
\begin{equation}\notag
\begin{aligned}
&\mathbb{E}\text{sin}\angle(\mathcal{R}({\bf Q}), \mathcal{R}({\bf u}_i))  \le \frac{\delta_i^{2q+2}\omega}{\sqrt{1 + \gamma^{4q+4}\omega^{2}}},\\
& \mathbb{E}\text{sin}\angle(\mathcal{R}({\bf P}), \mathcal{R}({\bf v}_i)) \le \frac{\delta_i^{2q+1}\omega}{\sqrt{1 + \gamma^{4q+2}\omega^{2}}},
\label{}
\end{aligned}
\end{equation}
where $\omega$ is defined in Theorem \ref{ThSigmR22Probab}.
\end{theorem}

Theorems \ref{ThSigmCanAngFNProbab} and \ref{ThSigmCanAngIndProbab} follow similarly as Theorem \ref{ThSigmCanAngiProbab}. 

\begin{theorem}
Under the hypotheses of Theorem \ref{ThLRAEB}, we have 
\begin{equation} \label{}
\mathbb{E} \|{\bf A}-{\bf Q}{\bf Q}^T{\bf A}\|_{2,F} \le (1 + C\delta_k^{2q+1})\|{\bf \Sigma}_\perp\|_{2,F},
\end{equation}
\begin{equation}  \label{}
\mathbb{E} \|{\bf A}-{\bf A}{\bf P}{\bf P}^T\|_{2,F} \le (1 + C\delta_k^{2q})\|{\bf \Sigma}_\perp\|_{2,F},
\end{equation}
where $C$ is defined in Theorem \ref{ThSigmR22Probab}.
\label{ThSigmLRAEProbab}
\end{theorem} 

The proof of this theorem is similar to that of \eqref{eqThSigR22Prob}.


\section{Numerical Simulations}
\label{secNuSim}

\textcolor{black}{We first investigate the performance of RU-QLP and several existing methods in terms of runtime for the decomposition of dense as well as sparse matrices.} Further, we empirically evaluate the effectiveness of the derived bounds; we utitlize several matrices with different properties. 
The simulations were performed in Python on an Intel Xeon Gold 6240 18-core processor @ 2.6 GHz with 251 GB of memory, equipped with a NVIDIA GeForce RTX 2080Ti GPU. 

We should add that a comprehensive comparison of RU-QLP with the state-of-the-art methods in approximation quality has been presented in \cite[Section V]{MFKC21}, demonstrating the high accuracy of RU-QLP on various classes of matrices. Hence we shall not conduct a similar study here.

\subsection{Comparison of Runtime}
\label{subRuntime}
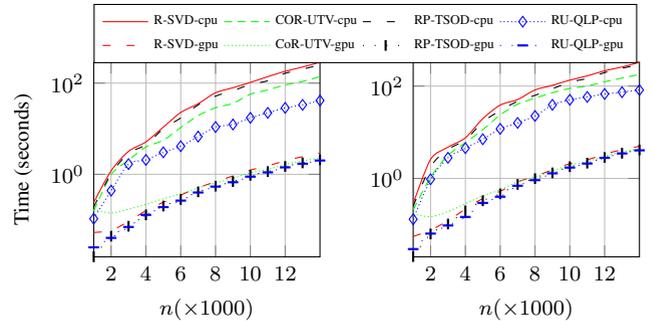
\begin{figure}[t]
\begin{center}       
%
%
%
\usetikzlibrary{positioning,calc}

\definecolor{mycolor1}{rgb}{0.00000,1.00000,1.00000}%
\definecolor{mycolor2}{rgb}{1.00000,0.00000,1.00000}%

\pgfplotsset{every axis label/.append style={font=\footnotesize},
every tick label/.append style={font=\footnotesize}
}

\begin{tikzpicture}[font=\footnotesize] 

\begin{axis}[%
name=ber,
ymode=log,
width  = 0.35\columnwidth,
height = 0.3\columnwidth,
scale only axis,
xmin  = 1,
xmax  = 14,
xlabel= {$n(\times 1000)$},
xmajorgrids,
ymin = 0.0155,
ymax = 286,
xtick       ={2, 4, 6, 8, 10, 12},
xticklabels ={$2$, $4$, $6$, $8$, $10$,$12$},
ylabel={Time (seconds)},
ymajorgrids,
]
\addplot+[smooth,color=red,solid, every mark/.append style={solid}, mark=none]
table[row sep=crcr]{
1	0.241104060000000  \\
2	1.19343245000000  \\
3	3.19015246000000  \\
4	5.13523473000000  \\
5	10.8544473600000  \\
6	22.7281131700000  \\
7	35.3174433100000  \\
8	61.4958522900000  \\
9	79.1443679300000  \\
10	103.389160510000  \\
11	138.491481420000  \\
12	181.375248790000  \\
13	225.257706050000  \\
14	285.692012430000  \\
};

\addplot+[smooth,color=green,densely dashed, every mark/.append style={thick}, mark=none]
table[row sep=crcr]{
1	0.154801010000000 \\
2	0.98253977000000 \\
3	2.03866029000000 \\
4	3.97045195000000 \\
5	5.58655310000000 \\
6	10.5743619200000 \\
7	19.3980237800000 \\
8	28.8666032600000 \\
9	35.9998547400000 \\
10	57.0625731300000 \\
11	70.7774292200000 \\
12	91.1938405000000 \\
13	109.012993750000 \\
14	138.701210920000 \\
};

\addplot+[smooth,color=black, loosely dashed, every mark/.append style={solid}, mark=none]
table[row sep=crcr]{
1	0.191713190000000  \\
2	1.02815589000000 \\
3	2.89121333000000 \\
4	4.44198756000000 \\
5	9.97143430000000 \\
6	17.2118681700000 \\
7	33.2912898700000 \\
8	51.8000698100000 \\
9	67.9649719600000 \\
10	93.9631518700000 \\
11	112.407884840000 \\
12	164.170593140000 \\
13	201.631197990000 \\
14	251.346982600000 \\
};

\addplot+[smooth,color=blue,densely dotted, every mark/.append style={solid}, mark=diamond]
table[row sep=crcr]{
1	0.106183830000000  \\
2	0.438319190000000 \\
3	1.67559960000000 \\
4	2.07013599000000 \\
5	3.03595620000000 \\
6	4.13540058000000 \\
7	6.71740681000000 \\
8	10.8919622300000 \\
9	12.3534156100000 \\
10	17.3295545600000 \\
11	22.2532701500000 \\
12	28.6345004400000 \\
13	33.3280490000000 \\
14	41.4701421900000 \\
};

\addplot+[smooth,color=red, loosely dashed, every mark/.append style={solid}, mark=none]
table[row sep=crcr]{
1	0.0532194400000000  \\
2	0.0583264800000000  \\
3	0.0927578900000000  \\
4	0.159852030000000  \\
5	0.243304130000000  \\
6	0.350067930000000  \\
7	0.509459420000000  \\
8	0.683178020000000  \\
9	0.923939560000000  \\
10	1.24045026000000  \\
11	1.52970505000000  \\
12	1.98626056000000  \\
13	2.43624265000000  \\
14	2.91972899000000  \\
};

\addplot+[smooth,color=green,densely dotted, every mark/.append style={solid}, mark=none]
table[row sep=crcr]{
1	0.171494030000000   \\
2	0.141719560000000  \\
3	0.172417780000000  \\
4	0.218248460000000  \\
5	0.292053940000000  \\
6	0.372118640000000  \\
7	0.502277830000000  \\
8	0.646128460000000  \\
9	0.822683500000000  \\
10	1.08362951000000  \\
11	1.28291085000000  \\
12	1.64693263000000  \\
13	1.95493033000000  \\
14	2.28909833000000  \\
};

\addplot+[smooth,color=black,loosely dotted, every mark/.append style={thick}, mark=|]
table[row sep=crcr]{
1	0.0155359500000000 \\
2	0.0439877300000000 \\
3	0.0725624300000000 \\
4	0.129248480000000 \\
5	0.201202700000000 \\
6	0.279899120000000 \\
7	0.412735200000000 \\
8	0.549839850000000 \\
9	0.692776300000000 \\
10	0.949797510000000 \\
11	1.14921699000000 \\
12	1.48631454000000 \\
13	1.78220520000000 \\
14	2.07973700000000 \\
};

\addplot+[smooth,color=blue,loosely dotted, every mark/.append style={thick}, mark=-]
table[row sep=crcr]{
1	0.0250478300000000  \\
2	0.0403722000000000  \\
3	0.0700651800000000  \\
4	0.128878640000000  \\
5	0.190710540000000  \\
6	0.264742010000000  \\
7	0.398910500000000  \\
8	0.537221980000000  \\
9	0.670624400000000  \\
10	0.887589620000000  \\
11	1.10407851000000  \\
12	1.42208323000000  \\
13	1.70775983000000  \\
14	2.00725150000000  \\
};

\end{axis}

\begin{axis}[%
name=SumRate,
at={($(ber.east)+(35,0em)$)},
		anchor= west,
ymode=log,
width  = 0.35\columnwidth,
height = 0.3\columnwidth,
scale only axis,
xmin   = 1,
xmax  = 14,
xlabel= {$n(\times 1000)$},
xmajorgrids,
ymin = 0.02,
ymax = 328,
xtick       ={2, 4, 6, 8, 10, 12},
xticklabels ={$2$, $4$, $6$, $8$, $10$,$12$},
ymajorgrids,
legend entries={R-SVD-cpu, COR-UTV-cpu, RP-TSOD-cpu, RU-QLP-cpu, R-SVD-gpu, CoR-UTV-gpu, RP-TSOD-gpu, RU-QLP-gpu}, 
legend style={at={(.99,1.3)},anchor=north east,draw=black,fill=white,legend cell align=left,font=\tiny, legend columns=4}
]

\addplot+[smooth,color=red,solid, every mark/.append style={solid}, mark=none]
table[row sep=crcr]{
1	0.283005910000000  \\
2	2.58282428000000 \\
3	4.75367098000000 \\
4	7.59234467000000 \\
5	19.7461503500000 \\
6	38.8318596800000 \\
7	51.1006019600000 \\
8	82.5368099200000 \\
9	103.869871850000 \\
10	134.775072720000 \\
11	166.535369590000 \\
12	215.337702470000 \\
13	259.154318570000 \\
14	327.812174940000 \\
};

\addplot+[smooth,color=green,densely dashed, every mark/.append style={solid}, mark=none]
table[row sep=crcr]{
1	0.207590630000000   \\
2	1.06037726000000  \\
3	3.16096668000000  \\
4	5.95985527000000  \\
5	11.3618715300000  \\
6	23.7403868700000  \\
7	40.5856289400000  \\
8	54.9580181600000  \\
9	71.7551997200000  \\
10	89.0564911800000  \\
11	100.361814120000  \\
12	124.392153790000  \\
13	149.205069110000  \\
14	178.397944120000  \\
};

\addplot+[smooth,color=black, loosely dashed, every mark/.append style={solid}, mark=none]
table[row sep=crcr]{
1	0.247594880000000  \\
2	1.80206146000000 \\
3	4.26774125000000 \\
4	6.80232677000000 \\
5	16.7154269200000 \\
6	29.5834717300000 \\
7	46.7346077900000 \\
8	63.0783029100000 \\
9	87.5729585600000 \\
10	121.163002200000 \\
11	147.122739740000 \\
12	196.633338360000 \\
13	243.521169470000 \\
14	291.003749850000 \\
};

\addplot+[smooth,color=blue,densely dotted, every mark/.append style={solid}, mark=diamond]
table[row sep=crcr]{
1	0.130138400000000  \\
2	0.957139780000000 \\
3	2.79854913000000 \\
4	4.52460041000000 \\
5	7.05830340000000 \\
6	11.8454018100000 \\
7	15.9411268700000 \\
8	22.5747904800000 \\
9	39.3493681400000 \\
10	50.5108012700000 \\
11	58.2002876800000 \\
12	67.4825069000000 \\
13	73.7502210100000 \\
14	82.2588456200000 \\
};

\addplot+[smooth,color=red, loosely dashed, every mark/.append style={solid}, mark=none]
table[row sep=crcr]{
1	0.0553027900000000   \\
2	0.0759458500000000  \\
3	0.121469860000000  \\
4	0.217722630000000  \\
5	0.366962360000000  \\
6	0.524031470000000  \\
7	0.838124040000000  \\
8	1.11087472000000  \\
9	1.53833385000000  \\
10	2.09558711000000  \\
11	2.52631495000000  \\
12	3.35166631000000  \\
13	4.19252875000000  \\
14	4.95670307000000  \\
};

\addplot+[smooth,color=green,densely dotted, every mark/.append style={thick}, mark=none]
table[row sep=crcr]{
1	0.169237850000000  \\
2	0.148000980000000 \\
3	0.191325070000000 \\
4	0.283211160000000 \\
5	0.420967200000000 \\
6	0.564162920000000 \\
7	0.841571780000000 \\
8	1.07525947000000 \\
9	1.44839747000000 \\
10	1.95292001000000 \\
11	2.30385056000000 \\
12	3.03417115000000 \\
13	3.73901036000000 \\
14	4.38754339000000 \\
};

\addplot+[smooth,color=black,loosely dotted, every mark/.append style={thick}, mark=|]
table[row sep=crcr]{
1	0.0202649400000000  \\ 
2	0.0626450800000000 \\
3	0.102748970000000 \\
4	0.197336150000000 \\
5	0.328228240000000 \\
6	0.462853190000000 \\
7	0.740112020000000 \\
8	0.987681910000000 \\
9	1.31770580000000 \\
10	1.79554017000000 \\
11	2.14236648000000 \\
12	2.84132388000000 \\
13	3.51615171000000 \\
14	4.14430850000000 \\
};

\addplot+[smooth,color=blue,loosely dotted, every mark/.append style={thick}, mark=-]
table[row sep=crcr]{
1	0.0294515800000000   \\
2	0.0637909200000000  \\
3	0.0955196000000000  \\
4	0.146150040000000  \\
5	0.294494180000000  \\
6	0.400873060000000  \\
7	0.692478360000000  \\
8	0.915888370000000  \\
9	1.28740785000000  \\
10	1.71691501000000  \\
11	2.09925728000000  \\
12	2.78892622000000  \\
13	3.47490020000000  \\
14	4.05788136000000  \\
};

\end{axis}

\end{tikzpicture}%
\captionsetup{justification=centering,font=scriptsize}  
\caption{Computational times for decomposition of dense matrices with $d=0.04n$. Left: $q=0$. Right: $q=2$.} 
\label{figSzzt}      
\end{center}
\end{figure}

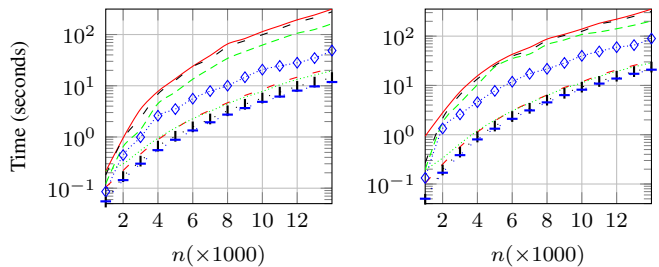
\begin{figure}[t]
\begin{center}       
%
%
%
\usetikzlibrary{positioning,calc}

\definecolor{mycolor1}{rgb}{0.00000,1.00000,1.00000}%
\definecolor{mycolor2}{rgb}{1.00000,0.00000,1.00000}%

\pgfplotsset{every axis label/.append style={font=\footnotesize},
every tick label/.append style={font=\footnotesize}
}

\begin{tikzpicture}[font=\footnotesize] 

\begin{axis}[%
name=ber,
ymode=log,
width  = 0.35\columnwidth,
height = 0.3\columnwidth,
scale only axis,
xmin  = 1,
xmax  = 14,
xlabel= {$n(\times 1000)$},
xmajorgrids,
ymin = 0.05,
ymax = 314,
xtick       ={2, 4, 6, 8, 10, 12},
xticklabels ={$2$, $4$, $6$, $8$, $10$,$12$},
ylabel={Time (seconds)},
ymajorgrids,
]
\addplot+[smooth,color=red,solid, every mark/.append style={solid}, mark=none]
table[row sep=crcr]{
1	0.212578830000000\\
2	0.958144960000000\\
3	3.49383181000000\\
4	7.74451935000000\\
5	13.7078262000000\\
6	24.0353952100000\\
7	38.8447393800000\\
8	65.5317179600000\\
9	83.5627927800000\\
10	114.003251670000\\
11	147.870708530000\\
12	193.786552610000\\
13	241.032141920000\\
14	313.642831980000  \\
};

\addplot+[smooth,color=green,densely dashed, every mark/.append style={solid}, mark=none]
table[row sep=crcr]{
1	0.122290670000000\\
2	0.64668649000000\\
3	1.40082687000000\\
4	4.48099726000000\\
5	7.51484877000000\\
6	13.4443974500000\\
7	21.5883550600000\\
8	34.2596207900000\\
9	46.5455942200000\\
10	62.1926757700000\\
11	81.7285962700000\\
12	105.856664120000\\
13	125.433694240000\\
14	162.919610920000 \\
};

\addplot+[smooth,color=black, loosely dashed, every mark/.append style={solid}, mark=none]
table[row sep=crcr]{
1	0.180877750000000 \\
2	0.931661310000000 \\
3	2.38377213000000 \\
4	6.90664476000000 \\
5	12.7158904100000 \\
6	23.4883832900000 \\
7	34.0712072300000 \\
8	57.6119287000000 \\
9	75.7209362400000 \\
10	99.0501080200000 \\
11	130.114932120000 \\
12	172.945967850000 \\
13	213.810632530000 \\
14	277.247590360000 \\
};

\addplot+[smooth,color=blue,densely dotted, every mark/.append style={solid}, mark=diamond]
table[row sep=crcr]{
1	0.0856093800000000 \\
2	0.442307890000000 \\
3	0.991510030000000 \\
4	2.59631413000000 \\
5	3.50279421000000 \\
6	5.59306484000000 \\
7	7.87030536000000 \\
8	10.0259648600000 \\
9	14.6768751100000 \\
10	20.6401033400000 \\
11	24.3819590200000 \\
12	28.1212095000000 \\
13	34.6876018000000 \\
14	48.6453001500000 \\
};

\addplot+[smooth,color=red, loosely dashed, every mark/.append style={solid}, mark=none]
table[row sep=crcr]{
1	0.104064080000000 \\
2	0.221635990000000 \\
3	0.464123630000000 \\
4	0.886364340000000 \\
5	1.46728306000000 \\
6	2.31496146000000 \\
7	3.32324336000000 \\
8	4.65022068000000 \\
9	6.47465920000000 \\
10	8.45082388000000 \\
11	11.4752701000000 \\
12	14.7001327500000 \\
13	18.0763857800000 \\
14	22.0970159500000  \\
};

\addplot+[smooth,color=green, densely dotted, every mark/.append style={thick}, mark=none]
table[row sep=crcr]{
1	0.219730090000000  \\
2	0.298905990000000  \\
3	0.54559703000000  \\
4	0.948496700000000  \\
5	1.50309842000000  \\
6	2.26904602000000  \\
7	3.16970389000000  \\
8	4.44667706000000  \\
9	6.01395769000000  \\
10	7.71594660000000  \\
11	10.4029058200000  \\
12	13.2971700400000  \\
13	16.1539357200000  \\
14	19.6118294000000  \\
};

\addplot+[smooth,color=black,loosely dotted, every mark/.append style={thick}, mark=|]
table[row sep=crcr]{
1	0.0536289900000000\\
2	0.168540360000000\\
3	0.337773300000000\\
4	0.663215540000000\\
5	1.03255782000000\\
6	1.57389002000000\\
7	2.24665453000000\\
8	3.24442992000000\\
9	4.35014646000000\\
10	5.75135450000000\\
11	7.50404899000000\\
12	9.72053285000000\\
13	11.7975395900000\\
14	14.4381717200000 \\
};

\addplot+[smooth,color=blue,loosely dotted, every mark/.append style={thick}, mark=-]
table[row sep=crcr]{
1	0.0553650900000000  \\
2	0.143886180000000  \\
3	0.301157500000000  \\
4	0.553499910000000  \\
5	0.882598570000000  \\
6	1.34130285000000  \\
7	1.93731356000000  \\
8	2.73465562000000  \\
9	3.66662006000000  \\
10	4.85425537000000  \\
11	6.15036509000000  \\
12	8.03906922000000  \\
13	9.70926111000000  \\
14	11.8265928000000  \\
};

\end{axis}

\begin{axis}[%
name=SumRate,
at={($(ber.east)+(35,0em)$)},
		anchor= west,
ymode=log,
width  = 0.35\columnwidth,
height = 0.3\columnwidth,
scale only axis,
xmin   = 1,
xmax  = 14,
xlabel= {$n(\times 1000)$},
xmajorgrids,
ymin = 0.04,
ymax = 354,
xtick       ={2, 4, 6, 8, 10, 12},
xticklabels ={$2$, $4$, $6$, $8$, $10$,$12$},
ymajorgrids,
]

\addplot+[smooth,color=red,solid, every mark/.append style={solid}, mark=none]
table[row sep=crcr]{
1	0.936912070000000 \\
2	2.78901302000000 \\
3	7.43788999000000 \\
4	15.6009934500000 \\
5	27.7189349500000 \\
6	42.7086105300000 \\
7	58.9808017000000 \\
8	87.1567489500000 \\ 
9	109.037361860000 \\
10	139.072740610000 \\
11	182.315093220000 \\
12	221.054141880000 \\
13	274.824175720000 \\
14	353.241991580000 \\
};

\addplot+[smooth,color=green, densely dashed, every mark/.append style={solid}, mark=none]
table[row sep=crcr]{
1	0.196017070000000  \\
2	1.76025341000000  \\
3	4.41899329000000  \\
4	10.3306545600000  \\
5	23.6342158300000  \\
6	35.1373862600000  \\
7	43.3229066700000  \\
8	67.4051970800000  \\
9	83.8496841200000  \\
10	108.092021700000  \\
11	121.326035560000  \\
12	141.218814670000  \\
13	166.735466000000  \\
14	203.090985600000  \\
};

\addplot+[smooth,color=black, loosely dashed, every mark/.append style={solid}, mark=none]
table[row sep=crcr]{
1	0.266327680000000  \\
2	2.14220406000000 \\
3	6.47391051000000 \\
4	13.5434455300000 \\
5	25.2817319000000 \\
6	39.1933133000000 \\
7	51.5138092000000 \\
8	81.3237144400000 \\
9	103.305128570000 \\ 
10	125.273586390000 \\
11	174.307920160000 \\
12	204.807871700000 \\
13	252.847365200000 \\
14	317.440730810000 \\
};

\addplot+[smooth,color=blue,densely dotted, every mark/.append style={solid}, mark=diamond]
table[row sep=crcr]{
1	0.132755130000000 \\
2	1.33716305000000\\
3	2.59446424000000\\
4	4.57630497000000\\
5	7.74871349000000\\
6	11.8869378600000\\
7	17.2257929400000\\
8	21.7343192700000\\
9	28.5887792700000\\
10	39.8569531400000\\
11	48.8213100400000\\
12	59.7063531300000\\
13	65.7487861500000\\
14	89.3493040200000 \\
};

\addplot+[smooth,color=red, loosely dashed, every mark/.append style={solid}, mark=none]
table[row sep=crcr]{
1	0.103448680000000  \\
2	0.252348730000000 \\
3	0.567618300000000 \\
4	1.13570659000000 \\
5	1.91444046000000 \\
6	3.05816934000000 \\
7	4.48179069000000 \\
8	6.38352895000000 \\
9	9.11101682000000 \\
10	11.7939940500000 \\
11	16.1930208700000 \\
12	20.5365016000000 \\
13	25.5360550200000 \\
14	30.9901761100000  \\
};

\addplot+[smooth,color=green, densely dotted, every mark/.append style={thick}, mark=none]
table[row sep=crcr]{
1	0.195008210000000 \\
2	0.334162450000000\\
3	0.656891420000000\\
4	1.20116484000000\\
5	1.96405425000000\\
6	3.03625207000000\\
7	4.38395052000000\\
8	6.23734031000000\\
9	8.73616550000000\\
10	11.2204221000000\\
11	15.2678508500000\\
12	19.2393909500000\\
13	23.8344912500000\\
14	28.8472734700000 \\
};

\addplot+[smooth,color=black,loosely dotted, every mark/.append style={thick}, mark=|]
table[row sep=crcr]{
1	0.0490143500000000  \\
2	0.192594430000000 \\
3	0.443645950000000 \\
4	0.898607520000000 \\
5	1.50036426000000 \\
6	2.33622451000000 \\
7	3.43635335000000 \\
8	5.02268364000000 \\
9	6.92697256000000 \\
10	9.22800896000000 \\
11	12.0172837300000 \\
12	15.3933273100000 \\
13	19.0791935900000 \\
14	23.5410854600000 \\
};

\addplot+[smooth,color=blue,loosely dotted, every mark/.append style={thick}, mark=-]
table[row sep=crcr]{
1	0.0502969500000000   \\
2	0.169319530000000  \\
3	0.385421630000000  \\
4	0.803157260000000  \\
5	1.32481883000000  \\
6	2.10176601000000  \\
7	3.11351035000000  \\
8	4.53065314000000  \\
9	6.36035843000000  \\
10	8.19150250000000  \\
11	10.8936596400000  \\
12	13.7898546500000  \\
13	17.2262982600000  \\
14	20.8581206100000  \\
};

\end{axis}

\end{tikzpicture}%
\captionsetup{justification=centering,font=scriptsize}  
\caption{Computational times for decomposition of dense matrices with $d=0.2n$. Left: $q=0$. Right: $q=2$.} 
\label{figSzo}      
\end{center}
\end{figure}

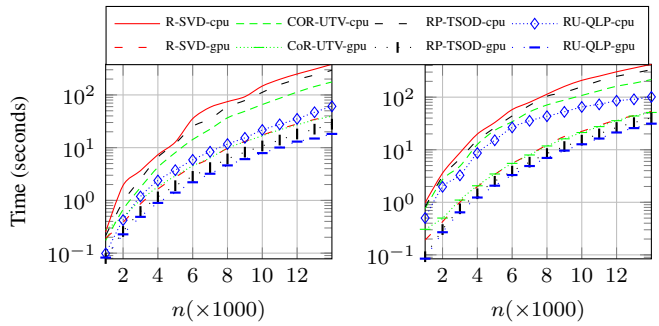
\begin{figure}[t]
\begin{center}       
%
%
%
\usetikzlibrary{positioning,calc}

\definecolor{mycolor1}{rgb}{0.00000,1.00000,1.00000}%
\definecolor{mycolor2}{rgb}{1.00000,0.00000,1.00000}%

\pgfplotsset{every axis label/.append style={font=\footnotesize},
every tick label/.append style={font=\footnotesize}
}

\begin{tikzpicture}[font=\footnotesize] 

\begin{axis}[%
name=ber,
ymode=log,
width  = 0.35\columnwidth,
height = 0.3\columnwidth,
scale only axis,
xmin  = 1,
xmax  = 14,
xlabel= {$n(\times 1000)$},
xmajorgrids,
ymin = 0.078,
ymax = 381,
xtick       ={2, 4, 6, 8, 10, 12},
xticklabels ={$2$, $4$, $6$, $8$, $10$,$12$},
ylabel={Time (seconds)},
ymajorgrids,
]
\addplot+[smooth,color=red,solid, every mark/.append style={solid}, mark=none]
table[row sep=crcr]{
1	0.252545710000000   \\
2	1.92591965000000  \\
3	3.65816343000000  \\
4	7.35687780000000  \\
5	12.7484203000000  \\
6	35.5050234800000  \\
7	58.4022012900000  \\
8	74.6871344400000  \\
9	92.5897178100000  \\
10	147.361810920000  \\
11	197.692769590000  \\
12	249.698409320000  \\
13	306.767103250000  \\
14	380.736653270000  \\
};

\addplot+[smooth,color=green,densely dashed, every mark/.append style={solid}, mark=none]
table[row sep=crcr]{
1	0.174327730000000  \\
2	0.671492220000000 \\
3	1.75803506000000 \\
4	4.25249189000000 \\
5	8.05291623000000 \\
6	14.2752064500000 \\
7	22.3520672300000 \\
8	36.8582477600000 \\
9	48.9711790100000 \\
10	65.7526674900000 \\
11	86.3508602400000 \\
12	113.152390720000 \\
13	143.992952050000 \\
14	176.036702630000 \\
};

\addplot+[smooth,color=black, loosely dashed, every mark/.append style={solid}, mark=none]
table[row sep=crcr]{
1	0.214387600000000  \\
2	0.942511680000000 \\
3	2.53328216000000 \\
4	6.60942256000000 \\
5	11.9720152600000 \\
6	25.2816160300000 \\ 
7	35.0302974600000 \\
8	59.4761960500000 \\
9	78.4303098300000 \\
10	111.971218230000 \\
11	156.836116550000 \\
12	196.266689480000 \\
13	238.321487430000 \\
14	293.101206660000 \\
};

\addplot+[smooth,color=blue,densely dotted, every mark/.append style={solid}, mark=diamond]
table[row sep=crcr]{
1	0.0986092100000000  \\
2	0.422449830000000 \\
3	1.18651700000000 \\
4	2.38593602000000 \\
5	3.74649149000000 \\
6	5.87524962000000 \\
7	8.35851097000000 \\
8	11.7123969200000 \\
9	15.3456697500000 \\
10	21.6686148000000 \\
11	27.6592721900000 \\
12	35.1572180400000 \\
13	46.9579979800000 \\
14	60.7423327000000 \\
};


\addplot+[smooth,color=red, loosely dashed, every mark/.append style={solid}, mark=none]
table[row sep=crcr]{
1	0.192280820000000  \\
2	0.384517620000000 \\
3	0.859930900000000 \\
4	1.64461927000000 \\
5	2.73433328000000 \\
6	4.35719161000000 \\
7	6.24108706000000 \\
8	9.73899026000000 \\
9	13.1265175800000 \\
10	17.3069884800000 \\
11	22.4547872100000 \\
12	28.1101202000000 \\
13	34.2245282200000 \\
14	40.6824640300000 \\
};

\addplot+[smooth,color=green,densely dotted, every mark/.append style={solid}, mark=none]
table[row sep=crcr]{
1	0.279640720000000   \\
2	0.489944460000000  \\
3	1.01246881000000  \\
4	1.80574689000000  \\
5	2.90223083000000  \\
6	4.51776195000000  \\
7	6.39602280000000  \\
8	9.80762711000000  \\
9	12.9057358300000  \\
10	16.8054046200000  \\
11	21.8084115000000  \\
12	27.0407146900000  \\
13	33.5762642400000  \\
14	39.3648717400000  \\
};

\addplot+[smooth,color=black,loosely dotted, every mark/.append style={thick}, mark=|]
table[row sep=crcr]{
1	0.0791002300000000  \\
2	0.257977390000000 \\
3	0.616081810000000 \\
4	1.16583810000000 \\
5	1.90090261000000 \\
6	2.85305800000000 \\
7	4.18880086000000 \\
8	6.11729484000000 \\
9	8.21759863000000 \\
10	10.7727110400000 \\
11	13.9435555500000 \\
12	18.1896062900000 \\
13	22.4321639500000 \\
14	27.8813355400000 \\
};

\addplot+[smooth,color=blue,loosely dotted, every mark/.append style={thick}, mark=-]
table[row sep=crcr]{
1	0.0823688500000000   \\
2	0.227393250000000  \\
3	0.487526130000000  \\
4	0.896560100000000  \\
5	1.40499730000000  \\
6	2.20685816000000  \\
7	3.20248318000000  \\
8	4.61178470000000  \\
9	6.06527996000000  \\
10	7.89050498000000  \\
11	10.0587375200000  \\
12	13.0287337800000  \\
13	14.9056294900000  \\
14	18.2330559700000  \\
};

\end{axis}

\begin{axis}[%
name=SumRate,
at={($(ber.east)+(35,0em)$)},
		anchor= west,
ymode=log,
width  = 0.35\columnwidth,
height = 0.3\columnwidth,
scale only axis,
xmin   = 1,
xmax  = 14,
xlabel= {$n(\times 1000)$},
xmajorgrids,
ymin = 0.084,
ymax = 422,
xtick       ={2, 4, 6, 8, 10, 12},
xticklabels ={$2$, $4$, $6$, $8$, $10$,$12$},
ymajorgrids,
legend entries={R-SVD-cpu, COR-UTV-cpu, RP-TSOD-cpu, RU-QLP-cpu, R-SVD-gpu, CoR-UTV-gpu, RP-TSOD-gpu, RU-QLP-gpu}, 
legend style={at={(0.99,1.3)},anchor=north east,draw=black,fill=white,legend cell align=left,font=\tiny, legend columns=4}
]

\addplot+[smooth,color=red,solid, every mark/.append style={solid}, mark=none]
table[row sep=crcr]{
1	0.982766320000000  \\
2	3.59141556000000 \\
3	8.75716013000000 \\
4	19.8793540000000 \\
5	32.9388737100000 \\
6	58.0617438000000 \\
7	79.8156583300000 \\
8	113.763719200000 \\
9	160.031530500000 \\ 
10	209.102366510000 \\
11	251.251093800000 \\
12	303.807775680000 \\
13	358.261550720000 \\
14	421.935185190000 \\
};

\addplot+[smooth,color=green,densely dashed, every mark/.append style={solid}, mark=none]
table[row sep=crcr]{
1	0.761842230000000  \\
2	2.70966870000000 \\
3	4.87184405000000 \\
4	12.6143130100000 \\
5	23.3713691200000 \\
6	34.0707124500000 \\
7	52.5963308800000 \\
8	71.0270363700000 \\
9	87.8583575500000 \\
10	107.857718290000 \\
11	132.173962120000 \\
12	161.626020610000 \\
13	184.176433920000 \\
14	216.554761110000  \\
};

\addplot+[smooth,color=black, loosely dashed, every mark/.append style={solid}, mark=none]
table[row sep=crcr]{
1	0.811817200000000 \\
2	3.10130177000000\\
3	6.75783521000000\\
4	15.1145870100000\\
5	26.2795419100000\\
6	43.9146975900000\\
7	69.7388044000000\\
8	104.208884890000\\
9	123.771015580000\\
10	157.894195440000\\
11	203.438465770000\\
12	249.473425630000\\
13	286.576174140000\\
14	333.7892638400000 \\
};

\addplot+[smooth,color=blue,densely dotted, every mark/.append style={solid}, mark=diamond]
table[row sep=crcr]{
1	0.502360190000000  \\
2	1.96781820000000 \\
3	3.26295447000000 \\
4	8.64994353000000 \\
5	15.1262269600000 \\
6	26.5744960900000 \\
7	35.5486564000000 \\
8	43.3870233900000 \\
9	52.4938296700000 \\
10	65.1122289300000 \\
11	72.9372338100000 \\
12	85.3779647900000 \\
13	91.0329215500000 \\
14	101.012458800000 \\
};

\addplot+[smooth,color=red, loosely dashed, every mark/.append style={solid}, mark=none]
table[row sep=crcr]{
1	0.191343930000000   \\
2	0.426389790000000  \\
3	1.00697699000000  \\
4	1.98511715000000  \\
5	3.38648195000000  \\
6	5.51664834000000  \\  
7	8.08434286000000  \\
8	12.3192606000000  \\
9	17.8679771400000  \\
10	22.3392354500000  \\
11	28.9710731000000  \\
12	37.0054572600000  \\
13	44.6639303700000  \\
14	53.8866620100000  \\
};

\addplot+[smooth,color=green,densely dotted, every mark/.append style={solid}, mark=-]
table[row sep=crcr]{
1	0.306318190000000  \\
2	0.501261620000000 \\
3	1.10272574000000 \\
4	2.06246333000000 \\
5	3.45424480000000 \\
6	5.50225167000000 \\
7	7.94650445000000 \\
8	11.7962917800000 \\
9	16.1249763500000 \\
10	21.2351052800000 \\
11	27.5247959100000 \\
12	35.6711874000000 \\
13	43.8453500800000 \\
14	51.2296555100000 \\
};

\addplot+[smooth,color=black,loosely dotted, every mark/.append style={thick}, mark=|]
table[row sep=crcr]{
1	0.0921213200000000  \\ 
2	0.302784630000000 \\
3	0.752704860000000 \\
4	1.45286765000000 \\
5	2.51084471000000 \\
6	3.92329121000000 \\
7	5.80143957000000 \\
8	8.44874983000000 \\
9	11.4845773700000 \\
10	15.2072114900000 \\
11	19.9957644500000 \\
12	26.4983895300000 \\
13	32.7638459700000 \\
14	40.6897659800000 \\
};

\addplot+[smooth,color=blue,loosely dotted, every mark/.append style={thick}, mark=-]
table[row sep=crcr]{
1	0.084651660000000 \\
2	0.269421670000000 \\
3	0.64398289000000 \\
4	1.23214092000000 \\
5	2.06879325000000 \\
6	3.30053544000000 \\
7	4.89524622000000 \\
8	6.97995214000000 \\
9	9.61917276000000 \\
10	12.7010971500000 \\
11	16.3194620100000 \\
12	21.2884031800000 \\
13	25.7718619300000 \\
14	31.3168162800000  \\
};

\end{axis}

\end{tikzpicture}%
\captionsetup{justification=centering,font=scriptsize}  
\caption{Computational times for decomposition of dense matrices with $d=0.3n$. Left: $q=0$. Right: $q=2$.} 
\label{figSzof}      
\end{center}
\end{figure}

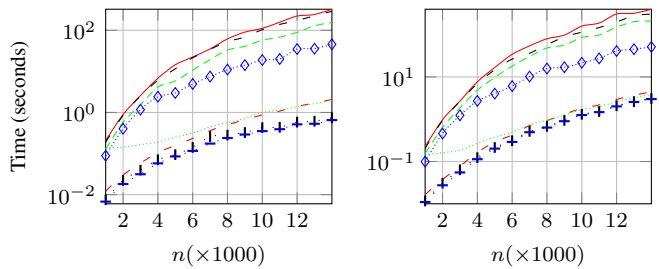
\begin{figure}[t]
\begin{center}       
%
%
%
\usetikzlibrary{positioning,calc}

\definecolor{mycolor1}{rgb}{0.00000,1.00000,1.00000}%
\definecolor{mycolor2}{rgb}{1.00000,0.00000,1.00000}%

\pgfplotsset{every axis label/.append style={font=\footnotesize},
every tick label/.append style={font=\footnotesize}
}

\begin{tikzpicture}[font=\footnotesize] 

\begin{axis}[%
name=ber,
ymode=log,
width  = 0.35\columnwidth,
height = 0.3\columnwidth,
scale only axis,
xmin  = 1,
xmax  = 14,
xlabel= {$n(\times 1000)$},
xmajorgrids,
ymin = 0.006,
ymax = 325,
xtick       ={2, 4, 6, 8, 10, 12},
xticklabels ={$2$, $4$, $6$, $8$, $10$,$12$},
ylabel={Time (seconds)},
ymajorgrids,
]
\addplot+[smooth,color=red,solid, every mark/.append style={solid}, mark=none]
table[row sep=crcr]{
1	0.208235343297323  \\
2	0.890387535095215 \\
3	2.58618259429932 \\
4	6.80056206385295 \\
5	14.4708412488302 \\
6	22.4127083619436 \\
7	35.0698771476746 \\
8	61.6804730892181 \\
9	90.3398737112681 \\ 
10	111.290565729141 \\
11	157.065581798553 \\
12	221.725306431452 \\
13	241.365882794062 \\
14	324.617106358210  \\
};

\addplot+[smooth,color=green,densely dashed, every mark/.append style={thick}, mark=none]
table[row sep=crcr]{
1	0.127088228861491 \\
2	0.578653891881307 \\
3	1.54356479644775 \\
4	4.14724373817444 \\
5	5.83091648419698 \\
6	10.6526358922323 \\
7	19.3609929084778 \\
8	33.1143843332926 \\ 
9	40.6807466348012 \\
10	57.0707100232442 \\
11	68.1817624568939 \\
12	88.5217675367991 \\
13	129.097711960475 \\
14	154.530131419500 \\
};

\addplot+[smooth,color=black, loosely dashed, every mark/.append style={solid}, mark=none]
table[row sep=crcr]{
1	0.188458363215129 \\
2	0.854038397471110 \\
3	2.52849260965983 \\
4	6.04660606384277 \\
5	11.5352809429169 \\
6	21.4113042354584 \\
7	34.4399495124817 \\
8	52.7580760320028 \\ 
9	68.4392337799072 \\
10	103.546522061030 \\
11	138.887496709824 \\
12	172.677163998286 \\
13	225.635247468948 \\
14	288.649590651194 \\
};

\addplot+[smooth,color=blue,densely dotted, every mark/.append style={solid}, mark=diamond]
table[row sep=crcr]{
1	0.0884248415629069 \\
2	0.401614030202230 \\
3	1.16963005065918 \\
4	2.41901723543803 \\
5	2.98030662536621 \\
6	4.94483645757039 \\
7	7.45220271746318 \\
8	11.0658488273621 \\
9	14.1630409666697 \\
10	18.8572448889414 \\
11	19.9861559867859 \\
12	34.9176481564840 \\
13	35.6027066707611 \\
14	45.5799497763316 \\ 
};

\addplot+[smooth,color=red, loosely dashed, every mark/.append style={solid}, mark=none]
table[row sep=crcr]{
1	0.0121626257896423   \\
2	0.0299673676490784  \\
3	0.0584239959716797  \\
4	0.105781614780426  \\
5	0.157014966011047  \\
6	0.233230054378510  \\
7	0.361199498176575  \\
8	0.503439486026764  \\
9	0.659185707569122  \\
10	0.875195860862732  \\
11	1.07869184017181  \\
12	1.39842611551285  \\
13	1.70798701047897  \\
14	2.06411194801331 \\
};

\addplot+[smooth,color=green,densely dotted, every mark/.append style={solid}, mark=none]
table[row sep=crcr]{
1	0.138586580753326  \\
2	0.142553985118866  \\
3	0.165307700634003  \\
4	0.189065933227539  \\
5	0.251366198062897  \\
6	0.313646972179413  \\
7	0.447465956211090  \\
8	0.562461853027344  \\
9	0.708173274993897  \\  
10	1.00174802541733  \\
11	1.13833558559418  \\
12	1.37094026803970  \\
13	1.64107972383499  \\
14	1.95342653989792  \\
};

\addplot+[smooth,color=black,loosely dotted, every mark/.append style={thick}, mark=|]
table[row sep=crcr]{
1	0.00751751661300659  \\
2	0.0221837759017944  \\
3	0.0394410490989685  \\
4	0.0707383751869202  \\
5	0.0943072438240051  \\
6	0.141171514987946  \\
7	0.259083719253540  \\
8	0.281415729522705  \\
9	0.324821710586548  \\
10	0.428015828132629  \\
11	0.436132013797760  \\
12	0.581888437271118  \\
13	0.612648844718933  \\
14	0.743159651756287 \\
};

\addplot+[smooth,color=blue,loosely dotted, every mark/.append style={thick}, mark=-]
table[row sep=crcr]{
1	0.00676435232162476  \\
2	0.0181699991226196  \\
3	0.0314189791679382  \\
4	0.0577749609947205  \\
5	0.0852021574974060  \\
6	0.115993738174438  \\
7	0.174886643886566  \\
8	0.238817632198334  \\
9	0.290941894054413  \\
10	0.351316988468170  \\
11	0.391824245452881  \\
12	0.513773202896118  \\
13	0.530351459980011  \\
14	0.652233839035034  \\
};

\end{axis}

\begin{axis}[%
name=SumRate,
at={($(ber.east)+(35,0em)$)},
		anchor= west,
ymode=log,
width  = 0.35\columnwidth,
height = 0.3\columnwidth,
scale only axis,
xmin   = 1,
xmax  = 14,
xlabel= {$n(\times 1000)$},
xmajorgrids,
ymin = 0.01,
ymax = 400,
xtick       ={2, 4, 6, 8, 10, 12},
xticklabels ={$2$, $4$, $6$, $8$, $10$,$12$},
ymajorgrids,
]

\addplot+[smooth,color=red,solid, every mark/.append style={solid}, mark=none]
table[row sep=crcr]{
1	0.224425713221232  \\
2	0.975392977396647 \\
3	2.69070577621460 \\
4	7.69545133908590 \\
5	16.7902293999990 \\
6	34.2808701992035 \\
7	54.3294597466787 \\
8	80.0724250475566 \\
9	97.6383267243703 \\
10	159.445385376612 \\
11	188.409300088882 \\
12	309.753629128138 \\
13	319.514836947123 \\
14	399.949247598648 \\
};

\addplot+[smooth,color=green,densely dashed, every mark/.append style={solid}, mark=none]
table[row sep=crcr]{
1	0.153662125269572 \\
2	0.644044796625773 \\
3	1.87810015678406 \\
4	4.85232694943746 \\
5	9.73664903640747 \\
6	18.1339797973633 \\ 
7	30.5792985757192 \\
8	48.4817812442780 \\
9	56.4554421106974 \\
10	86.1167678833008 \\
11	97.7045049667358 \\
12	150.784797588984 \\
13	202.913028637568 \\
14	210.462300300598 \\
};

\addplot+[smooth,color=black, loosely dashed, every mark/.append style={solid}, mark=none]
table[row sep=crcr]{
1	0.203855355580648 \\
2	0.939382632573446 \\
3	2.63807074228923 \\
4	7.00312709808350 \\
5	14.7636828422546 \\
6	26.4608079592387 \\
7	36.2562349637350 \\
8	57.8880894978841 \\
9	88.8277381261190 \\
10	114.142429033915 \\
11	147.578892230988 \\
12	199.229925155640 \\
13	282.716256697973 \\
14	306.264930725098 \\
};

\addplot+[smooth,color=blue,densely dotted, every mark/.append style={solid}, mark=diamond]
table[row sep=crcr]{
1	0.0991413593292236 \\
2	0.457987546920776 \\
3	1.23827107747396 \\
4	2.71337890625000 \\
5	4.03578774134318 \\
6	6.06899666786194 \\
7	10.2310938835144 \\
8	15.7998316287994 \\
9	16.9929495652517 \\
10	21.6503103574117 \\
11	27.4790658950806 \\
12	41.0532949765523 \\
13	45.4156056245168 \\
14	51.3189520835877 \\
};

\addplot+[smooth,color=red, loosely dashed, every mark/.append style={solid}, mark=none]
table[row sep=crcr]{
1	0.0165263414382935 \\
2	0.0376542806625366 \\
3	0.0772449374198914 \\
4	0.163626134395599 \\
5	0.291396856307983 \\
6	0.423526346683502 \\
7	0.707990646362305 \\
8	0.925164401531220 \\
9	1.31291455030441 \\
10	1.80865031480789 \\
11	2.17755341529846 \\
12	2.94207668304443 \\
13	3.69742816686630 \\
14	4.40182858705521  \\
};

\addplot+[smooth,color=green,densely dotted, every mark/.append style={thick}, mark=none]
table[row sep=crcr]{
1	0.145209491252899 \\
2	0.160396337509155 \\
3	0.186495065689087 \\
4	0.266667664051056 \\
5	0.377501428127289 \\
6	0.493679404258728 \\
7	0.741440713405609 \\
8	0.931863367557526 \\
9	1.25768387317657 \\
10	1.70675122737885 \\
11	1.97046220302582 \\
12	2.61434239149094 \\
13	3.24756431579590 \\
14	3.81688559055328 \\
};

\addplot+[smooth,color=black,loosely dotted, every mark/.append style={thick}, mark=|]
table[row sep=crcr]{
1	0.0116412639617920 \\
2	0.0306491255760193 \\
3	0.0637126564979553 \\
4	0.123219311237335 \\
5	0.214614748954773 \\
6	0.310050487518311 \\
7	0.515724062919617 \\
8	0.677040278911591 \\
9	0.935759782791138 \\
10	1.30775409936905 \\
11	1.51970732212067 \\
12	2.08196437358856 \\
13	2.61786413192749 \\
14	3.08890038728714 \\
};

\addplot+[smooth,color=blue,loosely dotted, every mark/.append style={thick}, mark=-]
table[row sep=crcr]{
1	0.0109807848930359 \\
2	0.0274227857589722 \\
3	0.0550214648246765 \\
4	0.113389968872070 \\
5	0.201908349990845 \\
6	0.290483474731445 \\ 
7	0.493421196937561 \\
8	0.635043799877167 \\
9	0.898339092731476 \\
10	1.25552260875702 \\ 
11	1.47170114517212 \\
12	2.01444232463837 \\
13	2.56072068214417 \\
14	2.99778193235397 \\
};

\end{axis}

\end{tikzpicture}%
\captionsetup{justification=centering,font=scriptsize}  
\caption{\textcolor{black}{Computational times for decomposition of sparse matrices with $d=0.04n$. Left: $q=0$. Right: $q=2$.}} 
\label{figSsp002}      
\end{center}
\end{figure}

\begin{figure}[t]
\begin{center}       
%
%
%
\usetikzlibrary{positioning,calc}

\definecolor{mycolor1}{rgb}{0.00000,1.00000,1.00000}%
\definecolor{mycolor2}{rgb}{1.00000,0.00000,1.00000}%

\pgfplotsset{every axis label/.append style={font=\footnotesize},
every tick label/.append style={font=\footnotesize}
}

\begin{tikzpicture}[font=\footnotesize] 

\begin{axis}[%
name=ber,
ymode=log,
width  = 0.35\columnwidth,
height = 0.3\columnwidth,
scale only axis,
xmin  = 1,
xmax  = 14,
xlabel= {$n(\times 1000)$},
xmajorgrids,
ymin = 0.04,
ymax = 432,
xtick       ={2, 4, 6, 8, 10, 12},
xticklabels ={$2$, $4$, $6$, $8$, $10$,$12$},
ylabel={Time (seconds)},
ymajorgrids,
]
\addplot+[smooth,color=red,solid, every mark/.append style={solid}, mark=none]
table[row sep=crcr]{
1	0.23111422856648  \\ 
2	1.12210321426392  \\ 
3	3.19470254580180  \\ 
4	7.63863476117452  \\ 
5	16.5813179016113  \\ 
6	23.0189038117727  \\ 
7	38.4196780522664  \\ 
8	64.6306711037954  \\ 
9	101.472688515981  \\ 
10	130.772794405619  \\ 
11	173.554800907771  \\ 
12	230.248995463053  \\ 
13	367.878076553345  \\ 
14	431.642583926519  \\
};

\addplot+[smooth,color=green,densely dashed, every mark/.append style={solid}, mark=none]
table[row sep=crcr]{
1	0.156687180201213  \\ 
2	0.960280021031698  \\ 
3	2.92720699310303  \\ 
4	4.55276989936829  \\ 
5	9.75968702634176  \\ 
6	13.4442013104757  \\ 
7	21.0232327779134  \\ 
8	40.7120460669200  \\ 
9	43.9327768484751  \\ 
10	63.4926010767619  \\ 
11	89.5302619139354  \\ 
12	134.587170918783  \\ 
13	167.095692475637  \\ 
14	220.173525412877 \\
};

\addplot+[smooth,color=black, loosely dashed, every mark/.append style={solid}, mark=none]
table[row sep=crcr]{
1	0.21700731913248  \\ 
2	1.06162118911743  \\ 
3	3.17253168423971  \\ 
4	6.54025149345398  \\ 
5	14.3351434866587  \\ 
6	22.2188227971395  \\ 
7	35.4310117562612  \\ 
8	53.6540784041087  \\ 
9	81.8360671997070  \\ 
10	107.793051799138  \\ 
11	141.227884531021  \\ 
12	197.803083578746  \\ 
13	240.480739752452  \\ 
14	303.556908925374 \\
};

\addplot+[smooth,color=blue,densely dotted, every mark/.append style={solid}, mark=diamond]
table[row sep=crcr]{
1	0.109987576802572  \\ 
2	0.595191791534424  \\ 
3	1.84293031692505  \\ 
4	2.70693771044413  \\ 
5	4.06979529062907  \\ 
6	5.64397390683492  \\ 
7	9.44030245145162  \\ 
8	14.2450714111328  \\ 
9	15.1673905849457  \\ 
10	22.5795931021373  \\ 
11	26.4160255591075  \\ 
12	39.3410653273265  \\ 
13	43.6978260676066  \\ 
14	59.0897785822550 \\
};

\addplot+[smooth,color=red, loosely dashed, every mark/.append style={solid}, mark=none]
table[row sep=crcr]{
1	0.245395541191101  \\ 
2	0.550196440219879  \\ 
3	0.871009171009064  \\ 
4	1.79278862476349   \\ 
5	3.73475694656372  \\ 
6	5.33216696977615  \\ 
7	7.75970399379730  \\ 
8	11.3005527257919  \\ 
9	13.9648556113243  \\ 
10	18.3594917654991  \\ 
11	25.0968192219734  \\ 
12	31.0597548484802  \\ 
13	37.7035253047943  \\ 
14	43.0544623732567  \\
};

\addplot+[smooth,color=green, densely dotted, every mark/.append style={thick}, mark=none]
table[row sep=crcr]{
1	0.333167402744293  \\ 
2	1.39640289545059  \\ 
3	1.94664639234543  \\ 
4	2.19278862476349  \\ 
5	2.69051325321198  \\ 
6	4.49012720584869  \\ 
7	6.20171004533768  \\ 
8	9.35624134540558  \\ 
9	12.1523290872574  \\ 
10	16.0577186942101  \\ 
11	20.4775619506836  \\ 
12	25.5900623202324  \\ 
13	30.7052639722824  \\ 
14	37.1392313838005  \\
};

\addplot+[smooth,color=black,loosely dotted, every mark/.append style={thick}, mark=|]
table[row sep=crcr]{
1	0.0756642222404480  \\ 
2	0.314431369304657  \\ 
3	0.751771688461304  \\ 
4	1.18167871236801  \\ 
5	2.17550402879715  \\ 
6	3.16158014535904  \\ 
7	4.76860356330872  \\ 
8	6.30832117795944  \\ 
9	7.16019552946091  \\ 
10	12.7948048114777  \\ 
11	13.9061313271523  \\ 
12	18.2098318338394  \\ 
13	21.5842366814613  \\ 
14	26.5475451350212 \\
};

\addplot+[smooth,color=blue,loosely dotted, every mark/.append style={thick}, mark=-]
table[row sep=crcr]{
1	0.0421631336212158  \\ 
2	0.161519110202789  \\ 
3	0.648253650665283  \\ 
4	0.671617536544800  \\ 
5	1.42878466844559  \\ 
6	2.00366634130478  \\ 
7	2.81753164529800  \\ 
8	3.53539955615997  \\ 
9	3.70717144012451  \\ 
10	4.60247898101807  \\ 
11	6.06118309497833  \\ 
12	7.47664237022400  \\ 
13	8.63035005331039  \\ 
14	10.1178569793701  \\  
};
\end{axis}

\begin{axis}[%
name=SumRate,
at={($(ber.east)+(35,0em)$)},
		anchor= west,
ymode=log,
width  = 0.35\columnwidth,
height = 0.3\columnwidth,
scale only axis,
xmin   = 1,
xmax  = 14,
xlabel= {$n(\times 1000)$},
xmajorgrids,
ymin = 0.05,
ymax = 473,
xtick       ={2, 4, 6, 8, 10, 12},
xticklabels ={$2$, $4$, $6$, $8$, $10$,$12$},
ymajorgrids,
]

\addplot+[smooth,color=red,solid, every mark/.append style={solid}, mark=none]
table[row sep=crcr]{
1	0.25099269549051  \\
2	1.26694544156392  \\
3	4.33468794822693  \\
4	9.50582583745321  \\
5	19.6696107387543  \\
6	35.9970933596293  \\
7	60.8242431481679  \\
8	101.595100482305  \\
9	133.822990735372  \\
10	190.082009394964  \\
11	244.314423163732  \\
12	329.381257692973  \\
13	396.833296537399  \\
14	472.576936324437 \\
};

\addplot+[smooth,color=green, densely dashed, every mark/.append style={solid}, mark=none]
table[row sep=crcr]{
1	0.193641980489095  \\
2	1.07567238807678  \\
3	3.42077763875326  \\
4	5.97618643442790  \\
5	10.1264930566152  \\
6	19.3680179119110  \\
7	32.5060058434804  \\
8	52.8395702044169  \\
9	66.7816189130147  \\
10	95.1668185393016  \\
11	129.608755985896  \\
12	165.137478748957  \\
13	207.691659291585  \\
14	265.454431533814  \\
};

\addplot+[smooth,color=black, loosely dashed, every mark/.append style={solid}, mark=none]
table[row sep=crcr]{
1	0.234918753306071  \\
2	1.21418937047323  \\
3	4.11791499455770  \\
4	8.78700812657674  \\
5	16.7504777908325  \\  
6	30.6614887714386  \\
7	60.5278569857279  \\
8	86.9537363847097  \\
9	115.700017770131  \\
10	163.268459399541  \\
11	230.721521933874  \\
12	319.677340745926  \\
13	367.451045433680  \\
14	431.037234147390 \\
};

\addplot+[smooth,color=blue,densely dotted, every mark/.append style={solid}, mark=diamond]
table[row sep=crcr]{
1	0.132465839385986  \\
2	0.695889234542847  \\
3	2.07429957389832  \\
4	2.98139405250549  \\
5	4.37069869041443  \\ 
6	7.16269421577454  \\
7	10.2984326680501  \\
8	16.5640488465627  \\
9	18.7837066650391  \\
10	26.1070335706075  \\
11	32.9518339633942  \\
12	44.0057687759399  \\
13	47.9322684605916  \\
14	65.3545192877452 \\
};

\addplot+[smooth,color=red, loosely dashed, every mark/.append style={solid}, mark=none]
table[row sep=crcr]{
1	0.259142875671387  \\
2	0.449689865112305  \\
3	1.29865723848343  \\
4	2.25757229328156  \\
5	4.29255384206772  \\  
6	6.53934776782990  \\
7	9.36613321304321  \\
8	13.1606768369675  \\
9	17.9273729920387  \\
10	22.4478815793991  \\
11	30.0332748889923  \\
12	38.2197369337082  \\
13	49.1457231044769  \\  
14	61.0274428129196  \\
};

\addplot+[smooth,color=green, densely dotted, every mark/.append style={thick}, mark=none]
table[row sep=crcr]{
1	0.285036504268646  \\
2	1.26860619068146  \\
3	1.38125609874725  \\
4	2.97145175933838  \\
5	5.24847275018692  \\
6	6.02830570936203  \\
7	9.10624164342880  \\
8	12.6070761680603  \\
9	17.0404813885689  \\
10	21.6354061365128  \\
11	29.0796950459480  \\
12	37.3066083192825  \\
13	46.4502519369125  \\
14	55.4786750078201 \\
};

\addplot+[smooth,color=black,loosely dotted, every mark/.append style={thick}, mark=|]
table[row sep=crcr]{
1	0.110137462615967  \\
2	0.381856560707092  \\
3	0.914300382137299  \\
4	2.05570876598358  \\
5	3.76062101125717  \\
6	5.10762107372284  \\
7	7.19917398691177  \\
8	9.69020313024521  \\
9	13.0282727479935  \\
10	17.0189938545227  \\
11	22.2759737968445  \\
12	29.8674979209900  \\
13	45.8060505390167  \\
14	54.9914402365685 \\  
};

\addplot+[smooth,color=blue,loosely dotted, every mark/.append style={thick}, mark=-]
table[row sep=crcr]{
1	0.053253650665283  \\
2	0.210659027099609  \\
3	0.593935072422028  \\
4	1.46585547924042  \\
5	2.41070955991745  \\
6	3.99012118577957  \\
7	5.41569393873215  \\
8	6.57656985521317  \\
9	9.11182826757431  \\  
10	11.4268012046814  \\
11	15.0735453367233  \\
12	20.3594583272934  \\
13	23.0932825803757  \\
14	28.7135926485062  \\
};
\end{axis}

\end{tikzpicture}%
\captionsetup{justification=centering,font=scriptsize}  
\caption{\textcolor{black}{Computational times for decomposition of sparse matrices with $d=0.3n$. Left: $q=0$. Right: $q=2$.}} 
\label{figSsp015}      
\end{center}
\end{figure}
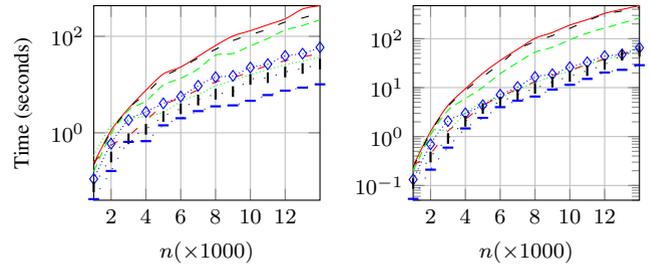

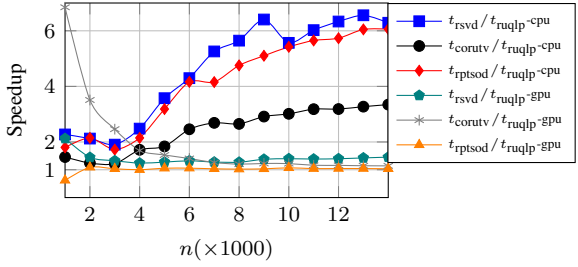
\begin{figure}[t]
\begin{center}       
%
%
%
\usetikzlibrary{positioning,calc}

\definecolor{mycolor1}{rgb}{0.00000,1.00000,1.00000}%
\definecolor{mycolor2}{rgb}{1.00000,0.00000,1.00000}%

\pgfplotsset{every axis label/.append style={font=\footnotesize},
every tick label/.append style={font=\footnotesize}
}

\begin{tikzpicture}[font=\footnotesize] 

\begin{axis}[%
name=ber,
width  = 0.5\columnwidth,
height = 0.3\columnwidth,
scale only axis,
xmin  = 1,
xmax  = 14,
xlabel= {$n(\times 1000)$},
xmajorgrids,
xtick       ={2, 4, 6, 8, 10, 12},
xticklabels ={$2$, $4$, $6$, $8$, $10$,$12$},
ymin = 0,
ymax = 7,
ylabel={Speedup},
ymajorgrids,
ytick       ={1, 2, 4, 6},
yticklabels ={$1$, $2$, $4$, $6$},
legend entries={$t_\text{rsvd}/t_\text{ruqlp}$-cpu,$t_\text{corutv}/t_\text{ruqlp}$-cpu, $t_\text{rptsod}/t_\text{ruqlp}$-cpu, $t_\text{rsvd}/t_\text{ruqlp}$-gpu, $t_\text{corutv}/t_\text{ruqlp}$-gpu, $t_\text{rptsod}/t_\text{ruqlp}$-gpu},
legend style={at={(1.58,1)},anchor=north east,draw=black,fill=white,legend cell align=left,font=\tiny},
]

\addplot[smooth,mark=square*,blue]
table[row sep=crcr]{
1	2.27060000000000  \\
2	2.12270000000000 \\
3	1.90390000000000 \\
4	2.48060000000000 \\
5	3.57530000000000 \\
6	4.29600000000000 \\
7	5.25760000000000 \\
8	5.64600000000000 \\
9	6.40670000000000 \\
10	5.56610000000000 \\
11	6.02340000000000 \\
12	6.33420000000000 \\
13	6.55880000000000 \\
14	6.28910000000000 \\
};

\addplot[smooth,mark=*,black]
table[row sep=crcr]{
 \\
1	1.45790000000000  \\
2	1.24160000000000 \\
3	1.21670000000000 \\
4	1.71800000000000 \\
5	1.84010000000000 \\
6	2.45700000000000 \\
7	2.68770000000000 \\
8	2.65030000000000 \\
9	2.91420000000000 \\
10	3.01320000000000 \\
11	3.18050000000000 \\
12	3.18480000000000 \\
13	3.27090000000000 \\
14	3.34460000000000 \\
};

\addplot[smooth,mark=diamond*,red]
table[row sep=crcr]{
1	1.80550000000000  \\
2	2.14570000000000 \\
3	1.72550000000000 \\
4	2.14570000000000 \\
5	3.18440000000000 \\
6	4.16210000000000 \\ 
7	4.15600000000000 \\
8	4.75580000000000 \\
9	5.10170000000000 \\
10	5.42210000000000 \\
11	5.65130000000000 \\
12	5.73330000000000 \\
13	6.04990000000000 \\
14	6.06090000000000 \\
};

\addplot[smooth,mark=pentagon*, teal]
table[row sep=crcr]{
1	2.12470000000000  \\
2	1.44470000000000 \\
3	1.32390000000000 \\
4	1.24030000000000 \\
5	1.27580000000000 \\
6	1.32230000000000 \\
7	1.27710000000000 \\
8	1.27170000000000 \\
9	1.37770000000000 \\
10	1.39750000000000 \\
11	1.38550000000000 \\
12	1.39670000000000 \\
13	1.42660000000000 \\
14	1.45460000000000 \\
};

\addplot[smooth,mark=asterisk, gray]
table[row sep=crcr]{
1	6.84670000000000   \\
2	3.51030000000000  \\
3	2.46080000000000  \\
4	1.69340000000000  \\
5	1.53140000000000  \\
6	1.40560000000000  \\
7	1.25910000000000  \\
8	1.20270000000000  \\
9	1.22670000000000  \\
10	1.22090000000000  \\
11	1.16200000000000  \\
12	1.15810000000000  \\
13	1.14470000000000  \\
14	1.14040000000000  \\
};

\addplot[smooth,mark=triangle*, orange]
table[row sep=crcr]{
1	0.62030000000000 \\
2	1.08960000000000 \\
3	1.03560000000000 \\
4	1.00290000000000 \\
5	1.05500000000000 \\
6	1.05730000000000 \\
7	1.03470000000000 \\
8	1.02350000000000 \\
9	1.03300000000000 \\
10	1.07010000000000 \\
11	1.04090000000000 \\
12	1.04520000000000 \\
13	1.04360000000000 \\
14	1.03610000000000 \\
};

\end{axis}

\end{tikzpicture}%
\captionsetup{justification=centering,font=scriptsize}  
\caption{Speedup for dense matrices with $d=0.04n$ with $q=0$.} 
\label{figSPUzzt}      
\end{center}
\end{figure}

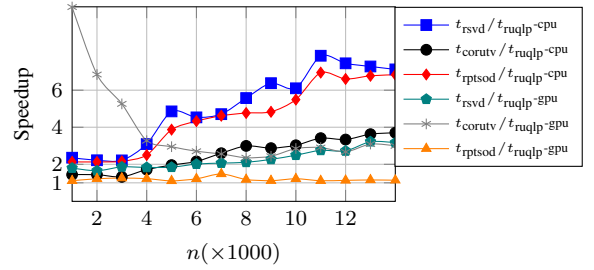
\begin{figure}[t]
\begin{center}       
%
%
%
\usetikzlibrary{positioning,calc}

\definecolor{mycolor1}{rgb}{0.00000,1.00000,1.00000}%
\definecolor{mycolor2}{rgb}{1.00000,0.00000,1.00000}%

\pgfplotsset{every axis label/.append style={font=\footnotesize},
every tick label/.append style={font=\footnotesize}
}

\begin{tikzpicture}[font=\footnotesize] 

\begin{axis}[%
name=ber,
width  = 0.5\columnwidth,
height = 0.3\columnwidth,
scale only axis,
xmin  = 1,
xmax  = 14,
xlabel= {$n(\times 1000)$},
xmajorgrids,
xtick       ={2, 4, 6, 8, 10, 12},
xticklabels ={$2$, $4$, $6$, $8$, $10$,$12$},
ymin = 0,
ymax = 10.5,
ylabel={Speedup},
ymajorgrids,
ytick       ={1, 2, 4, 6},
yticklabels ={$1$, $2$, $4$, $6$},
legend entries={$t_\text{rsvd}/t_\text{ruqlp}$-cpu,$t_\text{corutv}/t_\text{ruqlp}$-cpu, $t_\text{rptsod}/t_\text{ruqlp}$-cpu, $t_\text{rsvd}/t_\text{ruqlp}$-gpu, $t_\text{corutv}/t_\text{ruqlp}$-gpu, $t_\text{rptsod}/t_\text{ruqlp}$-gpu},
legend style={at={(1.58,1)},anchor=north east,draw=black,fill=white,legend cell align=left,font=\tiny},
]

\addplot[smooth,mark=square*,blue]
table[row sep=crcr]{
1	2.35494166100000 \\
2	2.21702298300000 \\
3	2.21111161800000 \\
4	3.09129128200000 \\
5	4.85548739400000 \\
6	4.53254795300000 \\
7	4.70597466000000 \\
8	5.57394864600000 \\
9	6.37856473900000 \\
10	6.10174049200000 \\
11	7.85871890000000 \\
12	7.44994961400000 \\
13	7.27942508800000 \\
14	7.12192768900000 \\
};

\addplot[smooth,mark=*,black]
table[row sep=crcr]{
 \\
1	1.43724576300000  \\
2	1.44082090900000 \\
3	1.31970343600000 \\
4	1.71443331500000 \\
5	1.95648207300000 \\
6	2.15429488600000 \\
7	2.59802284500000 \\
8	2.99248479300000 \\
9	2.87231723300000 \\
10	3.02646067100000 \\
11	3.41144953000000 \\
12	3.33515835700000 \\
13	3.62606453400000 \\
14	3.7030938300000 \\
};

\addplot[smooth,mark=diamond*,red]
table[row sep=crcr]{
1	2.13128301800000  \\
2	2.12651534400000 \\
3	2.16178834300000 \\
4	2.49961264200000 \\
5	3.87050139200000 \\
6	4.33003283700000 \\
7	4.62144560700000 \\
8	4.76764836200000 \\
9	4.83224146100000 \\
10	5.49107373200000 \\ 
11	6.94918506600000 \\
12	6.60526902900000 \\
13	6.7758690200000 \\
14	6.83281941000000 \\
};

\addplot[smooth,mark=pentagon*, teal]
table[row sep=crcr]{
1	1.79804735300000  \\
2	1.64927733100000 \\
3	1.85951286500000 \\
4	1.83092490200000 \\
5	1.84285199600000 \\
6	2.01071245800000 \\
7	2.06533495200000 \\
8	2.10804990100000 \\
9	2.26569538800000 \\
10	2.49118570900000 \\
11	2.75299921500000 \\
12	2.72187437500000 \\
13	3.22048139700000 \\
14	3.16468086200000  \\
};

\addplot[smooth,mark=asterisk, gray]
table[row sep=crcr]{
1	10.4877827400000  \\
2	6.84556918000000 \\
3	5.26139629700000 \\
4	3.27245453700000 \\
5	2.95023278100000 \\
6	2.70399917400000 \\
7	2.55860565600000 \\
8	2.35519399400000 \\
9	2.43407116500000 \\
10	2.85140786900000 \\
11	2.90521987600000 \\
12	2.66837635800000 \\
13	3.09432489100000 \\
14	2.99497883000000  \\
};

\addplot[smooth,mark=triangle*, orange]
table[row sep=crcr]{
1	1.11134314900000  \\
2	1.22090131900000 \\
3	1.25532560700000 \\
4	1.22437772300000 \\
5	1.10686450400000 \\
6	1.21706151700000 \\
7	1.48143799600000 \\
8	1.17837082200000 \\
9	1.11644873800000 \\
10	1.21831804900000 \\
11	1.11308072100000 \\
12	1.13257841000000 \\
13	1.15517518300000 \\
14	1.13940677000000 \\
};

\end{axis}

\end{tikzpicture}%
\captionsetup{justification=centering,font=scriptsize}  
\caption{\textcolor{black}{Speedup for sparse matrices with $d=0.04n$ with $q=0$.}}
\label{figSspup002}      
\end{center}
\end{figure}

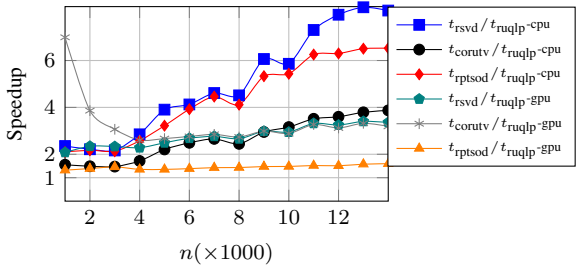
\begin{figure}[t]
\begin{center}       
%
%
%
\usetikzlibrary{positioning,calc}

\definecolor{mycolor1}{rgb}{0.00000,1.00000,1.00000}%
\definecolor{mycolor2}{rgb}{1.00000,0.00000,1.00000}%

\pgfplotsset{every axis label/.append style={font=\footnotesize},
every tick label/.append style={font=\footnotesize}
}

\begin{tikzpicture}[font=\footnotesize] 

\begin{axis}[%
name=ber,
width  = 0.5\columnwidth,
height = 0.3\columnwidth,
scale only axis,
xmin  = 1,
xmax  = 14,
xlabel= {$n(\times 1000)$},
xmajorgrids,
xtick       ={2, 4, 6, 8, 10, 12},
xticklabels ={$2$, $4$, $6$, $8$, $10$,$12$},
ymin = 0.0,
ymax = 8.3,
ylabel={Speedup},
ymajorgrids,
ytick       ={1, 2, 4, 6},
yticklabels ={$1$, $2$, $4$, $6$},
legend entries={$t_\text{rsvd}/t_\text{ruqlp}$-cpu,$t_\text{corutv}/t_\text{ruqlp}$-cpu, $t_\text{rptsod}/t_\text{ruqlp}$-cpu, $t_\text{rsvd}/t_\text{ruqlp}$-gpu, $t_\text{corutv}/t_\text{ruqlp}$-gpu, $t_\text{rptsod}/t_\text{ruqlp}$-gpu},
legend style={at={(1.58,1)},anchor=north east,draw=black,fill=white,legend cell align=left,font=\tiny},
]

\addplot[smooth,mark=square*,blue]
table[row sep=crcr]{
1	2.35167217000000  \\
2	2.21355531100000 \\
3	2.16230501300000 \\
4	2.85144543300000 \\
5	3.90428684600000 \\
6	4.12544735700000 \\
7	4.61256508800000 \\
8	4.51392008400000 \\
9	6.06411148500000 \\
10	5.85990166300000 \\
11	7.30270553900000 \\
12	7.95417905000000 \\
13	8.27183803200000 \\
14	8.13221428000000 \\
};

\addplot[smooth,mark=*,black]
table[row sep=crcr]{
 \\
1	1.55467784000000 \\
2	1.48824993900000 \\
3	1.47384696100000 \\
4	1.71436449300000 \\
5	2.20571136400000 \\
6	2.48906788900000 \\
7	2.65629089900000 \\
8	2.44425397400000 \\
9	2.95091969200000 \\
10	3.16045034700000 \\
11	3.51390889100000 \\
12	3.59987744900000 \\
13	3.79058073900000 \\
14	3.87272799500000 \\
};

\addplot[smooth,mark=diamond*,red]
table[row sep=crcr]{
1	2.09438686800000 \\
2	2.17184611100000 \\
3	2.12993967900000 \\
4	2.57693023700000 \\
5	3.21965870900000 \\
6	3.92280201900000 \\
7	4.46149370300000 \\
8	4.11734529200000 \\
9	5.33235780700000 \\
10	5.43359332100000 \\
11	6.25472112300000 \\
12	6.30082733600000 \\
13	6.50549649900000 \\
14	6.52633350900000 \\
};

\addplot[smooth,mark=pentagon*, teal]
table[row sep=crcr]{
1	2.07827791900000 \\
2	2.34372833800000 \\
3	2.33140249700000 \\
4	2.27928916400000 \\
5	2.48899060100000 \\
6	2.68359343500000 \\
7	2.76581532100000 \\
8	2.65959451700000 \\
9	2.98401139400000 \\
10	2.94465126400000 \\
11	3.31777441600000 \\
12	3.24694699700000 \\
13	3.40227039100000 \\
14	3.35889504200000 \\
};

\addplot[smooth,mark=asterisk, gray]
table[row sep=crcr]{
1	6.99047026300000  \\
2	3.86169178300000 \\
3	3.06619978800000 \\
4	2.59945546400000 \\
5	2.66870332800000 \\
6	2.76008148100000 \\
7	2.85014103600000 \\
8	2.73445992300000 \\
9	2.98979310100000 \\
10	2.88340112700000 \\
11	3.27144144900000 \\
12	3.12160336200000 \\
13	3.32174673900000 \\
14	3.19272834900000  \\
};

\addplot[smooth,mark=triangle*, orange]
table[row sep=crcr]{
1	1.32088252200000  \\
2	1.38994336000000 \\
3	1.46189478400000 \\
4	1.35596967100000 \\
5	1.35468391800000 \\
6	1.38662853100000 \\
7	1.42892536100000 \\
8	1.43369477800000 \\
9	1.47585673300000 \\
10	1.47973540900000 \\
11	1.51957032500000 \\ 
12	1.50986892900000 \\
13	1.57670886700000 \\
14	1.59782637200000 \\
};

\end{axis}

\end{tikzpicture}%
\captionsetup{justification=centering,font=scriptsize}  
\caption{\textcolor{black}{Speedup for sparse matrices with $d=0.2n$ with $q=0$.}} 
\label{figSspup01}      
\end{center}
\end{figure}

\begin{figure}[t]
\begin{center}       
%
%
%
\usetikzlibrary{positioning,calc}

\definecolor{mycolor1}{rgb}{0.00000,1.00000,1.00000}%
\definecolor{mycolor2}{rgb}{1.00000,0.00000,1.00000}%

\pgfplotsset{every axis label/.append style={font=\footnotesize},
every tick label/.append style={font=\footnotesize}
}

\begin{tikzpicture}[font=\footnotesize] 

\begin{axis}[%
name=ber,
width  = 0.5\columnwidth,
height = 0.3\columnwidth,
scale only axis,
xmin  = 1,
xmax  = 14,
xlabel= {$n(\times 1000)$},
xmajorgrids,
xtick       ={2, 4, 6, 8, 10, 12},
xticklabels ={$2$, $4$, $6$, $8$, $10$,$12$},
ymin = 0.0,
ymax = 8.5,
ylabel={Speedup},
ymajorgrids,
ytick       ={1, 2, 4, 6},
yticklabels ={$1$, $2$, $4$, $6$},
legend entries={$t_\text{rsvd}/t_\text{ruqlp}$-cpu,$t_\text{corutv}/t_\text{ruqlp}$-cpu, $t_\text{rptsod}/t_\text{ruqlp}$-cpu, $t_\text{rsvd}/t_\text{ruqlp}$-gpu, $t_\text{corutv}/t_\text{ruqlp}$-gpu, $t_\text{rptsod}/t_\text{ruqlp}$-gpu},
legend style={at={(1.58,1)},anchor=north east,draw=black,fill=white,legend cell align=left,font=\tiny},
]

\addplot[smooth,mark=square*,blue]
table[row sep=crcr]{
1	2.10127575600000  \\
2	1.88528005600000 \\
3	1.73349069000000 \\
4	2.82187311900000 \\
5	4.07423880500000 \\
6	4.07849224500000 \\
7	4.06975075700000 \\
8	4.53705490400000 \\
9	6.69018760700000 \\
10	5.79163644900000 \\
11	6.57005727500000 \\
12	7.85263753100000 \\
13	8.41868142300000 \\
14	8.0486040500000 \\
};

\addplot[smooth,mark=*,black]
table[row sep=crcr]{
 \\
1	1.42458980100000 \\
2	1.61339594200000 \\
3	1.58834382700000 \\
4	1.68188942100000 \\
5	2.39807811700000 \\
6	2.38204526300000 \\
7	2.52696602000000 \\
8	2.85797416500000 \\
9	2.89652835200000 \\
10	2.81194620200000 \\
11	3.38924043300000 \\
12	3.42103524100000 \\
13	3.82389028300000 \\
14	3.92608479300000 \\
};

\addplot[smooth,mark=diamond*,red]
table[row sep=crcr]{
1	1.97301663900000 \\ 
2	1.78366234900000 \\
3	1.72146046700000 \\
4	2.41610712700000 \\
5	3.52232544000000 \\
6	3.93673379100000 \\
7	3.75316489500000 \\
8	3.76650118900000 \\
9	5.38552711700000 \\
10	5.397391471600000 \\
11	5.4629572500000 \\
12	5.7790359000000 \\
13	5.80326552600000 \\
14	6.093721520400000 \\ 
};

\addplot[smooth,mark=pentagon*, teal]
table[row sep=crcr]{
1	5.82014476000000  \\
2	3.40638602800000 \\
3	1.34362401200000 \\
4	2.66935946000000 \\
5	2.61393968500000 \\
6	2.66120504200000 \\
7	2.75407873700000 \\
8	3.19640044800000 \\
9	3.76698403000000 \\
10	3.98904413100000 \\
11	4.14058094400000 \\
12	4.15423840100000 \\
13	4.36871332800000 \\
14	4.25529462000000 \\
};

\addplot[smooth,mark=asterisk, gray]
table[row sep=crcr]{
1	7.90186530600000 \\
2	5.24543454800000 \\
3	3.00290849200000 \\
4	2.26493652300000 \\
5	1.88307819400000 \\
6	2.24095554900000 \\
7	2.30111460200000 \\
8	2.64644524500000 \\
9	3.27805964300000 \\
10	3.48892819700000 \\
11	3.37847605500000 \\
12	3.42266769700000 \\
13	3.55782370200000 \\
14	3.67066182700000 \\
};

\addplot[smooth,mark=triangle*, orange]
table[row sep=crcr]{
1	1.79455879400000  \\
2	1.94671311000000 \\
3	1.8268755100000 \\
4	1.75945184300000 \\
5	1.52262554100000 \\ 
6	1.57789751700000 \\
7	1.69247560100000 \\
8	1.78433047700000 \\
9	1.93144440300000 \\
10	2.77998115000000 \\
11	2.29429322800000 \\
12	2.43556277400000 \\
13	2.50096885400000 \\
14	2.62383083600000 \\
};

\end{axis}

\end{tikzpicture}%
\captionsetup{justification=centering,font=scriptsize}  
\caption{\textcolor{black}{Speedup for sparse matrices with $d=0.3n$ with $q=0$.}} 
\label{figSspup015}      
\end{center}
\end{figure}
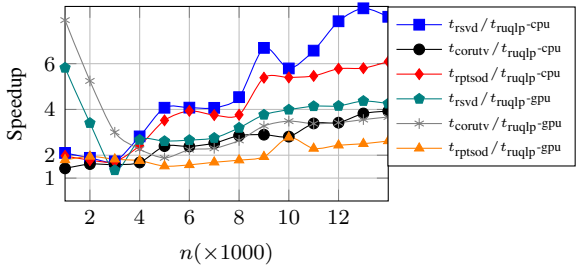

\textcolor{black}{We generate both dense and sparse matrices of size $n\times n$. A dense matrix is a matrix in which most of its elements are non-zero, whereas most of the elements of a sparse are zero. In our simulations, the sparse matrices have $0.1\times n^2$ non-zero entries.} We use multiple randomized methods to factorize the matrices into three components. Methods considered include R-SVD (randomized SVD) \cite{HMT2011}, CoR-UTV (compressed randomized UTV) \cite{MFKDeJSTSP18}, RP-TSOD (randomized pivoted two-sided orthogonal decomposition) \cite{MFKC20}, and RU-QLP. For the sampling size parameter $d$ (i.e., the dimension of the reduced matrix), we consider three cases: $d=0.04n$, $d=0.2n$ and $d=0.3n$. The results averaged over 10 independent  trials are displayed in Figures \ref{figSzzt}--\ref{figSsp015}. Furthermore, Figures \ref{figSPUzzt}--\ref{figSspup015} display the speedups offered by RU-QLP. We make three observations: 
\begin{enumerate}
\item On the CPU, RU-QLP is substantially faster than other methods. \textcolor{black}{For the decomposition of the dense matrices, RU-QLP achieves speedups of up to 7.1 times, 3.6 times, and 6.5 times over R-SVD, CoR-UTV, and RP-TSOD, respectively, whereas on the sparse matrices, it achieves speedups of up to 8.5 times, 10.5 times, and 6.9 times over R-SVD, CoR-UTV, and RP-TSOD, respectively.} The discrepancy is pronounced particularly when $q=0$ (no power iteration case). This is due to the fact that for $q\ge 1$, matrix-matrix multiplications form the dominant cost; see Section \ref{subsubImOr}. In addition, by increasing the dimension of the input matrix as well as the parameter $d$, RU-QLP confers more advantage. 

\item On the CPU-GPU architecture, in spite of the reduced gaps in computational time, RU-QLP still outperforms other methods. \textcolor{black}{For the decomposition of the dense matrices, RU-QLP achieves speedups of up to 2.3 times, 6.8 times, and 1.5 times over R-SVD, CoR-UTV, and RP-TSOD, respectively, while on the sparse matrices, it achieves speedups of up to 5.8 times, 8.6 times, and 2.6 times over R-SVD, CoR-UTV, and RP-TSOD, respectively.} Moreover, by increasing the dimension of the input matrix and particularly $d$, RU-QLP starts to become more efficient.

\item RU-QLP is faster than CoR-UTV, RP-TSOD and R-SVD, because it only makes use of unpivoted QR, which leverages the parallel architecture of system better than pivoted QR and the SVD utilized in other methods. CoR-UTV and RP-TSOD outperform R-SVD, because they utilize pivoted QR, which imposes less communication cost upon the system than the SVD; see Section \ref{subsubComCost}.
\end{enumerate}

\begin{figure}[t]
\begin{center}       
%
%
%
\usetikzlibrary{positioning,calc}

\definecolor{mycolor1}{rgb}{0.00000,1.00000,1.00000}%
\definecolor{mycolor2}{rgb}{1.00000,0.00000,1.00000}%

\pgfplotsset{every axis label/.append style={font=\footnotesize},
every tick label/.append style={font=\footnotesize}
}

\begin{tikzpicture}[font=\footnotesize] 

\begin{axis}[%
name=ber,
ymode=log,
width  = 0.35\columnwidth,
height = 0.3\columnwidth,
scale only axis,
xmin  = 1,
xmax  = 14,
xlabel= {$n(\times 1000)$},
xmajorgrids,
ymin = 0.02,
ymax = 392,
xtick       ={2, 4, 6, 8, 10, 12},
xticklabels ={$2$, $4$, $6$, $8$, $10$,$12$},
ylabel={Time (seconds)},
ymajorgrids,
]
\addplot+[smooth,color=red,solid, every mark/.append style={solid}, mark=none]
table[row sep=crcr]{
1	0.216599623362223 \\
2	0.906644503275553 \\
3	2.63580242792765 \\
4	7.11749211947123 \\
5	15.4540585676829 \\
6	22.6437668005625 \\
7	36.6397174199422 \\
8	62.4679939746857 \\
9	94.3449877103170 \\
10	113.808926423391 \\
11	171.484635035197 \\
12	229.434520403544 \\
13	313.106295188268 \\
14	391.822401205699  \\
};

\addplot+[smooth,color=green,densely dashed, every mark/.append style={solid}, mark=none]
table[row sep=crcr]{
1	0.143192847569784 \\
2	0.609568516413371 \\
3	1.79658714930216 \\
4	4.27922471364339 \\
5	8.73070907592773 \\
6	12.5642421245575 \\
7	21.1001354058584 \\ 
8	33.8259516557058 \\
9	43.7460095882416 \\
10	62.4468058745066 \\
11	78.2298350334168 \\
12	129.983141740163 \\
13	160.376160144806 \\
14	190.633227427801 \\
};

\addplot+[smooth,color=black, loosely dashed, every mark/.append style={solid}, mark=none]
table[row sep=crcr]{
1	0.192902485529582 \\
2	0.889560937881470 \\
3	2.59634979565938 \\
4	6.43227481842041 \\
5	12.7441441218058 \\
6	21.5314865112305 \\
7	35.4396882057190 \\
8	53.3816766738892 \\
9	79.0497201283773 \\
10	107.361454884211 \\
11	139.248289267222 \\
12	191.400457541148 \\
13	237.628370920817 \\
14	301.054857889811 \\
};

\addplot+[smooth,color=blue,densely dotted, every mark/.append style={solid}, mark=diamond]
table[row sep=crcr]{
1	0.0921045144399007 \\
2	0.409587462743123 \\
3	1.21897808710734 \\
4	2.49609971046448 \\
5	3.95822827021281 \\ 
6	5.48880275090535 \\
7	7.94345808029175 \\
8	13.8389676411947 \\ 
9	14.8245340983073 \\
10	19.7588315010071 \\
11	22.2629093329112 \\
12	36.1076574325562 \\
13	36.5273228486379 \\
14	56.5219690799713 \\
};

\addplot+[smooth,color=red, loosely dashed, every mark/.append style={solid}, mark=none]
table[row sep=crcr]{
1	0.0515531897544861   \\
2	0.172107100486755  \\
3	0.383600473403931  \\
4	0.758249163627625  \\
5	1.23849165439606  \\
6	1.93883186578751  \\
7	2.78125238418579  \\
8	3.86143469810486  \\
9	5.36160111427307  \\
10	7.01116669178009  \\
11	9.59535765647888  \\
12	12.0736199617386  \\
13	14.7113977074623  \\
14	18.1352038383484  \\
};

\addplot+[smooth,color=green, densely dotted, every mark/.append style={thick}, mark=none]
table[row sep=crcr]{
1	0.173403680324554  \\
2	0.283575773239136  \\
3	0.504501342773438  \\
4	0.864758610725403  \\
5	1.32791453599930  \\
6	1.99409264326096  \\
7	2.86604875326157  \\
8	3.97013092041016  \\  
9	5.37198954820633  \\
10	6.86533111333847  \\
11	9.46135777235031  \\
12	11.6075355410576  \\
13	14.3632138967514  \\
14	17.2380436658859  \\
};

\addplot+[smooth,color=black,loosely dotted, every mark/.append style={thick}, mark=|]
table[row sep=crcr]{
1	0.0327654480934143  \\
2	0.102067768573761  \\
3	0.240534842014313  \\
4	0.451089262962341  \\
5	0.674074351787567  \\
6	1.00180584192276  \\
7	1.43690073490143  \\
8	2.08156496286392  \\
9	2.65178447961807  \\
10	3.52322590351105  \\
11	4.39475953578949  \\  
12	5.61437672376633  \\
13	6.81768012046814  \\
14	8.62691646814346 \\
};

\addplot+[smooth,color=blue,loosely dotted, every mark/.append style={thick}, mark=-]
table[row sep=crcr]{
1	0.0248057246208191  \\
2	0.0734330415725708  \\
3	0.164536356925964  \\
4	0.332669138908386  \\
5	0.497587919235230  \\
6	0.722476005554199  \\
7	1.00558137893677  \\
8	1.45188850164413  \\
9	1.79677635431290  \\
10	2.38098371028900  \\
11	2.89210671186447  \\
12	3.71845304965973  \\
13	4.32399427890778  \\
14	5.39915764331818  \\  
};

\end{axis}

\begin{axis}[%
name=SumRate,
at={($(ber.east)+(35,0em)$)},
		anchor= west,
ymode=log,
width  = 0.35\columnwidth,
height = 0.3\columnwidth,
scale only axis,
xmin   = 1,
xmax  = 14,
xlabel= {$n(\times 1000)$},
xmajorgrids,
ymin = 0.1,
ymax = 455,
xtick       ={2, 4, 6, 8, 10, 12},
xticklabels ={$2$, $4$, $6$, $8$, $10$,$12$},
ymajorgrids,
]

\addplot+[smooth,color=red,solid, every mark/.append style={solid}, mark=none]
table[row sep=crcr]{
1	0.228997151056925   \\
2	0.996509075164795  \\
3	3.49695301055908  \\
4	9.34402028719584  \\
5	17.7295206387838  \\
6	34.6534392833710  \\
7	57.5101151466370  \\
8	95.3679498831431  \\
9	119.269742488861  \\
10	174.422742764155  \\
11	225.808374961217  \\
12	322.020395199458  \\
13	389.524038076401  \\
14	454.049478292465 \\
};

\addplot+[smooth,color=green, densely dashed, every mark/.append style={solid}, mark=none]
table[row sep=crcr]{
1	0.155032237370809  \\
2	0.674880425135295  \\
3	1.90560619036357  \\
4	5.27578409512838  \\
5	9.84230931599935  \\
6	19.2831125259399  \\
7	32.1862001419067  \\
8	52.2680186430613  \\
9	63.4753268559774  \\
10	92.4619375069936  \\
11	123.581212679545  \\
12	157.236850579580  \\
13	203.187651952108  \\
14	235.085053602854  \\
};

\addplot+[smooth,color=black, loosely dashed, every mark/.append style={solid}, mark=none]
table[row sep=crcr]{
1	0.213663101196289  \\
2	0.973623116811117  \\
3	3.35575588544210  \\
4	8.74584611256917  \\
5	16.6035617192586  \\
6	27.9081540107727  \\  
7	40.0139559110006  \\
8	74.0913818677266  \\
9	95.6090792020162  \\
10	143.660519202550  \\
11	176.292346398036  \\
12	219.006923039754  \\
13	330.966747522354  \\
14	398.342339992523 \\
};

\addplot+[smooth,color=blue,densely dotted, every mark/.append style={solid}, mark=diamond]
table[row sep=crcr]{
1	0.102909167607625  \\
2	0.494200785954793  \\
3	1.26625212033590  \\
4	2.87498704592387  \\
5	4.12778099377950  \\
6	6.99381677309672  \\
7	10.2325121561686  \\
8	16.0326174100240  \\
9	17.5572438240051  \\
10	24.9373449484507  \\
11	28.3816045125326  \\
12	43.1034887631734  \\
13	45.4502123991648  \\
14	62.4847609996796 \\
};

\addplot+[smooth,color=red, loosely dashed, every mark/.append style={solid}, mark=none]
table[row sep=crcr]{
1	0.399742782115936  \\
2	0.972331643104553  \\
3	1.85018938779831  \\
4	3.05579656362534  \\
5	4.82367408275604  \\
6	7.75003057718277  \\
7	9.22953844070435  \\
8	13.3200529813766  \\
9	15.2533316016197  \\
10	19.5369037985802  \\
11	24.6606346964836  \\
12	27.8683789968491  \\
13	32.4273526072502  \\
14	33.8216350078583  \\
};

\addplot+[smooth,color=green, densely dotted, every mark/.append style={thick}, mark=none]
table[row sep=crcr]{
1	0.449600398540497  \\
2	1.33798772096634  \\
3	2.82037693262100  \\
4	4.91147118806839  \\
5	5.92923688888550  \\
6	6.93352544307709  \\
7	7.99820232391357  \\
8	11.0746768116951  \\
9	13.9560215473175  \\
10	16.5712271928787  \\
11	21.9431757330894  \\
12	26.3862940073013  \\
13	31.0420438051224  \\
14	33.0896491408348 \\
};

\addplot+[smooth,color=black,loosely dotted, every mark/.append style={thick}, mark=|]
table[row sep=crcr]{
1	0.249248385429382  \\
2	0.311456441879272  \\
3	0.698092401027679  \\
4	1.34766018390656  \\
5	2.17104423046112  \\
6	3.35436999797821  \\
7	4.68715542554855  \\
8	6.66282045841217  \\
9	9.32107752561569  \\
10	12.2330427169800  \\
11	15.8780269026756  \\
12	20.3394302725792  \\
13	25.3035895824432  \\
14	30.5183458328247 \\  
};

\addplot+[smooth,color=blue,loosely dotted, every mark/.append style={thick}, mark=-]
table[row sep=crcr]{
1	0.193117439746857  \\
2	0.287205815315247  \\
3	0.461214363574982  \\
4	0.885101735591888  \\
5	1.26876288652420  \\
6	1.83378326892853  \\
7	2.68335968255997  \\
8	4.10972064733505  \\  
9	5.37851417064667  \\
10	7.23080635070801  \\    
11	9.15518110990524  \\
12	11.7156784534454  \\
13	14.3694746494293  \\
14	18.1045593619347  \\
};

\end{axis}

\end{tikzpicture}%
\captionsetup{justification=centering,font=scriptsize}  
\caption{\textcolor{black}{Computational times for decomposition of sparse matrices with $d=0.2n$. Left: $q=0$. Right: $q=2$.}}
\label{figSsp01}      
\end{center}
\end{figure}
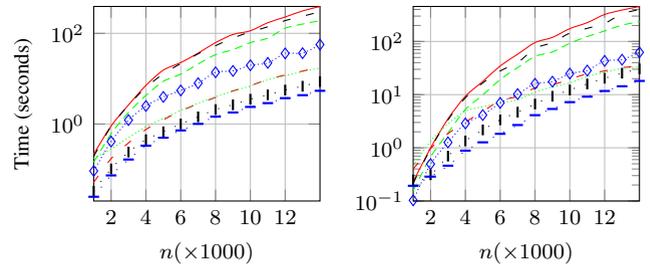

\begin{figure}[t]
\begin{center}       
%
%
%
\usetikzlibrary{positioning,calc}

\definecolor{mycolor1}{rgb}{0.00000,1.00000,1.00000}%
\definecolor{mycolor2}{rgb}{1.00000,0.00000,1.00000}%

\pgfplotsset{every axis label/.append style={font=\footnotesize},
every tick label/.append style={font=\footnotesize}
}

\begin{tikzpicture}[font=\footnotesize] 

\begin{axis}[%
name=ber,
width  = 0.5\columnwidth,
height = 0.3\columnwidth,
scale only axis,
xmin  = 1,
xmax  = 14,
xlabel= {$n(\times 1000)$},
xmajorgrids,
xtick       ={2, 4, 6, 8, 10, 12},
xticklabels ={$2$, $4$, $6$, $8$, $10$,$12$},
ymin = 0.0,
ymax = 7,
ylabel={Speedup},
ymajorgrids,
ytick       ={1, 2, 4, 6},
yticklabels ={$1$, $2$, $4$, $6$},
legend entries={$t_\text{rsvd}/t_\text{ruqlp}$-cpu,$t_\text{corutv}/t_\text{ruqlp}$-cpu, $t_\text{rptsod}/t_\text{ruqlp}$-cpu, $t_\text{rsvd}/t_\text{ruqlp}$-gpu, $t_\text{corutv}/t_\text{ruqlp}$-gpu, $t_\text{rptsod}/t_\text{ruqlp}$-gpu},
legend style={at={(1.58,1)},anchor=north east,draw=black,fill=white,legend cell align=left,font=\tiny},
]

\addplot[smooth,mark=square*,blue]
table[row sep=crcr]{
1	2.48310000000000  \\
2	2.76620000000000 \\
3	3.02370000000000 \\
4	2.98290000000000 \\
5	3.91340000000000 \\
6	4.49740000000000 \\
7	5.93560000000000 \\
8	6.33620000000000 \\
9	5.69350000000000 \\
10	5.52340000000000 \\
11	6.26480000000000 \\
12	6.89110000000000 \\
13	6.74870000000000 \\
14	6.44750000000000 \\
};

\addplot[smooth,mark=*,black]
table[row sep=crcr]{
 \\
1	1.46850000000000 \\
2	1.46210000000000\\
3	1.41280000000000\\
4	1.72590000000000\\
5	2.14540000000000\\
6	2.40380000000000\\
7	2.74300000000000\\
8	3.11710000000000\\
9	3.17140000000000\\
10	3.09280000000000\\
11	3.35200000000000\\
12	3.26430000000000\\
13	3.41610000000000\\
14	3.44910000000000 \\
};

\addplot[smooth,mark=diamond*,red]
table[row sep=crcr]{
1	2.11280000000000  \\
2	2.20640000000000 \\
3	2.10420000000000 \\
4	2.66020000000000 \\
5	3.23020000000000 \\
6	4.19960000000000 \\
7	4.32910000000000 \\
8	5.04630000000000 \\
9	5.15920000000000 \\
10	5.49890000000000 \\
11	5.83650000000000 \\
12	6.15000000000000 \\
13	6.16390000000000 \\
14	6.29940000000000 \\
};

\addplot[smooth,mark=pentagon*, teal]
table[row sep=crcr]{
1	1.87960000000000  \\
2	1.54040000000000 \\
3	1.54110000000000 \\
4	1.60140000000000 \\
5	1.66250000000000 \\
6	1.72590000000000 \\
7	1.71540000000000 \\
8	1.70050000000000 \\
9	1.76580000000000 \\
10	1.74090000000000 \\
11	1.86580000000000 \\
12	1.82860000000000 \\
13	1.86180000000000 \\
14	1.86840000000000 \\
};

\addplot[smooth,mark=asterisk, gray]
table[row sep=crcr]{
1	3.96870000000000   \\
2	2.07740000000000  \\
3	1.81170000000000  \\
4	1.71360000000000  \\
5	1.70300000000000  \\
6	1.69170000000000  \\
7	1.63610000000000  \\
8	1.62600000000000  \\
9	1.64020000000000  \\
10	1.58950000000000  \\
11	1.69140000000000  \\
12	1.65410000000000  \\
13	1.66380000000000  \\
14	1.65830000000000  \\
};

\addplot[smooth,mark=triangle*, orange]
table[row sep=crcr]{
1	0.96860000000000 \\
2	1.17130000000000 \\
3	1.12160000000000 \\
4	1.19820000000000 \\
5	1.16990000000000 \\
6	1.17340000000000 \\
7	1.15970000000000 \\
8	1.18640000000000 \\
9	1.18640000000000 \\
10	1.18480000000000 \\
11	1.22010000000000 \\
12	1.20920000000000 \\
13	1.21510000000000 \\
14	1.22080000000000 \\
};

\end{axis}

\end{tikzpicture}%
\captionsetup{justification=centering,font=scriptsize}  
\caption{Speedup for dense matrices with $d=0.2n$ with $q=0$.} 
\label{figSPUz1}      
\end{center}
\end{figure}

\begin{figure}[t]
\begin{center}       
%
%
%
\usetikzlibrary{positioning,calc}

\definecolor{mycolor1}{rgb}{0.00000,1.00000,1.00000}%
\definecolor{mycolor2}{rgb}{1.00000,0.00000,1.00000}%

\pgfplotsset{every axis label/.append style={font=\footnotesize},
every tick label/.append style={font=\footnotesize}
}

\begin{tikzpicture}[font=\footnotesize] 

\begin{axis}[%
name=ber,
width  = 0.5\columnwidth,
height = 0.3\columnwidth,
scale only axis,
xmin  = 1,
xmax  = 14,
xlabel= {$n(\times 1000)$},
xmajorgrids,
xtick       ={2, 4, 6, 8, 10, 12},
xticklabels ={$2$, $4$, $6$, $8$, $10$,$12$},
ymin = 0.0,
ymax = 7.2,
ylabel={Speedup},
ymajorgrids,
ytick       ={1, 2, 4, 6},
yticklabels ={$1$, $2$, $4$, $6$},
legend entries={$t_\text{rsvd}/t_\text{ruqlp}$-cpu,$t_\text{corutv}/t_\text{ruqlp}$-cpu, $t_\text{rptsod}/t_\text{ruqlp}$-cpu, $t_\text{rsvd}/t_\text{ruqlp}$-gpu, $t_\text{corutv}/t_\text{ruqlp}$-gpu, $t_\text{rptsod}/t_\text{ruqlp}$-gpu},
legend style={at={(1.58,1)},anchor=north east,draw=black,fill=white,legend cell align=left,font=\tiny},
]

\addplot[smooth,mark=square*,blue]
table[row sep=crcr]{
1	2.56110000000000   \\
2	3.55890000000000 \\
3	3.58310000000000 \\
4	3.08340000000000 \\
5	4.40280000000000 \\
6	6.04320000000000 \\
7	6.98720000000000 \\
8	6.57680000000000 \\
9	6.03360000000000 \\
10	6.80070000000000 \\
11	7.14740000000000 \\
12	7.10230000000000 \\
13	6.93280000000000 \\
14	6.86810000000000 \\
};

\addplot[smooth,mark=*,black]
table[row sep=crcr]{
 \\
1	1.76790000000000 \\
2	1.58950000000000 \\
3	1.48170000000000 \\
4	1.98230000000000 \\
5	2.14950000000000 \\
6	2.5297000000000 \\
7	2.87420000000000 \\
8	3.24690000000000 \\
9	3.49120000000000 \\
10	3.23450000000000 \\
11	3.42190000000000 \\
12	3.71850000000000 \\
13	3.56640000000000 \\
14	3.62320000000000 \\
};

\addplot[smooth,mark=diamond*,red]
table[row sep=crcr]{
1	2.17410000000000  \\
2	2.33110000000000 \\
3	2.43510000000000 \\
4	2.77020000000000 \\
5	3.69550000000000 \\
6	4.30310000000000 \\
7	4.99100000000000 \\
8	5.17810000000000 \\
9	5.21090000000000 \\
10	5.56740000000000 \\
11	5.97030000000000 \\
12	6.58250000000000 \\
13	6.21520000000000 \\
14	6.42530000000000 \\
};

\addplot[smooth,mark=pentagon*, teal]
table[row sep=crcr]{
1	2.33440000000000 \\
2	1.69100000000000 \\
3	1.76390000000000 \\
4	1.83440000000000 \\
5	1.94610000000000 \\
6	1.97440000000000 \\
7	1.94880000000000 \\
8	2.11180000000000 \\
9	2.16420000000000 \\
10	2.19340000000000 \\
11	2.23240000000000 \\
12	2.15750000000000 \\
13	2.29610000000000 \\
14	2.23120000000000 \\
};

\addplot[smooth,mark=asterisk, gray]
table[row sep=crcr]{
1	3.39500000000000   \\
2	2.15460000000000  \\
3	2.07670000000000  \\
4	2.01410000000000  \\
5	2.06560000000000  \\
6	2.04710000000000  \\
7	1.99720000000000  \\
8	2.12660000000000  \\
9	2.12780000000000  \\
10	2.12980000000000  \\
11	2.16810000000000  \\
12	2.07550000000000  \\
13	2.25260000000000  \\
14	2.15900000000000  \\
};

\addplot[smooth,mark=triangle*, orange]
table[row sep=crcr]{
1	0.960300000000000  \\
2	1.13450000000000 \\
3	1.26370000000000 \\ 
4	1.30030000000000 \\
5	1.35300000000000 \\
6	1.29280000000000 \\
7	1.30800000000000 \\
8	1.32640000000000 \\
9	1.35490000000000 \\
10	1.36531000000000 \\
11	1.38620000000000 \\
12	1.39610000000000 \\
13	1.50490000000000 \\
14	1.52920000000000 \\
};

\end{axis}

\end{tikzpicture}%
\captionsetup{justification=centering,font=scriptsize}  
\caption{Speedup for dense matrices with $d=0.3n$ with $q=0$.} 
\label{figSPUzof}      
\end{center}
\end{figure}


\subsection{Evaluation of Error Bounds}
To empirically assess the quality of derived bounds, we construct four $n\times n$ matrices, and consider one matrix from an application. The first two matrices are formed as:
 \begin{equation}\label{eq_A1PlusA2}
  {\bf A} = {\bf U}{\bf \Sigma}{\bf V}^T + \mu \sigma_k{\bf A}_N,
 \end{equation}
where $\bf U$ and $\bf V$ have random orthonormal columns, the entries of ${\bf \Sigma} = \text{diag}(\sigma_i)$ decrease linearly from 1 to $10^{-10}$, and $\sigma_{i+1}=0$ for $i\ge k$, and ${\bf A}_N$ is a normalized Gaussian matrix. Two cases for $\mu$ are considered: i) $\mu = 0.005$, giving the matrix \texttt{LowRankLargeGap}; and ii) $\mu = 0.01$, giving the matrix \texttt{LowRankMediumGap}. 

The second two matrices are formed as $ {\bf A} = {\bf U}{\bf \Sigma}{\bf V}^T$, with 
 \begin{equation}
{\bf \Sigma}= \text{diag}(\underbrace{
1, ..., 1}_k, 2^{-z}, 3^{-z}, ..., (n-k+1)^{-z}).
 \end{equation}

Two cases for $z$ are considered: i) $z = 1$, giving the matrix \texttt{LowRankSlowDecay}; and ii) $z=2$, giving the matrix \texttt{LowRankFastDecay}. These matrices were used in \cite{TrYUC17}. We set for the four matrices $n=800$, $k=16$, and $p=16$, hence $d=32$.

The fifth matrix \texttt{impcol\_e} is from a hydrocarbon separation problem taken from SuiteSparse Matrix Collection \cite{DavisH11}. This matrix is of order 225 and has a well-defined gap between $\sigma_{10}$ and $\sigma_{11}$.

The singular values of the matrices with their estimations by RU-QLP are displayed in Figures \ref{figSVMatSte}--\ref{figSVMatimpcole}. They show clearly the rank-revealing property and high accuracy of RU-QLP.

\subsubsection{Principal angles} We compute the sines of principal angles between subspaces, sin$\theta_i$ and sin$\phi_i$ Section \ref{subBoPABS}, and compare them to the theoretical upper bounds given in Theorem \ref{ThCanABk}. The results are shown in Figures \ref{figAnglQSte}--\ref{figAnglVTrop}. We make several observations: 
\begin{enumerate}
\item If the spectrum has a well-defined gap, the bounds are qualitatively accurate. In addition, the bounds become more accurate as the gap gets larger.

\item The bounds are quantitatively informative for polynomially decaying spectrums.  The principal angles become smaller as the spectrum decays faster.

\item By increasing $q$, the angles become smaller, and hence the approximate subspaces become more accurate.

\item $\theta_i$ are smaller than $\phi_i$ or, in other words, $\bf Q$ is a better approximation to ${\bf U}_k$ than $\bf P$ is to ${\bf V}_k$. This is due to the construction of these matrices by RU-QLP (see the explanation in Theorem \ref{ThCanABk})
\end{enumerate}

\subsubsection{Low-rank approximation} We compute in both spectral and Frobenius norms approximation errors incurred by RU-QLP's $\bf Q$ and $\bf P$, and compare them to the theoretical upper bounds given in Theorem \ref{ThLRAEB} as well as to minimum errors by the SVD. The results are displayed in Figures \ref{figQLRSteLarge}--\ref{figPLRimpcole}. We make three observations: 
\begin{enumerate}
\item When $q=0$, the bounds are qualitatively informative, particularly for the spectral norm cases.

\item When $q\ge 1$, the effect on the theoretical bounds is pronounced; they become qualitatively tight in most cases.

\item When $q=2$, the approximation errors closely match those of the optimal SVD in all cases, demonstrating the high quality of approximate left and right subspaces.
\end{enumerate}

\section{Conclusion}
\label{secConclusion}

We presented the rank-revealing RU-QLP decomposition, which furnishes an approximation to the SVD and pivoted QLP. We presented a new error analysis, which may be viewed as a systematic treatment of randomized low-rank matrix factorization methods. We furnished bounds in 2- and Frobenius norm on the rank-revealing property, principal angles between subspaces, and on the errors of low-rank approximations. We investigated the accuracy of the bounds using five matrices. We further investigated the runtime performance of RU-QLP and several existing methods on a hybrid GPU-accelerated multicore machine. Our results showed that RU-QLP harnesses best the advanced architecture, thus outperforming other randomized methods. 

\begin{figure}[t]
	\begin{center}
		\input{Plots/ASVSteMat}
	\captionsetup{justification=centering,font=scriptsize}
		\caption{Singular values of \texttt{LowRankLargeGap} (left), \texttt{LowRankMediumGap} (right)}
		\label{figSVMatSte}       
	\end{center}
\end{figure}
\begin{figure}[t]
	\begin{center}
		\input{Plots/ASVTropMat}
	\captionsetup{justification=centering,font=scriptsize}
		\caption{Singular values of \texttt{LowRankSlowDecay} (left), \texttt{LowRankFastDecay} (right).}
		\label{figSVMatTrop}       
	\end{center}
\end{figure}
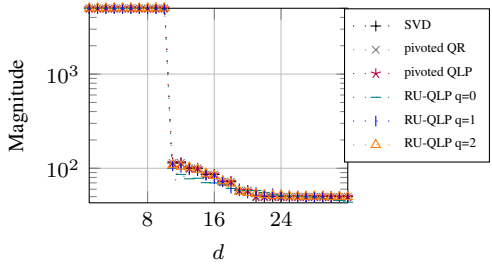
\begin{figure}[t]
\begin{center}     
%
%
%
\usetikzlibrary{positioning,calc}

\definecolor{mycolor1}{rgb}{0.00000,1.00000,1.00000}%
\definecolor{mycolor2}{rgb}{1.00000,0.00000,1.00000}%

\pgfplotsset{every axis label/.append style={font=\footnotesize},
every tick label/.append style={font=\footnotesize}
}

\begin{tikzpicture}[font=\footnotesize] 

\begin{axis}[%
name=ber,
ymode=log,
width  = 0.4\columnwidth,
height = 0.3\columnwidth,
scale only axis,
xmin  = 1,
xmax  = 32,
xlabel= {$d$},
xmajorgrids,
xtick       ={8, 16, 24},
xticklabels ={$8$, $16$, $24$},
ymin = 43,
ymax = 5.e+03,
ylabel={Magnitude},
ymajorgrids,
legend entries={SVD,pivoted QR,pivoted QLP, RU-QLP q=0, RU-QLP q=1, RU-QLP q=2},
legend style={at={(1.55,1)},anchor=north east,draw=black,fill=white,legend cell align=left,font=\tiny},
]

\addplot+[smooth,color=black,loosely dotted, every mark/.append style={solid}, mark=+]
table[row sep=crcr]{
1	5000.00032500008 \\
2	5000.00012500001 \\
3	5000.00010000001 \\
4	5000.00010000000 \\
5	5000.00010000000 \\
6	5000.00010000000 \\
7	5000.00010000000 \\
8	5000.00010000000 \\
9	5000.00010000000 \\
10	5000.00010000000 \\
11	114.013169754649 \\
12	114.008777902283 \\
13	100.015010421027 \\
14	100.010004991460 \\
15	86.0174518981040 \\
16	86.0116321768124 \\
17	72.0208436921943 \\
18	72.0138925058439 \\
19	58.0258740096420 \\
20	58.0172442883791 \\
21	50.0413761366775 \\
22	50.0325837540161 \\
23	50.0211955469256 \\
24	50.0125034542160 \\
25	50.0100097704246 \\
26	50.0100094813719 \\
27	50.0100094807702 \\
28	50.0100094803536 \\
29	50.0100094801164 \\ 
30	50.0100094800555 \\
31	50.0100094800543 \\
32	50.0100030018126 \\
};

\addplot+[smooth,color=gray,loosely dotted, every mark/.append style={solid}, mark=x]
table[row sep=crcr]{
1	5000 \\
2	5000 \\
3	5000 \\
4	5000 \\ 
5	5000 \\
6	5000 \\
7	5000 \\
8	5000 \\ 
9	5000 \\
10	5000 \\
11	114.008778459487 \\
12	114.004391853194 \\
13	100.010004886732 \\
14	100.005004866513 \\
15	86.0116310843995 \\
16	86.0058180264825 \\
17	72.0138896450997 \\
18	72.0069477250896 \\
19	58.0172379090825 \\
20	58.008622952832  \\
21	50.0099990002000 \\
22	50.0099950025986 \\
23	50 \\
24	50 \\
25	50 \\
26	50 \\
27	50 \\
28	50 \\ 
29	50 \\
30	50 \\
31	50 \\
32	50 \\
};

\addplot+[smooth,color=purple, loosely dotted, every mark/.append style={solid}, mark=star]
table[row sep=crcr]{
1	5000.00012500000 \\
2	5000.00012500000 \\
3	5000.00012500000 \\
4	5000.00012500000 \\
5	5000.00012500000 \\
6	5000.00012500000 \\
7	5000.00012500000 \\
8	5000.00012500000 \\
9	5000.00012500000 \\
10	5000.00012500000 \\
11	114.013167043205 \\
12	114.008777564613 \\
13	100.015007966817 \\
14	100.010004491110 \\
15	86.017448871388 \\ 
16	86.0116313899742 \\
17	72.0208392137189 \\
18	72.0138911644422 \\
19	58.0258661168358 \\
20	58.0172417204075 \\
21	50.0224979372287 \\
22	50.0224817077433 \\
23	50.0124984378905 \\
24	50.0124984378905 \\
25	50.0124981880779 \\
26	50.0124979391132 \\
27	50.0124979383152 \\
28	50.0124976894003 \\
29	50.0124976886024 \\
30	50.0124974397373 \\
31	50.0124974389396 \\
32	50.0124971901241 \\
};

\addplot+[smooth,color=teal,loosely dotted, every mark/.append style={solid}, mark=-]
table[row sep=crcr]{
1	4999.85728763061 \\
2	4999.86145123612\\
3	4999.81074678247\\
4	4999.83060628832\\
5	4999.82941950286\\
6	4999.79523415983\\
7	4999.78690385101\\
8	4999.73718827296\\
9	4999.77559909767\\
10	4793.18041435361\\
11	99.7678388666491\\
12	86.1448206261379\\
13	77.3487245340685\\
14	78.3810172234385\\
15	70.3627054216385\\
16	69.9880611941918\\
17	66.7446361945858\\
18	61.0816036162881\\
19	58.0464159985167\\
20	58.0526613003022\\
21	58.1270143326394\\
22	55.3782481181573\\
23	52.5807273193926\\
24	51.2211761499652\\
25	50.2575311375561\\
26	50.4058416141325\\
27	50.0461204007306\\
28	49.0889438713319\\
29	47.7470304672554\\
30	46.9097592661891\\
31	45.5194072262536\\
32	43.6326383324191 \\
};

\addplot+[smooth,color=blue,loosely dotted, every mark/.append style={solid}, mark=|]
table[row sep=crcr]{
1	5000.00012412011  \\
2	5000.00012512174  \\
3	5000.0001229937  \\
4	5000.00012976889  \\
5	5000.00012447146  \\
6	5000.00011859382  \\
7	5000.00012952545  \\
8	5000.00012032255  \\
9	5000.00011119723  \\
10	5000.00014371622  \\
11	106.908374191917  \\
12	102.497004751313  \\
13	96.5530563194041  \\
14	94.4269984889860  \\
15	84.5461141673070  \\
16	77.7688873616808  \\
17	73.3270722669887  \\
18	61.6020472375821  \\
19	60.4223832230741  \\
20	54.6101907897730  \\
21	55.2021016009807  \\
22	52.339625072054  \\
23	53.5894039917060  \\
24	51.9018367165120  \\
25	50.885738722385  \\
26	50.5191877494430  \\
27	49.7933760339270  \\
28	49.9715888367394  \\
29	50.4596719831712  \\
30	49.6123566696670  \\
31	49.8262579087530  \\
32	49.4727531113512  \\
};

\addplot+[smooth,color=orange,loosely dotted, every mark/.append style={solid}, mark=triangle]
table[row sep=crcr]{
1	5000.00011612006  \\
2	5000.00011789077 \\
3	5000.00012176322 \\
4	5000.00013276354 \\
5	5000.00014029104 \\
6	5000.00012999359 \\
7	5000.00012113191 \\
8	5000.00012861364 \\
9	5000.00011684973 \\
10	5000.00012458261 \\
11	108.265557562993 \\
12	104.296766229538 \\
13	102.376947115600 \\
14	96.8794894162268 \\
15	87.7701079410757 \\
16	82.9734519766162 \\
17	72.994649037806 \\
18	69.0778407717175 \\
19	56.8631790796786 \\
20	54.9709617245314 \\
21	53.6140687499118 \\
22	52.025484709294 \\
23	51.1898542425641 \\
24	51.6237798866328 \\
25	50.575392019008 \\
26	50.2912517681965 \\
27	50.1204813674829 \\
28	50.3715765539698 \\
29	50.3098589631210 \\
30	49.642809646504 \\
31	50.0531110904213 \\
32	49.7408856159650 \\
};

\end{axis}

\end{tikzpicture}%
\captionsetup{justification=centering,font=scriptsize}
\caption{Singular values of \texttt{impcol\_e}}     
\label{figSVMatimpcole}     
\end{center}  
\end{figure} 

\begin{figure}[t]
\begin{center}       
%
%
%
\usetikzlibrary{positioning,calc}

\definecolor{mycolor1}{rgb}{0.00000,1.00000,1.00000}%
\definecolor{mycolor2}{rgb}{1.00000,0.00000,1.00000}%

\pgfplotsset{every axis label/.append style={font=\footnotesize},
every tick label/.append style={font=\footnotesize}
}

\begin{tikzpicture}[font=\footnotesize] 

\begin{axis}[%
name=ber,
ymode=log,
width  = 0.34\columnwidth,
height = 0.3\columnwidth,
scale only axis,
xmin  = 1,
xmax  = 16,
xlabel= {Index},
xmajorgrids,
ymin = 0.0,
ymax = 0.001,
xtick       ={5, 10},
xticklabels ={$5$, $10$},
ylabel={sin$\theta_i$},
ymajorgrids,
]
\addplot+[smooth,color=green,solid, every mark/.append style={solid}, mark=none]
table[row sep=crcr]{
1	0.000306201692778056  \\
2	0.000307110549006419 \\
3	0.000307747498952399 \\
4	0.000308572133472086 \\
5	0.000309234149093798 \\
6	0.000310238580229057 \\
7	0.000310894140482281 \\
8	0.000311693768784512 \\
9	0.000312612010570418 \\
10	0.000313306843156637 \\
11	0.000314081901949417 \\
12	0.000314902057465304 \\
13	0.000315496337481174 \\
14	0.000316537517696999 \\
15	0.000317369714768109 \\
16	0.000318101093843612 \\
};

\addplot+[smooth,color=green, dotted, every mark/.append style={solid}, mark=|]
table[row sep=crcr]{
1	6.81470784360175e-09  \\
2	6.85522219860501e-09 \\
3	6.88368726988044e-09 \\
4	6.92062750213845e-09 \\
5	6.95035461058144e-09 \\
6	6.99557918292032e-09 \\
7	7.02517491809997e-09 \\
8	7.06135928265395e-09 \\
9	7.10302573443145e-09 \\
10	7.13463615854232e-09 \\
11	7.16997916442391e-09 \\
12	7.20747369256227e-09 \\
13	7.23470310923936e-09 \\
14	7.28253288694770e-09 \\
15	7.32087570277325e-09 \\
16	7.35465652768496e-09 \\
};

\addplot+[smooth,color=green, loosely dashed, every mark/.append style={solid}, mark=none]
table[row sep=crcr]{
1	1.62487033362823e-13  \\
2	1.63938195787142e-13 \\
3	1.64960340618698e-13 \\
4	1.66289971297550e-13 \\
5	1.67362552788594e-13 \\
6	1.68998699282790e-13 \\
7	1.70072290244130e-13 \\
8	1.71387960377477e-13 \\
9	1.72907140718190e-13 \\
10	1.74062650641220e-13 \\
11	1.75357638129220e-13 \\
12	1.76734951785692e-13 \\
13	1.77737438695096e-13 \\
14	1.79502924551245e-13 \\ 
15	1.80922423597859e-13 \\
16	1.82176116575342e-13 \\
};

\addplot+[smooth,color=blue,densely dashed, every mark/.append style={solid}, mark=none]
table[row sep=crcr]{
1	0.000112805466781591  \\
2	9.05967038021850e-05  \\
3	7.54794140689154e-05  \\
4	7.11218087906559e-05  \\
5	5.73125698514657e-05  \\
6	5.61110230100237e-05  \\
7	4.67796103324712e-05  \\
8	4.38106079500859e-05  \\
9	4.31401611640833e-05  \\
10	3.92806418578699e-05  \\
11	3.82603686157752e-05  \\
12	3.65070028232955e-05  \\
13	3.10904921262284e-05  \\
14	2.86186508345777e-05  \\
15	2.58528863313819e-05  \\
16	2.22764159490343e-05  \\
};

\addplot+[smooth,color=blue,densely dotted, every mark/.append style={solid}, mark=none]
table[row sep=crcr]{
1	1.55326165009224e-09  \\
2	1.35837316232742e-09 \\
3	1.27782782308245e-09 \\
4	9.60846667147340e-10 \\
5	8.63907390449709e-10 \\
6	8.16769204162813e-10 \\
7	6.99289918779351e-10 \\
8	6.46132329810481e-10 \\
9	5.57332442967314e-10 \\
10	5.05088686255812e-10 \\
11	4.59844483914494e-10 \\
12	4.40626816804476e-10 \\
13	4.09137591059872e-10 \\
14	3.92924002046674e-10 \\
15	3.42844572519874e-10 \\
16	2.82015171978676e-10 \\
};

\addplot+[smooth,color=blue,loosely dotted, every mark/.append style={solid}, mark=-]
table[row sep=crcr]{
1	2.20951796356791e-14  \\
2	2.03428031568097e-14 \\
3	1.99131089472608e-14 \\
4	1.79194518858946e-14 \\
5	1.61191000882652e-14 \\
6	1.26154219573265e-14 \\
7	1.18606234615548e-14 \\
8	1.11214407663117e-14 \\
9	1.00231139509196e-14 \\
10	8.87896974453336e-15 \\
11	8.49176509834063e-15 \\
12	7.77271856866792e-15 \\
13	7.29625214103300e-15 \\
14	6.58678118872868e-15 \\
15	6.01411179508135e-15 \\
16	5.27698395537874e-15 \\
};

\end{axis}

\begin{axis}[%
name=SumRate,
at={($(ber.east)+(35,0em)$)},
		anchor= west,
ymode=log,
width  = 0.34\columnwidth,
height = 0.3\columnwidth,
scale only axis,
xmin   = 1,
xmax  = 16,
xlabel= {Index},
xmajorgrids,
ymin = 0.0,
ymax = 0.01,
xtick       ={5, 10},
xticklabels ={$5$, $10$},
ymajorgrids,
legend entries={Th. bound q=0, Th. bound q=1,  Th. bound q=2,  Computed q=0, Computed q=1, Computed q=2}, 
legend style={at={(0.7,1.3)},anchor=north east,draw=black,fill=white,legend cell align=left,font=\tiny, legend columns=3}
]

\addplot+[smooth,color=green,solid, every mark/.append style={thick}, mark=none]
table[row sep=crcr]{
1	0.00110529577114235  \\
2	0.00110920908655944 \\
3	0.00111139193121813 \\
4	0.00111463680994620 \\
5	0.00111738405441990 \\
6	0.00111990921842394 \\
7	0.00112320499246392 \\
8	0.00112604092921753 \\ 
9	0.00112879336122020 \\
10	0.00113010068999228 \\
11	0.00113356739798721 \\
12	0.00113685785737616 \\
13	0.00113943501877848 \\
14	0.00114328243735472 \\
15	0.00114546180575322 \\
16	0.00114954529841195 \\
};

\addplot+[smooth,color=green, dotted, every mark/.append style={solid}, mark=|]
table[row sep=crcr]{
1	1.23163695869076e-07  \\
2	1.24037366468752e-07 \\
3	1.24526040882866e-07 \\
4	1.25254248938731e-07 \\
5	1.25872438696637e-07 \\
6	1.26441997749082e-07 \\
7	1.27187304531143e-07 \\
8	1.27830376747645e-07 \\
9	1.28456064276555e-07 \\ 
10	1.28753783552647e-07 \\
11	1.29544928918279e-07 \\
12	1.30298093742851e-07 \\
13	1.30889513798762e-07 \\
14	1.31774930876298e-07 \\
15	1.32277799127203e-07 \\
16	1.33222603995393e-07 \\
};

\addplot+[smooth,color=green, loosely dashed, every mark/.append style={solid}, mark=none]
table[row sep=crcr]{
1	1.07819742940987e-11  \\
2	1.08969015787763e-11 \\
3	1.09613613280872e-11 \\
4	1.10576521496503e-11 \\
5	1.11396152946448e-11 \\
6	1.12153090759146e-11 \\
7	1.13146172769559e-11 \\
8	1.14005374735727e-11 \\
9	1.14843426202187e-11 \\
10	1.15242911877109e-11 \\
11	1.16306730830667e-11 \\
12	1.17322502020511e-11 \\
13	1.18122193968542e-11 \\
14	1.19322795290925e-11 \\
15	1.20006470711617e-11 \\
16	1.21294497720563e-11 \\
};

\addplot+[smooth,color=blue,densely dashed, every mark/.append style={solid}, mark=none]
table[row sep=crcr]{
1	0.000382221376974906   \\
2	0.000346935702104776  \\
3	0.000339248961871418  \\
4	0.000226722245236292  \\
5	0.000219190028316944  \\
6	0.000203997104844051  \\
7	0.000184969317968262  \\
8	0.000172068277738557  \\
9	0.000164211450406647  \\
10	0.000150094419111759  \\
11	0.000142122006950635  \\
12	0.000139312510703960  \\
13	0.000128989159345413  \\
14	0.000113819089115866  \\  
15	0.000104022537527378  \\
16	0.000100679201239130  \\
};

\addplot+[smooth,color=blue,densely dotted, every mark/.append style={solid}, mark=none]
table[row sep=crcr]{
1	2.84237139305243e-08  \\
2	2.54661434414976e-08 \\
3	2.14136343003489e-08 \\
4	1.67042934692263e-08 \\
5	1.38123153840570e-08 \\
6	1.21881159220996e-08 \\
7	1.12774585464041e-08 \\
8	9.82362877189974e-09 \\
9	9.53196519669044e-09 \\
10	8.81669721762740e-09 \\
11	8.38606178118513e-09 \\
12	7.58588174030224e-09 \\
13	6.96108074165756e-09 \\
14	6.53423839456359e-09 \\
15	5.41078641144113e-09 \\
16	5.03791022069499e-09 \\
};

\addplot+[smooth,color=blue,loosely dotted, every mark/.append style={solid}, mark=-]
table[row sep=crcr]{
1	2.08851397941688e-12  \\
2	1.81789745838962e-12 \\
3	1.69426909326623e-12 \\
4	1.33451930726355e-12 \\
5	1.07328424661128e-12 \\
6	9.45733577190082e-13 \\
7	8.40671553367334e-13 \\
8	7.23172351654413e-13 \\
9	6.58513925946306e-13 \\
10	6.03828056386166e-13 \\
11	5.52804264322720e-13 \\
12	5.01678813441479e-13 \\
13	4.71204409310682e-13 \\
14	4.55934986459685e-13 \\
15	3.89079409832542e-13 \\
16	3.40584926328537e-13 \\
};

\end{axis}

\end{tikzpicture}%
\captionsetup{justification=centering,font=scriptsize}  
\caption{Principal angles between $\mathcal{R}({\bf Q})$ and $\mathcal{R}({\bf U}_k)$ for  \texttt{LowRankLargeGap} (left), and \texttt{LowRankMediumGap} (right).} 
\label{figAnglQSte}      
\end{center}
\end{figure}
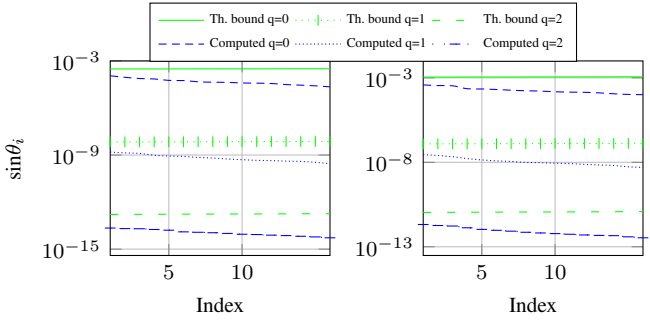
\begin{figure}[t]
\begin{center}        
%
%
%
\usetikzlibrary{positioning,calc}

\definecolor{mycolor1}{rgb}{0.00000,1.00000,1.00000}%
\definecolor{mycolor2}{rgb}{1.00000,0.00000,1.00000}%

\pgfplotsset{every axis label/.append style={font=\footnotesize},
every tick label/.append style={font=\footnotesize}
}

\begin{tikzpicture}[font=\footnotesize] 

\begin{axis}[%
name=ber,
ymode=log,
width  = 0.34\columnwidth,
height = 0.3\columnwidth,
scale only axis,
xmin  = 1,
xmax  = 16,
xlabel= {Index},
xmajorgrids,
ymin = 0.0,
ymax = 1,
xtick       ={5, 10},
xticklabels ={$5$, $10$},
ylabel={sin$\theta_i$},
ymajorgrids,
]

\addplot+[smooth,color=green,solid, every mark/.append style={solid}, mark=none]
table[row sep=crcr]{
1	0.959260556618083  \\
2	0.959260556618083 \\
3	0.959260556618083 \\
4	0.959260556618083 \\
5	0.959260556618083 \\
6	0.959260556618083 \\
7	0.959260556618083 \\
8	0.959260556618083 \\
9	0.959260556618083 \\
10	0.959260556618083 \\
11	0.959260556618083 \\
12	0.959260556618083 \\ 
13	0.959260556618083 \\
14	0.959260556618083 \\
15	0.959260556618083 \\
16	0.959260556618083 \\
};

\addplot+[smooth,color=green, dotted, every mark/.append style={solid}, mark=|]
table[row sep=crcr]{
1	0.629184090006384  \\
2	0.629184090006385 \\
3	0.629184090006385 \\
4	0.629184090006385 \\
5	0.629184090006385 \\
6	0.629184090006385 \\
7	0.629184090006385 \\
8	0.629184090006386 \\
9	0.629184090006386 \\
10	0.629184090006387 \\
11	0.629184090006387 \\
12	0.629184090006387 \\
13	0.629184090006387 \\
14	0.629184090006387 \\
15	0.629184090006387 \\
16	0.629184090006387 \\
};

\addplot+[smooth,color=green, loosely dashed, every mark/.append style={solid}, mark=none]
table[row sep=crcr]{
1	0.215359007786829  \\
2	0.215359007786830 \\
3	0.215359007786830 \\
4	0.215359007786830 \\
5	0.215359007786830 \\
6	0.215359007786831 \\
7	0.215359007786831 \\
8	0.215359007786831 \\
9	0.215359007786831 \\
10	0.215359007786831 \\
11	0.215359007786831 \\
12	0.215359007786832 \\
13	0.215359007786832 \\
14	0.215359007786832 \\
15	0.215359007786832 \\
16	0.215359007786832 \\
};

\addplot+[smooth,color=blue,densely dashed, every mark/.append style={solid}, mark=none]
table[row sep=crcr]{
1	0.0127390446619115     \\
2	0.00688951165229769  \\
3	0.00654494397147856  \\
4	0.00395643544111637  \\
5	0.00315332406322982  \\
6	0.00295989132394129  \\
7	0.00236972101403484  \\
8	0.00181175314176161  \\
9	0.00168424136794331  \\
10	0.00134826585742103  \\
11	0.00112605956325564  \\
12	0.000886498147799171  \\
13	0.00082796277647884  \\
14	0.000617286004960489  \\
15	0.00057373040852750  \\
16	0.000534939196017148  \\
};

\addplot+[smooth,color=blue,densely dotted, every mark/.append style={solid}, mark=none]
table[row sep=crcr]{
1	3.60639545602986e-05   \\
2	1.66412112351073e-05  \\
3	1.31417817666395e-05  \\
4	1.01243910100312e-05  \\
5	6.42263930129291e-06  \\
6	4.31343940364975e-06  \\
7	3.52620759186033e-06  \\
8	2.50901155808060e-06  \\
9	1.39656540102354e-06  \\
10	1.21233919749366e-06  \\
11	1.08245801574808e-06  \\
12	7.32094809580126e-07  \\
13	6.98188650702384e-07  \\
14	5.95388980263785e-07  \\
15	4.43872085633153e-07  \\
16	3.14548438941405e-07  \\
};

\addplot+[smooth,color=blue,loosely dotted, every mark/.append style={solid}, mark=-]
table[row sep=crcr]{
1	1.23394917320304e-07  \\
2	5.09515783688814e-08 \\
3	3.22734264367093e-08 \\
4	1.76125377139816e-08 \\
5	1.18246515968802e-08 \\
6	7.75315439426797e-09 \\
7	5.53412067340126e-09 \\
8	4.09218558706767e-09 \\
9	3.22890952913717e-09 \\
10	2.44058920855690e-09 \\
11	1.24774897272222e-09 \\
12	9.60300992355948e-10 \\
13	6.00791284502430e-10 \\
14	4.63865288004610e-10 \\
15	3.69679212419444e-10 \\
16	2.47031953623827e-10 \\
};

\end{axis}

\begin{axis}[%
name=SumRate,
at={($(ber.east)+(35,0em)$)},
		anchor= west,
ymode=log,
width  = 0.34\columnwidth,
height = 0.3\columnwidth,
scale only axis,
xmin   = 1,
xmax  = 16,
xlabel= {Index},
xmajorgrids,
ymin = 0.0,
ymax = 0.7,
xtick       ={5, 10},
xticklabels ={$5$, $10$},
ymajorgrids,
]

\addplot+[smooth,color=green,solid, every mark/.append style={solid}, mark=none]
table[row sep=crcr]{
1	0.618510910037144  \\
2	0.618510910037145  \\
3	0.618510910037145  \\
4	0.618510910037145  \\
5	0.618510910037145  \\
6	0.618510910037145  \\
7	0.618510910037145  \\
8	0.618510910037145  \\
9	0.618510910037145  \\
10	0.618510910037145  \\
11	0.618510910037145  \\
12	0.618510910037145  \\
13	0.618510910037145  \\
14	0.618510910037145  \\
15	0.618510910037145  \\
16	0.618510910037145 \\
};

\addplot+[smooth,color=green, dotted, every mark/.append style={solid}, mark=|]
table[row sep=crcr]{
1	0.0482969076210420   \\
2	0.0482969076210422  \\
3	0.0482969076210422  \\
4	0.0482969076210422  \\
5	0.0482969076210422  \\
6	0.0482969076210422  \\
7	0.0482969076210422  \\
8	0.0482969076210422  \\
9	0.0482969076210423  \\
10	0.0482969076210423  \\
11	0.0482969076210423  \\
12	0.0482969076210423  \\
13	0.0482969076210423  \\
14	0.0482969076210423  \\
15	0.0482969076210423  \\
16	0.0482969076210424 \\
};

\addplot+[smooth,color=green, loosely dashed, every mark/.append style={solid}, mark=none]
table[row sep=crcr]{
1	0.00298753573987305  \\
2	0.00298753573987307  \\
3	0.00298753573987307  \\
4	0.00298753573987307  \\
5	0.00298753573987307  \\
6	0.00298753573987307  \\
7	0.00298753573987307  \\
8	0.00298753573987307  \\
9	0.00298753573987308  \\
10	0.00298753573987308  \\
11	0.00298753573987308  \\
12	0.00298753573987308  \\
13	0.00298753573987308  \\
14	0.00298753573987308  \\
15	0.00298753573987308  \\
16	0.00298753573987309  \\
};

\addplot+[smooth,color=blue,densely dashed, every mark/.append style={solid}, mark=none]
table[row sep=crcr]{
1	2.52226286932301e-05   \\
2	1.66599999894212e-05  \\
3	1.44988291738515e-05  \\
4	1.07733442932559e-05  \\
5	8.39430248868299e-06  \\
6	6.36752366946986e-06  \\
7	4.03052205463384e-06  \\
8	3.70561105744759e-06  \\
9	2.70395894573033e-06  \\
10	2.58261726187410e-06  \\
11	1.67828264320503e-06  \\
12	1.29259596875147e-06  \\
13	1.13841321539607e-06  \\
14	8.69920911012297e-07  \\
15	4.19896593278258e-07  \\
16	3.15232342220929e-07  \\
};

\addplot+[smooth,color=blue,densely dotted, every mark/.append style={solid}, mark=none]
table[row sep=crcr]{
1	6.63106602060602e-10  \\
2	1.86366333699631e-10  \\
3	7.23615213345091e-11  \\
4	4.90704181021862e-11  \\
5	3.76711751282161e-11  \\
6	2.54759415999490e-11  \\
7	1.12901175644565e-11  \\
8	9.31259920568175e-12  \\
9	5.27092879656303e-12  \\
10	3.86976129661483e-12  \\
11	2.27286792244674e-12  \\
12	1.75093516282329e-12  \\
13	6.35479085733978e-13  \\
14	3.69012757703135e-13  \\
15	2.16992775608385e-13  \\
16	9.98417121342052e-14 \\
};

\addplot+[smooth,color=blue,loosely dotted, every mark/.append style={solid}, mark=-]
table[row sep=crcr]{
1	7.55425711480307e-15  \\
2	2.95851162053857e-15  \\
3	1.42279020530679e-15  \\
4	1.36002313041919e-15  \\
5	9.24257614836360e-16  \\
6	8.43175398941687e-16  \\
7	8.08758117709373e-16  \\
8	7.69638466005760e-16  \\
9	7.42698741067945e-16  \\
10	7.13125095861965e-16  \\
11	6.80217185136679e-16  \\
12	6.60198499294612e-16  \\
13	6.55392493034180e-16  \\
14	6.29540443843959e-16  \\
15	5.92337006895775e-16  \\
16	5.85037648397538e-16  \\
};

\end{axis}

\end{tikzpicture}%
\captionsetup{justification=centering,font=scriptsize}  
\caption{Principal angles between $\mathcal{R}({\bf Q})$ and $\mathcal{R}({\bf U}_k)$ for  \texttt{LowRankSlowDecay} (left), and \texttt{LowRankFastDecay} (right).} 
\label{figAnglQTrop}      
\end{center}
\end{figure}
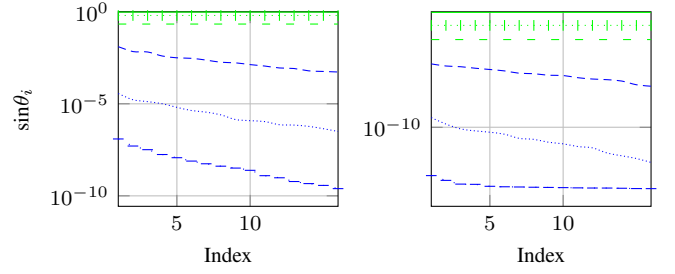
\begin{figure}[t]
\begin{center}       
%
%
%
\usetikzlibrary{positioning,calc}

\definecolor{mycolor1}{rgb}{0.00000,1.00000,1.00000}%
\definecolor{mycolor2}{rgb}{1.00000,0.00000,1.00000}%

\pgfplotsset{every axis label/.append style={font=\footnotesize},
every tick label/.append style={font=\footnotesize}
}

\begin{tikzpicture}[font=\footnotesize] 

\begin{axis}[%
name=ber,
ymode=log,
width  = 0.34\columnwidth,
height = 0.3\columnwidth,
scale only axis,
xmin  = 1,
xmax  = 10,
xlabel= {Index},
xmajorgrids,
ymin = 0.0,
ymax = 0.01,
xtick       ={2, 4, 6, 8},
xticklabels ={$2$,$4$,$6$,$8$},
ylabel={sin$\theta_i$/sin$\phi_i$},
ymajorgrids,
]

\addplot+[smooth,color=green, solid, every mark/.append style={solid}, mark=none]
table[row sep=crcr]{
1	0.00480939028041113   \\
2	0.00480939066515358  \\
3	0.00480939071324638  \\
4	0.00480939071324639  \\
5	0.00480939071324639  \\
6	0.00480939071324639  \\
7	0.00480939071324639  \\
8	0.00480939071324639  \\
9	0.00480939071324639  \\
10	0.00480939071324639 \\
};

\addplot+[smooth,color=green, dotted, every mark/.append style={solid}, mark=|]
table[row sep=crcr]{
1	2.62873250469841e-06  \\
2	2.62873292529576e-06  \\
3	2.62873297787043e-06  \\
4	2.62873297787043e-06  \\
5	2.62873297787044e-06  \\
6	2.62873297787044e-06  \\
7	2.62873297787044e-06  \\
8	2.62873297787044e-06  \\
9	2.62873297787044e-06  \\
10	2.62873297787044e-06 \\
};

\addplot+[smooth,color=green, loosely dashed, every mark/.append style={solid}, mark=none]
table[row sep=crcr]{
1	1.33562087547431e-09  \\
2	1.33562119602345e-09  \\
3	1.33562123609209e-09  \\
4	1.33562123609210e-09  \\
5	1.33562123609210e-09  \\
6	1.33562123609210e-09  \\
7	1.33562123609210e-09  \\
8	1.33562123609210e-09  \\
9	1.33562123609210e-09  \\
10	1.33562123609210e-09  \\
};

\addplot+[smooth,color=blue,densely dashed, every mark/.append style={solid}, mark=none]
table[row sep=crcr]{
1	0.000507303863724995  \\
2	0.000388100447525928  \\
3	0.000260355617851734  \\
4	0.000223690164246055  \\
5	0.000169099164404334  \\
6	0.000138281655782456  \\
7	9.37471450426147e-05  \\
8	8.80015897352973e-05  \\
9	6.33232884533455e-05  \\
10	3.62998760379301e-05  \\
};

\addplot+[smooth,color=blue,densely dotted, every mark/.append style={solid}, mark=none]
table[row sep=crcr]{
1	8.82097205637245e-08   \\
2	5.04866437220056e-08  \\
3	2.70578295096308e-08  \\
4	2.34211327279771e-08  \\
5	1.62408109855647e-08  \\
6	1.37660490855917e-08  \\
7	1.03656071486616e-08  \\
8	8.14495912212779e-09  \\
9	5.23371987621152e-09  \\
10	3.83031448621008e-09  \\
};

\addplot+[smooth,color=blue,loosely dotted, every mark/.append style={solid}, mark=-]
table[row sep=crcr]{
1	7.81071750774594e-12   \\
2	4.42075523594691e-12  \\
3	3.98951197163496e-12  \\
4	3.02577528598522e-12  \\
5	2.26369127999466e-12  \\
6	1.55775371115325e-12  \\
7	1.07783570260459e-12  \\
8	8.47051574460522e-13  \\  
9	6.86176411333840e-13  \\
10	4.24394518848567e-13  \\
};

\end{axis}

\begin{axis}[%
name=SumRate,
at={($(ber.east)+(35,0em)$)},
		anchor= west,
ymode=log,
width  = 0.34\columnwidth,
height = 0.3\columnwidth,
scale only axis,
xmin   = 1,
xmax  = 10,
xlabel= {Index},
xmajorgrids,
ymin = 0.0,
ymax = 0.22,
xtick       ={2, 4, 6, 8},
xticklabels ={$2$,$4$,$6$,$8$},
ymajorgrids,
]

\addplot+[smooth,color=green,solid, every mark/.append style={solid}, mark=none]
table[row sep=crcr]{
1	0.216924747387137   \\
2	0.216924755655822  \\
3	0.216924756689408  \\
4	0.216924756689408  \\
5	0.216924756689408  \\
6	0.216924756689408  \\
7	0.216924756689408  \\
8	0.216924756689408  \\  
9	0.216924756689408  \\
10	0.216924756689408 \\
};

\addplot+[smooth,color=green, dotted, every mark/.append style={solid}, mark=|]
table[row sep=crcr]{
1	0.000109664337074891   \\
2	0.000109664350234616  \\
3	0.000109664351879582  \\
4	0.000109664351879582  \\
5	0.000109664351879582  \\
6	0.000109664351879582  \\
7	0.000109664351879582  \\
8	0.000109664351879582  \\
9	0.000109664351879582  \\
10	0.000109664351879582 \\
};

\addplot+[smooth,color=green, loosely dashed, every mark/.append style={solid}, mark=none]
table[row sep=crcr]{
1	6.48656509910222e-08   \\
2	6.48656639641575e-08  \\
3	6.48656655857994e-08  \\
4	6.48656655857996e-08  \\
5	6.48656655857996e-08  \\
6	6.48656655857996e-08  \\
7	6.48656655857997e-08  \\
8	6.48656655857997e-08  \\
9	6.48656655857997e-08  \\
10	6.48656655857998e-08  \\
};

\addplot+[smooth,color=blue,densely dashed, every mark/.append style={solid}, mark=none]
table[row sep=crcr]{
1	0.0427852007408759   \\
2	0.0224113117590083  \\
3	0.0205717747989818  \\
4	0.0184601692111982  \\
5	0.0114722922632097  \\
6	0.0100760928022697  \\
7	0.00890296393100922  \\
8	0.00867797141459245  \\
9	0.00702186023717126  \\
10	0.00599283807094799  \\
};

\addplot+[smooth,color=blue,densely dotted, every mark/.append style={solid}, mark=none]
table[row sep=crcr]{
1	5.60215764914524e-06   \\
2	3.96033268980407e-06  \\
3	2.57114975228928e-06  \\
4	1.70294624895916e-06  \\
5	1.52099850980358e-06  \\
6	1.30017458293844e-06  \\
7	1.22514197940268e-06  \\
8	7.96964484061435e-07  \\
9	6.35410153805043e-07  \\
10	6.05623397001690e-07  \\
};

\addplot+[smooth,color=blue,loosely dotted, every mark/.append style={solid}, mark=-]
table[row sep=crcr]{
1	8.58427737681405e-10   \\
2	4.35680587685420e-10  \\
3	3.26393427774063e-10  \\
4	2.18513130448331e-10  \\
5	1.73386594494103e-10  \\
6	1.48149691716408e-10  \\
7	1.29060455262526e-10  \\
8	9.95852037675281e-11  \\
9	6.76495999476791e-11  \\
10	4.34152384934031e-11  \\
};

\end{axis}

\end{tikzpicture}%
\captionsetup{justification=centering,font=scriptsize}   
\caption{Principal angles between $\mathcal{R}({\bf Q})$ and $\mathcal{R}({\bf U}_k)$ (left) and between $\mathcal{R}({\bf P})$ and $\mathcal{R}({\bf V}_k)$ (right) for \texttt{impcol\_e}.} 
\label{figAnglQimpcole}       
\end{center}
\end{figure}
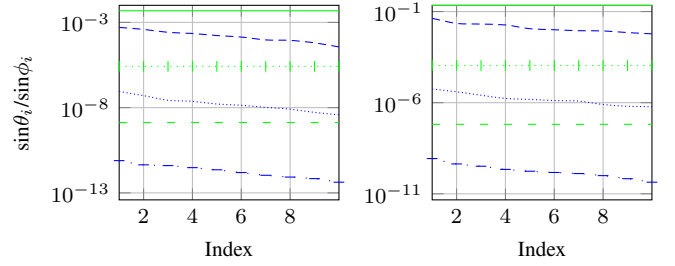
\begin{figure}[t]
\begin{center}       
%
%
%
\usetikzlibrary{positioning,calc}

\definecolor{mycolor1}{rgb}{0.00000,1.00000,1.00000}%
\definecolor{mycolor2}{rgb}{1.00000,0.00000,1.00000}%

\pgfplotsset{every axis label/.append style={font=\footnotesize},
every tick label/.append style={font=\footnotesize}
}

\begin{tikzpicture}[font=\footnotesize] 

\begin{axis}[%
name=ber,
ymode=log,
width  = 0.34\columnwidth,
height = 0.3\columnwidth,
scale only axis,
xmin  = 1,
xmax  = 16,
xlabel= {Index},
xmajorgrids,
ymin = 0.0,
ymax = 0.1,
xtick       ={5, 10},
xticklabels ={$5$, $10$},
ylabel={sin$\phi_i$},
ymajorgrids,
]
\addplot+[smooth,color=green,solid, every mark/.append style={solid}, mark=none]
table[row sep=crcr]{
1	0.0767853533961513   \\
2	0.0768881605134613  \\
3	0.0769602287005758  \\
4	0.0770642840489565  \\
5	0.0771811483742199  \\
6	0.0772517773953626  \\
7	0.0773469414843904  \\
8	0.0774382486238043  \\
9	0.0775547378638385  \\
10	0.0776558011363944  \\
11	0.0777335474657682  \\
12	0.0778427673737887  \\
13	0.0779572697072519  \\
14	0.0780472971019224  \\
15	0.0781635523940218  \\
16	0.0782434020502158  \\
};

\addplot+[smooth,color=green, dotted, every mark/.append style={solid}, mark=|]
table[row sep=crcr]{
1	1.80568609405433e-06  \\
2	1.81299187168547e-06 \\
3	1.81812508576227e-06 \\
4	1.82555388821149e-06 \\
5	1.83392146891403e-06 \\
6	1.83899104985703e-06 \\
7	1.84583659389299e-06 \\
8	1.85242078439573e-06 \\
9	1.86084377581083e-06 \\
10	1.86817220053145e-06 \\
11	1.87382303239291e-06 \\
12	1.88178086282176e-06 \\
13	1.89014795196562e-06 \\
14	1.89674411885253e-06 \\
15	1.90528485041952e-06 \\
16	1.91116599449705e-06 \\
};

\addplot+[smooth,color=green, loosely dashed, every mark/.append style={solid}, mark=none]
table[row sep=crcr]{
1	3.00885044965636e-11  \\
2	3.02916740409905e-11 \\
3	3.04347527793311e-11 \\
4	3.06422940219649e-11 \\
5	3.08767374029974e-11 \\
6	3.10191247351628e-11 \\
7	3.12118084295344e-11 \\
8	3.13975857484472e-11 \\
9	3.16358885434414e-11 \\
10	3.18438098508353e-11 \\ 
11	3.20045064802173e-11 \\
12	3.22313571476502e-11 \\
13	3.24705650785892e-11 \\
14	3.26596422675670e-11 \\
15	3.29051117675441e-11 \\
16	3.30745691505898e-11 \\
};

\addplot+[smooth,color=blue,densely dashed, every mark/.append style={solid}, mark=none]
table[row sep=crcr]{
1	0.000112805466781591  \\
2	9.05967038021850e-05  \\
3	7.54794140689154e-05  \\
4	7.11218087906559e-05  \\
5	5.73125698514657e-05  \\
6	5.61110230100237e-05  \\
7	4.67796103324712e-05  \\
8	4.38106079500859e-05  \\
9	4.31401611640833e-05  \\
10	3.92806418578699e-05  \\
11	3.82603686157752e-05  \\
12	3.65070028232955e-05  \\
13	3.10904921262284e-05  \\
14	2.86186508345777e-05  \\
15	2.58528863313819e-05  \\
16	2.22764159490343e-05  \\
};

\addplot+[smooth,color=blue,densely dotted, every mark/.append style={solid}, mark=none]
table[row sep=crcr]{
1	1.55326165009224e-09  \\
2	1.35837316232742e-09 \\
3	1.27782782308245e-09 \\
4	9.60846667147340e-10 \\
5	8.63907390449709e-10 \\
6	8.16769204162813e-10 \\
7	6.99289918779351e-10 \\
8	6.46132329810481e-10 \\
9	5.57332442967314e-10 \\
10	5.05088686255812e-10 \\
11	4.59844483914494e-10 \\
12	4.40626816804476e-10 \\
13	4.09137591059872e-10 \\
14	3.92924002046674e-10 \\
15	3.42844572519874e-10 \\
16	2.82015171978676e-10 \\
};

\addplot+[smooth,color=blue,loosely dotted, every mark/.append style={solid}, mark=-]
table[row sep=crcr]{
1	2.20951796356791e-14  \\
2	2.03428031568097e-14 \\
3	1.99131089472608e-14 \\
4	1.79194518858946e-14 \\
5	1.61191000882652e-14 \\
6	1.26154219573265e-14 \\
7	1.18606234615548e-14 \\
8	1.11214407663117e-14 \\
9	1.00231139509196e-14 \\
10	8.87896974453336e-15 \\
11	8.49176509834063e-15 \\
12	7.77271856866792e-15 \\
13	7.29625214103300e-15 \\
14	6.58678118872868e-15 \\
15	6.01411179508135e-15 \\
16	5.27698395537874e-15 \\
};

\end{axis}

\begin{axis}[%
name=SumRate,
at={($(ber.east)+(35,0em)$)},
		anchor= west,
ymode=log,
width  = 0.34\columnwidth,
height = 0.3\columnwidth,
scale only axis,
xmin   = 1,
xmax  = 16,
xlabel= {Index},
xmajorgrids,
ymin = 0.0,
ymax = 0.15,
xtick       ={5, 10},
xticklabels ={$5$, $10$},
ymajorgrids,
legend entries={Th. bound q=0, Th. bound q=1,  Th. bound q=2,  Computed q=0, Computed q=1, Computed q=2}, 
legend style={at={(0.8,1.3)},anchor=north east,draw=black,fill=white,legend cell align=left,font=\tiny, legend columns=3}
]

\addplot+[smooth,color=green,solid, every mark/.append style={solid}, mark=none]
table[row sep=crcr]{
1	0.123042468055957   \\
2	0.123120989069821  \\
3	0.123353500999150  \\
4	0.123493486457197  \\
5	0.123610742804532  \\
6	0.123796758936231  \\
7	0.123895703694608  \\
8	0.124092398750813  \\
9	0.124210871323918  \\
10	0.124410796222347  \\
11	0.124578663836926  \\
12	0.124717929573235  \\
13	0.124845998361352  \\
14	0.125020437889278  \\
15	0.125188405763068  \\
16	0.125367205893023 \\
};

\addplot+[smooth,color=green, dotted, every mark/.append style={solid}, mark=|]
table[row sep=crcr]{
1	1.33389711610771e-05  \\
2	1.33649181811220e-05  \\
3	1.34419528083205e-05  \\
4	1.34884777597324e-05  \\
5	1.35275329779057e-05  \\
6	1.35896484921060e-05  \\
7	1.36227677763915e-05  \\
8	1.36887699415555e-05  \\
9	1.37286291616586e-05  \\
10	1.37960718112648e-05  \\
11	1.38528745218753e-05  \\
12	1.39001199025594e-05  \\
13	1.39436636770897e-05  \\
14	1.40031232236514e-05  \\
15	1.40605401037707e-05  \\
16	1.41218360914270e-05 \\
};

\addplot+[smooth,color=green, loosely dashed, every mark/.append style={solid}, mark=none]
table[row sep=crcr]{
1	1.19388806835431e-09  \\
2	1.19776117058548e-09  \\
3	1.20928964244145e-09  \\
4	1.21627362618953e-09  \\
5	1.22214872249380e-09  \\
6	1.23151610325872e-09  \\
7	1.23652237168035e-09  \\
8	1.24652338484641e-09  \\
9	1.25257867196700e-09  \\
10	1.26285105982702e-09  \\
11	1.27152886425617e-09  \\
12	1.27876468301485e-09  \\
13	1.28544812231382e-09  \\
14	1.29459693180478e-09  \\
15	1.30345605836416e-09  \\
16	1.31294036350387e-09  \\
};

\addplot+[smooth,color=blue,densely dashed, every mark/.append style={solid}, mark=none]
table[row sep=crcr]{
1	0.0609748026822137  \\
2	0.0588447552116144 \\
3	0.0501846282403032 \\
4	0.0414325092595386 \\
5	0.0342425519853107 \\
6	0.0334124993654204 \\
7	0.0320718287823502 \\
8	0.0289253566729413 \\
9	0.0256027451336627 \\
10	0.0235779786067386 \\
11	0.0224215162584153 \\
12	0.0195925509664705 \\
13	0.0186762087795896 \\
14	0.0180726301408849 \\
15	0.0162902850875792 \\
16	0.0132287391188049  \\
};

\addplot+[smooth,color=blue,densely dotted, every mark/.append style={solid}, mark=none]
table[row sep=crcr]{
1	3.88923431936958e-06 \\
2	2.68577868131170e-06 \\
3	2.47834165661338e-06 \\
4	2.22176511482827e-06 \\
5	1.92291821396799e-06 \\
6	1.75602235508932e-06 \\
7	1.56254170341411e-06 \\
8	1.39050380819486e-06 \\
9	1.29942483585857e-06 \\
10	1.16807005573540e-06 \\
11	1.13561291237179e-06 \\
12	1.03367855754780e-06 \\
13	9.57188605713535e-07 \\
14	8.80701270900169e-07 \\
15	8.15249466055335e-07 \\
16	7.24889869514828e-07 \\
};

\addplot+[smooth,color=blue,loosely dotted, every mark/.append style={solid}, mark=-]
table[row sep=crcr]{
1	2.29562897658282e-10 \\
2	1.74614150693074e-10 \\
3	1.54112551087877e-10 \\
4	1.39687542674672e-10 \\
5	1.28261792281965e-10 \\
6	1.03149986119100e-10 \\
7	9.82729378439990e-11 \\
8	8.87571365116902e-11 \\
9	7.76423838203923e-11 \\
10	7.34997439090173e-11 \\
11	6.52695040480409e-11 \\
12	6.37746747037090e-11 \\
13	5.76231600386949e-11 \\
14	5.00619985467685e-11 \\
15	4.83552721081098e-11 \\
16	4.35426339898576e-11 \\
};

\end{axis}

\end{tikzpicture}%
\captionsetup{justification=centering,font=scriptsize}  
\caption{Principal angles between $\mathcal{R}({\bf P})$ and $\mathcal{R}({\bf V}_k)$ for  \texttt{LowRankLargeGap} (left), and \texttt{LowRankMediumGap} (right).} 
\label{figAnglVSte}      
\end{center}
\end{figure}
\begin{figure}[t]
\begin{center}       
%
%
%
\usetikzlibrary{positioning,calc}

\definecolor{mycolor1}{rgb}{0.00000,1.00000,1.00000}%
\definecolor{mycolor2}{rgb}{1.00000,0.00000,1.00000}%

\pgfplotsset{every axis label/.append style={font=\footnotesize},
every tick label/.append style={font=\footnotesize}
}

\begin{tikzpicture}[font=\footnotesize] 

\begin{axis}[%
name=ber,
ymode=log,
width  = 0.35\columnwidth,
height = 0.3\columnwidth,
scale only axis,
xmin  = 1,
xmax  = 16,
xlabel= {Index},
xmajorgrids,
ymin = 0.0,
ymax = 1.2,
xtick       ={5, 10},
xticklabels ={$5$, $10$},
ylabel={sin$\phi_i$},
ymajorgrids,
]
\addplot+[smooth,color=green,solid, every mark/.append style={solid}, mark=none]
table[row sep=crcr]{
1	0.988976894552712  \\
2	0.988976894552712 \\
3	0.988976894552712 \\
4	0.988976894552712 \\
5	0.988976894552712 \\
6	0.988976894552712 \\
7	0.988976894552712 \\
8	0.988976894552712 \\
9	0.988976894552712 \\ 
10	0.988976894552712 \\
11	0.988976894552712 \\
12	0.988976894552712 \\
13	0.988976894552712 \\
14	0.988976894552712 \\
15	0.988976894552712 \\
16	0.988976894552712 \\
};

\addplot+[smooth,color=green, dotted, every mark/.append style={solid}, mark=|]
table[row sep=crcr]{
1	0.881011353607130  \\
2	0.881011353607130 \\
3	0.881011353607130 \\ 
4	0.881011353607130 \\
5	0.881011353607130 \\ 
6	0.881011353607130 \\
7	0.881011353607130 \\
8	0.881011353607130 \\
9	0.881011353607130 \\
10	0.881011353607131 \\
11	0.881011353607131 \\
12	0.881011353607131 \\
13	0.881011353607131 \\
14	0.881011353607131 \\
15	0.881011353607131 \\
16	0.881011353607131  \\
};

\addplot+[smooth,color=green, loosely dashed, every mark/.append style={solid}, mark=none]
table[row sep=crcr]{
1	0.411219469314689  \\
2	0.411219469314689 \\
3	0.411219469314689 \\
4	0.411219469314689 \\
5	0.411219469314689 \\
6	0.411219469314689 \\
7	0.411219469314689 \\
8	0.411219469314689 \\
9	0.411219469314689 \\ 
10	0.411219469314690 \\ 
11	0.411219469314690 \\
12	0.411219469314690 \\
13	0.411219469314690 \\
14	0.411219469314691 \\
15	0.411219469314691 \\
16	0.411219469314691 \\
};

\addplot+[smooth,color=blue,densely dashed, every mark/.append style={solid}, mark=none]
table[row sep=crcr]{
1	0.158775028174782  \\
2	0.148990306297723 \\
3	0.135689140577601 \\
4	0.127537746146303 \\
5	0.100680653127888 \\
6	0.0963638337434248 \\
7	0.0848609324022176 \\
8	0.0738190022223206 \\
9	0.0659580159884620 \\
10	0.0622423578931712 \\
11	0.0539210499799164 \\
12	0.0430146993409694 \\
13	0.0409039825573083 \\
14	0.0352370495961822 \\
15	0.0303107065840316 \\
16	0.0248882715232639  \\
};

\addplot+[smooth,color=blue,densely dotted, every mark/.append style={solid}, mark=none]
table[row sep=crcr]{
1	0.000539621367669739  \\
2	0.000513038426971617 \\
3	0.000284869734173024 \\
4	0.000194939288870757 \\
5	0.000139944296982922 \\
6	0.000126296050274233 \\
7	8.52085690899776e-05 \\
8	7.32574051354941e-05 \\
9	5.17155543902550e-05 \\
10	3.70233363493988e-05 \\
11	3.61626656649455e-05 \\
12	2.42558270204331e-05 \\
13	2.24244306812919e-05 \\
14	1.97261975188281e-05 \\
15	1.37823185045091e-05 \\
16	1.09990331217238e-05  \\
};

\addplot+[smooth,color=blue,loosely dotted, every mark/.append style={solid}, mark=-]
table[row sep=crcr]{
1	3.09551503040988e-06  \\
2	1.29639740859655e-06 \\
3	8.60003423792043e-07 \\
4	7.14804754967445e-07 \\
5	4.35728972998496e-07 \\
6	2.59319795154871e-07 \\
7	1.58265919790107e-07 \\
8	1.30718732925067e-07 \\
9	1.13920100500136e-07 \\
10	7.05597710484807e-08 \\
11	5.35479205232701e-08 \\
12	3.55107765036015e-08 \\
13	3.05053486135328e-08 \\
14	1.94772229634728e-08 \\
15	1.33712068916843e-08 \\
16	7.82754000985577e-09 \\
};

\end{axis}

\begin{axis}[%
name=SumRate,
at={($(ber.east)+(35,0em)$)},
		anchor= west,
ymode=log,
width  = 0.35\columnwidth,
height = 0.3\columnwidth,
scale only axis,
xmin   = 1,
xmax  = 16,
xlabel= {Index},
xmajorgrids,
ymin = 0.0,
ymax = 1,
xtick       ={5, 10},
xticklabels ={$5$, $10$},
ymajorgrids,
]

\addplot+[smooth,color=green,solid, every mark/.append style={solid}, mark=none]
table[row sep=crcr]{
1	0.947442387384770   \\
2	0.947442387384770  \\
3	0.947442387384770  \\
4	0.947442387384770  \\
5	0.947442387384770  \\
6	0.947442387384770  \\
7	0.947442387384770  \\
8	0.947442387384770  \\  
9	0.947442387384770  \\
10	0.947442387384770  \\
11	0.947442387384770  \\
12	0.947442387384770  \\
13	0.947442387384770  \\
14	0.947442387384770  \\
15	0.947442387384770  \\
16	0.947442387384770 \\
};

\addplot+[smooth,color=green, dotted, every mark/.append style={solid}, mark=|]
table[row sep=crcr]{
1	0.186222997314121  \\
2	0.186222997314121  \\
3	0.186222997314121  \\
4	0.186222997314121  \\
5	0.186222997314121  \\
6	0.186222997314121  \\
7	0.186222997314121  \\
8	0.186222997314121  \\
9	0.186222997314121  \\
10	0.186222997314121  \\
11	0.186222997314121  \\
12	0.186222997314121  \\
13	0.186222997314121  \\
14	0.186222997314121  \\
15	0.186222997314121  \\
16	0.186222997314121  \\
};

\addplot+[smooth,color=green, loosely dashed, every mark/.append style={solid}, mark=none]
table[row sep=crcr]{
1	0.0131349451299846  \\
2	0.0131349451299846  \\
3	0.0131349451299846  \\
4	0.0131349451299846  \\
5	0.0131349451299846  \\
6	0.0131349451299846  \\
7	0.0131349451299846  \\  
8	0.0131349451299846  \\
9	0.0131349451299846  \\
10	0.0131349451299846  \\
11	0.0131349451299846  \\
12	0.0131349451299846  \\
13	0.0131349451299846  \\
14	0.0131349451299846  \\
15	0.0131349451299846  \\
16	0.0131349451299846  \\
};

\addplot+[smooth,color=blue,densely dashed, every mark/.append style={solid}, mark=none]
table[row sep=crcr]{
1	0.00957326964760631   \\
2	0.00585758767199018  \\
3	0.00421620590625491  \\
4	0.0040942986487033  \\
5	0.00351552145017214  \\
6	0.00301694274413711  \\
7	0.00222316503907576  \\
8	0.00200411982613269  \\
9	0.00180258865764310  \\
10	0.00144907642296092  \\
11	0.00128799557183006  \\
12	0.00107093568131187  \\
13	0.000826260567774307  \\
14	0.000695106542541868  \\
15	0.000600199194183523  \\
16	0.000594147807909450  \\
};

\addplot+[smooth,color=blue,densely dotted, every mark/.append style={solid}, mark=none]
table[row sep=crcr]{
1	1.27719855476397e-07  \\
2	7.65470912349157e-08  \\
3	3.80189111225090e-08  \\
4	2.50478163091774e-08  \\
5	1.24530071314874e-08  \\
6	1.08987486253865e-08   \\
7	9.44354887250065e-09  \\
8	4.90819574967344e-09  \\
9	3.32638828736924e-09  \\
10	2.10564759657277e-09  \\
11	1.42049063150163e-09  \\
12	1.00364841157686e-09  \\
13	7.03751583733237e-10  \\
14	5.98109219934374e-10  \\  
15	4.19684948441838e-10  \\
16	2.22420509196479e-10 \\
};

\addplot+[smooth,color=blue,loosely dotted, every mark/.append style={solid}, mark=-]
table[row sep=crcr]{
1	2.07084373280074e-12   \\
2	7.53028964089242e-13  \\
3	2.99278302813739e-13  \\
4	1.27472515863318e-13  \\
5	5.28274008363234e-14  \\
6	4.05197402125666e-14  \\
7	1.57877179524728e-14  \\
8	1.01882501596512e-14  \\
9	6.87747879165717e-15  \\
10	3.27864323979449e-15  \\
11	1.88748744252507e-15  \\
12	1.54162295142396e-15  \\
13	1.19110536452824e-15  \\
14	8.88964233015538e-16  \\
15	8.23207642603307e-16  \\
16	7.67380848552581e-16  \\
};

\end{axis}

\end{tikzpicture}%
\captionsetup{justification=centering,font=scriptsize}  
\caption{Principal angles between $\mathcal{R}({\bf P})$ and $\mathcal{R}({\bf V}_k)$ for  \texttt{LowRankSlowDecay} (left), and \texttt{LowRankFastDecay} (right).} 
\label{figAnglVTrop}      
\end{center}
\end{figure}


\begin{figure}[t]
\begin{center}       
%
%
%
\usetikzlibrary{positioning,calc}

\definecolor{mycolor1}{rgb}{0.00000,1.00000,1.00000}%
\definecolor{mycolor2}{rgb}{1.00000,0.00000,1.00000}%

\pgfplotsset{every axis label/.append style={font=\footnotesize},
every tick label/.append style={font=\footnotesize}
}

\begin{tikzpicture}[font=\footnotesize]

\begin{axis}[%
name=ber,
ymode=log,
width  = 0.33\columnwidth,
height = 0.3\columnwidth,
scale only axis,
xmin  = 17,
xmax  = 32,
xlabel= {$d$},
xmajorgrids,
ymin = 0.0 ,
ymax = 0.0072,
xtick       ={20,30},
xticklabels ={$20$,$30$},
ylabel={Magnitude},
ymajorgrids,
]

\addplot+[smooth,color=red,solid, every mark/.append style={solid}, mark=none]
table[row sep=crcr]{
17	0.00712083652325025   \\
18	0.00676094431693713  \\
19	0.00634236495946405  \\
20	0.00580488902435197  \\
21	0.00595466418412648  \\
22	0.00553068266106969  \\
23	0.00550227508996604  \\
24	0.00545865322575386  \\
25	0.00539078974546222  \\
26	0.00532867916544025  \\
27	0.00529941584949673  \\
28	0.00530446844386628  \\
29	0.00521724554976698  \\
30	0.00522706207894395  \\
31	0.00520138221742183  \\
32	0.00517057184770554 \\
};

\addplot+[smooth,color=red, loosely dashed, every mark/.append style={solid}, mark=none]
table[row sep=crcr]{
17	0.00483613900241773  \\
18	0.00483612326839021  \\
19	0.00483611351755948  \\
20	0.00483611148921960  \\
21	0.00483610672469349  \\
22	0.00483610131217257  \\
23	0.00483610010536568  \\
24	0.00483610032443428  \\
25	0.00483609801839351  \\
26	0.00483609677305874  \\
27	0.00483609567879053  \\
28	0.00483609531955303  \\
29	0.00483609513151047   \\
30	0.00483609468871020  \\
31	0.00483609355689491  \\
32	0.00483609300183010  \\
};

\addplot+[smooth,color=red, dotted, every mark/.append style={solid}, mark=|]
table[row sep=crcr]{
17	0.00483608511384454   \\
18	0.00483608511360121  \\
19	0.00483608511338164  \\
20	0.00483608511304233  \\
21	0.00483608511323290  \\
22	0.00483608511288172  \\
23	0.00483608511283051  \\
24	0.00483608511287852  \\
25	0.00483608511279321  \\
26	0.00483608511278197  \\
27	0.00483608511277842  \\
28	0.00483608511274503  \\
29	0.00483608511274216  \\
30	0.00483608511268794  \\
31	0.00483608511272109  \\
32	0.00483608511268990 \\
  };

\addplot+[smooth,color=blue,densely dashed, every mark/.append style={solid}, mark=none]
table[row sep=crcr]{
17	0.00483413381336874  \\
18	0.00482867046007228 \\
19	0.00482759553889071 \\
20	0.00482532216494997 \\
21	0.00481818403399306 \\
22	0.00480388402233870 \\
23	0.00480665208064950 \\
24	0.00480398077323274 \\
25	0.00479559343102103 \\
26	0.00479077117760516 \\
27	0.00479405724454778 \\
28	0.00477413626514107 \\
29	0.00478201063774573 \\
30	0.00478756886793069 \\
31	0.00477176371826010 \\
32	0.00476694065012709 \\
};

\addplot+[smooth,color=blue,densely dotted, every mark/.append style={solid}, mark=none]
table[row sep=crcr]{
17	0.00482724228547699 \\
18	0.00482253680289874\\
19	0.00481617764970981\\
20	0.00480666068470465\\
21	0.00479657940409838\\
22	0.00479985395448390\\
23	0.00479900149820058\\
24	0.00479104514592011\\
25	0.00478009232343638\\
26	0.00478063972190925\\
27	0.00477354804878865\\
28	0.00478397077184620\\
29	0.00477717799487420\\
30	0.00475860549617293\\
31	0.00474817234547512\\
32	0.00474561726507805 \\
};

\addplot+[smooth,color=blue,loosely dotted, every mark/.append style={solid}, mark=-]
table[row sep=crcr]{
17	0.00483004857977552  \\
18	0.00481520167793226 \\
19	0.00481536297465729 \\
20	0.00480814440898576 \\
21	0.00479784417254417 \\
22	0.00478674887011269 \\
23	0.00478894950036403 \\
24	0.00478334865513329 \\
25	0.00478097815767377 \\
26	0.00476042652735751 \\
27	0.00477081866775583 \\
28	0.00475539114800699 \\
29	0.00474750249147167 \\
30	0.00473006955076277 \\
31	0.00473682056964091 \\
32	0.00472724779970453 \\
};

\addplot+[smooth,color=cyan,loosely dotted, every mark/.append style={solid}, mark=triangle]
table[row sep=crcr]{
17	0.00479617817667860  \\
18	0.00476592516062344 \\
19	0.00473201653357643 \\
20	0.00472493852824964 \\
21	0.00470743709254221 \\
22	0.00469649318969371 \\
23	0.00468175881675619 \\
24	0.00467627998310049 \\
25	0.0046448593189550 \\
26	0.00461644476579934 \\
27	0.00459506091762965 \\
28	0.00457581711958483 \\
29	0.00456428979207927 \\
30	0.00454949537031759 \\ 
31	0.00454658437870722 \\
32	0.00453265555260447 \\
};

\end{axis}

\begin{axis}[%
name=SumRate,
at={($(ber.east)+(35,0em)$)},
		anchor= west,
ymode=log,
width  = 0.33\columnwidth,
height = 0.3\columnwidth,
scale only axis,
xmin  = 17,
xmax  = 32,
xlabel= {$d$},
xmajorgrids,
ymin = 0.0 ,
ymax = 0.10,
xtick       ={20,30},
xticklabels ={$20$,$30$},
ylabel={},
ymajorgrids,
legend entries = {Th. bound q=0, Th. bound q=1, Th. bound q=2, SVD, Computed q=0, Computed q=1, Computed q=2},
legend style={at={(0.97,1.34)},anchor=north east,draw=black,fill=white,legend cell align=left,font=\tiny, legend columns=4}
]

\addplot+[smooth,color=red,solid, every mark/.append style={solid}, mark=none]
table[row sep=crcr]{
17	0.100152294806159   \\
18	0.0950905256969505 \\
19	0.0892033405224928 \\
20	0.0816439129006393 \\
21	0.0837504527583370 \\
22	0.0777872878477394 \\
23	0.0773877444195747 \\
24	0.0767742168108018 \\
25	0.0758197386027222 \\
26	0.0749461731022881 \\
27	0.0745345938958548 \\
28	0.0746056570243485 \\
29	0.0733788948349404 \\
30	0.0735169611872425 \\
31	0.0731557821244531 \\
32	0.0727224441000739 \\
};
\addplot+[smooth,color=red, loosely dashed, every mark/.append style={solid}, mark=none]
table[row sep=crcr]{
17	0.0680187527844869  \\
18	0.0680185314904075 \\
19	0.0680183943480829 \\
20	0.0680183658201294 \\
21	0.0680182988085890 \\
22	0.0680182226832084 \\
23	0.0680182057098544 \\
24	0.0680182087909844 \\
25	0.0680181763572563 \\
26	0.0680181588420199 \\
27	0.0680181434514866 \\
28	0.0680181383989257 \\
29	0.0680181357541671 \\
30	0.0680181295263224 \\
31	0.0680181136077013 \\
32	0.0680181058008918 \\
};

\addplot+[smooth,color=red, dotted, every mark/.append style={solid}, mark=|]
table[row sep=crcr]{
17	0.0680179948588904  \\
18	0.0680179948554681 \\
19	0.0680179948523799 \\
20	0.0680179948476076 \\
21	0.0680179948502879 \\
22	0.0680179948453487 \\
23	0.0680179948446284 \\
24	0.0680179948453036 \\
25	0.0680179948441038 \\
26	0.0680179948439457 \\
27	0.0680179948438958 \\
28	0.0680179948434262 \\
29	0.0680179948433858 \\
30	0.0680179948426232 \\
31	0.0680179948430895 \\
32	0.0680179948426508 \\
};

\addplot+[smooth,color=cyan,loosely dotted, every mark/.append style={solid}, mark=triangle]
table[row sep=crcr]{
17	0.0678458539838310   \\
18	0.0676761152674477  \\
19	0.0675080923671837  \\
20	0.0673420415088685  \\
21	0.0671760784095531  \\
22	0.0670109360217135  \\
23	0.0668461554483528  \\
24	0.0666820038136758  \\
25	0.0665178324821752  \\
26	0.0663554618703980  \\
27	0.0661946807361299  \\
28	0.0660349996056746  \\
29	0.0658762709221579  \\
30	0.0657179604773571  \\
31	0.0655602960729963  \\
32	0.0654024540186852  \\
};

\addplot+[smooth,color=blue,densely dashed, every mark/.append style={solid}, mark=none]
table[row sep=crcr]{
17	0.0679186370168970
18	0.0678058118689234
19	0.0676966011480053
20	0.0675823084443789
21	0.0674772848156234
22	0.0673609164228069
23	0.0672568668787803
24	0.0671475257271150
25	0.0670379007483928
26	0.0669272969790482
27	0.0668266007451107
28	0.0667005071175799
29	0.0666075496886779
30	0.0664949374300568
31	0.0663923212760850
32	0.0662776121944603 \\
};

\addplot+[smooth,color=blue,densely dotted, every mark/.append style={solid}, mark=none]
table[row sep=crcr]{
17	0.0678869102576619  \\
18	0.0677591832351234  \\
19	0.0676299950791068  \\
20	0.0674974494641324  \\
21	0.067365000739007  \\
22	0.0672358140166899  \\
23	0.0671071961307210  \\
24	0.0669733967802615  \\
25	0.066845743873200  \\
26	0.0667213806372160  \\
27	0.0665858606917689  \\
28	0.066459041582450  \\
29	0.0663293816056834  \\
30	0.0661992522290929  \\
31	0.066066952753395  \\
32	0.0659368011256660 \\
};

\addplot+[smooth,color=blue,loosely dotted, every mark/.append style={solid}, mark=-]
table[row sep=crcr]{
17	0.0678756036942463   \\
18	0.0677360469549489  \\
19	0.0675953562714042  \\
20	0.0674539549893508  \\
21	0.0673065848999143  \\
22	0.0671686123311853  \\
23	0.0670299590122107  \\
24	0.0668887758369341  \\
25	0.0667507399435806  \\
26	0.0666116787645386  \\
27	0.0664745987529607  \\
28	0.0663327189137631  \\
29	0.0661986510912477  \\
30	0.0660475009794952  \\
31	0.0659146804898712  \\
32	0.0657805136407269 \\
};

\end{axis}

\end{tikzpicture}%
\captionsetup{justification=centering,font=scriptsize}  
\caption{{\bf Q}-based low-rank approximation error for \texttt{LowRankLargeGap}. Left: Spectral norm. Right: Frobenius norn.} 
\label{figQLRSteLarge}      
\end{center}
\end{figure}
\begin{figure}[t]
\begin{center}       
%
%
%
\usetikzlibrary{positioning,calc}

\definecolor{mycolor1}{rgb}{0.00000,1.00000,1.00000}%
\definecolor{mycolor2}{rgb}{1.00000,0.00000,1.00000}%

\pgfplotsset{every axis label/.append style={font=\footnotesize},
every tick label/.append style={font=\footnotesize}
}

\begin{tikzpicture}[font=\footnotesize]

\begin{axis}[%
name=ber,
ymode=log,
width  = 0.34\columnwidth,
height = 0.3\columnwidth,
scale only axis,
xmin  = 17,
xmax  = 32,
xlabel= {$d$},
xmajorgrids,
ymin = 0.0 ,
ymax = 0.023,
xtick       ={20,30},
xticklabels ={$20$,$30$},
ylabel={Magnitude},
ymajorgrids,
]

\addplot+[smooth,color=red,solid, every mark/.append style={solid}, mark=none]
table[row sep=crcr]{
17	0.0223295085307701  \\
18	0.0160713172830099 \\
19	0.0173982251455226 \\
20	0.0138833841876508 \\
21	0.0133793291086175 \\
22	0.0127650337740649 \\
23	0.0123022191016307 \\
24	0.0120137358065462 \\
25	0.0117409603629494 \\
26	0.0117582819072071 \\
27	0.0115646393410957 \\
28	0.0112981274967772 \\
29	0.0113052960306909 \\
30	0.0110978192469036 \\
31	0.0111147627632177 \\
32	0.0111397498393876 \\
};

\addplot+[smooth,color=red, loosely dashed, every mark/.append style={solid}, mark=none]
table[row sep=crcr]{
17	0.00966550288680286  \\
18	0.00966506790497258 \\
19	0.00966493646591272 \\
20	0.00966470641880802 \\
21	0.00966472027558544 \\
22	0.00966462203960746 \\
23	0.00966460570813845 \\
24	0.00966457179926848 \\
25	0.00966453409248194 \\
26	0.00966450941161476 \\
27	0.00966449362579606 \\
28	0.00966448395788840 \\
29	0.00966447378257479 \\
30	0.00966446279699610 \\
31	0.00966445969846994 \\
32	0.00966445996654491  \\
};

\addplot+[smooth,color=red, dotted, every mark/.append style={solid}, mark=|]
table[row sep=crcr]{
17	0.00966431927215745  \\
18	0.00966431895912460 \\
19	0.00966431893918776 \\
20	0.00966431893213533 \\
21	0.00966431892469008 \\
22	0.00966431892164021 \\
23	0.00966431891679000 \\
24	0.00966431891420571 \\
25	0.00966431891476494 \\
26	0.00966431891275942 \\
27	0.00966431891259128 \\
28	0.00966431891005024 \\
29	0.00966431890769136 \\
30	0.00966431890812681 \\
31	0.00966431890732722 \\
32	0.00966431890745407 \\
  };

\addplot+[smooth,color=blue,densely dashed, every mark/.append style={solid}, mark=none]
table[row sep=crcr]{
17	0.0102780869425159   \\
18	0.00966299020968616  \\
19	0.00974511436268644  \\
20	0.00963901762236972  \\
21	0.00962771031182827  \\
22	0.00961852126745979  \\
23	0.00961914559394954  \\
24	0.0096009098237081  \\
25	0.0095987882088165  \\
26	0.00958267319654213  \\
27	0.00958580079027226  \\
28	0.00957744134493390  \\
29	0.00956638137925900  \\
30	0.00957811552543239  \\
31	0.00955179464728145  \\
32	0.00955846581526356  \\
};

\addplot+[smooth,color=blue,densely dotted, every mark/.append style={solid}, mark=none]
table[row sep=crcr]{
17	0.00966114697592875   \\
18	0.00964619373635689  \\
19	0.00964416348635501  \\
20	0.00962467732672953  \\
21	0.00962319706179646  \\
22	0.00961504749553757  \\
23	0.00959387183519250  \\
24	0.0095903670394442  \\
25	0.00960320048601551  \\
26	0.0095743337456485  \\
27	0.00957568030289069  \\
28	0.00955044572675784  \\
29	0.00954476713550320  \\
30	0.00952605507939073  \\
31	0.00952621561093939  \\
32	0.00951075785115297  \\
};

\addplot+[smooth,color=blue,loosely dotted, every mark/.append style={solid}, mark=-]
table[row sep=crcr]{
17	0.0096540544548457  \\
18	0.00964547676806721  \\
19	0.00963746711347087  \\
20	0.00962449223393988  \\
21	0.00961040201263568  \\
22	0.0096063668340874  \\
23	0.00958179467409565  \\
24	0.00958049497822815  \\
25	0.00956102717213725  \\
26	0.00956230271623304  \\
27	0.00956702159568664  \\
28	0.00953812450028670  \\
29	0.00951378154129883  \\
30	0.00952997635880917  \\
31	0.00950096907254317  \\
32	0.00947979345126659  \\
};

\addplot+[smooth,color=cyan,loosely dotted, every mark/.append style={solid}, mark=triangle]
table[row sep=crcr]{
17	0.00961403801113749  \\
18	0.00958622717919707 \\
19	0.00954741827476937 \\
20	0.0094744622320338 \\
21	0.00943065256704664 \\
22	0.00939014848068508 \\
23	0.00937336607467642 \\
24	0.00935018396980432 \\
25	0.00929105790837784 \\
26	0.00927596169153734 \\
27	0.00924101960926558 \\
28	0.00923337283246877 \\
29	0.00921226549953755 \\
30	0.00916792219010894 \\
31	0.00913438771632940 \\
32	0.00911566525384084 \\
};

\end{axis}

\begin{axis}[%
name=SumRate,
at={($(ber.east)+(35,0em)$)},
		anchor= west,
ymode=log,
width  = 0.34\columnwidth,
height = 0.3\columnwidth,
scale only axis,
xmin  = 17,
xmax  = 32,
xlabel= {$d$},
xmajorgrids,
ymin = 0.0 ,
ymax = 0.32,
xtick       ={20,30},
xticklabels ={$20$,$30$},
ylabel={},
ymajorgrids,
]

\addplot+[smooth,color=red,solid, every mark/.append style={solid}, mark=none]
table[row sep=crcr]{
17	0.314523590615222  \\
18	0.226373473952784 \\
19	0.245063711794694 \\
20	0.195555215134865 \\
21	0.188455318014114 \\
22	0.179802625365038 \\
23	0.173283622389096 \\
24	0.169220174164185 \\
25	0.165377979794635 \\
26	0.165621963413327 \\
27	0.162894399790268 \\
28	0.159140431712510 \\
29	0.159241404513720 \\
30	0.156318978213286 \\
31	0.156557637097400 \\
32	0.156909594011503 \\
};
\addplot+[smooth,color=red, loosely dashed, every mark/.append style={solid}, mark=none]
table[row sep=crcr]{
17	0.136144002850303  \\
18	0.136137875888444 \\
19	0.136136024495924 \\
20	0.136132784154056 \\
21	0.136132979334500 \\
22	0.136131595625921 \\
23	0.136131365588065 \\
24	0.136130887962705 \\
25	0.136130356841571 \\
26	0.136130009197782 \\
27	0.136129786845715 \\
28	0.136129650667839 \\
29	0.136129507342868 \\
30	0.136129352604861 \\
31	0.136129308960390 \\
32	0.136129312736376 \\
};

\addplot+[smooth,color=red, dotted, every mark/.append style={solid}, mark=|]
table[row sep=crcr]{
17	0.136127330977402  \\
18	0.136127326568159 \\
19	0.136127326287338 \\
20	0.136127326188000 \\
21	0.136127326083130 \\
22	0.136127326040171 \\
23	0.136127325971853 \\
24	0.136127325935452 \\
25	0.136127325943329 \\
26	0.136127325915080 \\
27	0.136127325912711 \\
28	0.136127325876919 \\
29	0.136127325843693 \\
30	0.136127325849827 \\ 
31	0.136127325838564 \\
32	0.136127325840351 \\
};

\addplot+[smooth,color=cyan,loosely dotted, every mark/.append style={solid}, mark=triangle]
table[row sep=crcr]{
17	0.135783834566655   \\
18	0.135443050773176  \\
19	0.135103383566860  \\
20	0.134765615256640  \\
21	0.134432159920587  \\
22	0.134100963505390  \\
23	0.133771796446728  \\
24	0.133442997317262  \\
25	0.133115016405910  \\
26	0.132790375538624  \\
27	0.132465997034658  \\
28	0.132143270456604  \\
29	0.131820289610908    \\
30	0.131497995868643  \\
31	0.131178016909030  \\
32	0.130859600646012  \\
};

\addplot+[smooth,color=blue,densely dashed, every mark/.append style={solid}, mark=none]
table[row sep=crcr]{
17	0.136031745619151    \\
18	0.135717583766316  \\
19	0.135516873134317  \\
20	0.135270121220135  \\
21	0.135042190206720  \\
22	0.134823810916549  \\
23	0.134598832054316  \\
24	0.134385609647074  \\
25	0.134158971337449  \\
26	0.13394853687393  \\
27	0.133738193492314  \\
28	0.13350703756229  \\
29	0.133288571429403  \\
30	0.133062452082680  \\
31	0.13284392688052  \\
32	0.132631420520855  \\
};

\addplot+[smooth,color=blue,densely dotted, every mark/.append style={solid}, mark=none]
table[row sep=crcr]{
17	0.135867476273159   \\
18	0.135604462279871  \\
19	0.135341913183973  \\
20	0.13507830449056  \\
21	0.134816488881717  \\
22	0.134548014995797  \\
23	0.134299843407920  \\
24	0.134036628119211  \\
25	0.133775237572560  \\
26	0.133499552178134  \\
27	0.133242598584468  \\
28	0.132985635728807  \\
29	0.132723986781743  \\
30	0.132460084523542  \\
31	0.132204217112395  \\
32	0.131934874945964 \\
};

\addplot+[smooth,color=blue,loosely dotted, every mark/.append style={solid}, mark=-]
table[row sep=crcr]{
17	0.135845712788748   \\
18	0.135563464840949  \\
19	0.135275549181456  \\
20	0.134991703849653  \\
21	0.134709140842719  \\
22	0.134439566867335  \\
23	0.134139224405708  \\
24	0.133849445985594  \\
25	0.133579752413638  \\
26	0.133303829837086  \\
27	0.133005274469443  \\
28	0.132731452592974  \\
29	0.132448308039694  \\
30	0.132171431424238  \\
31	0.131891519261222  \\
32	0.131607891853417  \\
};

\end{axis}

\end{tikzpicture}%
\captionsetup{justification=centering,font=scriptsize}  
\caption{{\bf Q}-based low-rank approximation error for \texttt{LowRankMediumGap}. Left: Spectral norm. Right: Frobenius norn.} 
\label{figQLRSteMedium}      
\end{center}
\end{figure}
\begin{figure}[t]
\begin{center}       
%
%
%
\usetikzlibrary{positioning,calc}

\definecolor{mycolor1}{rgb}{0.00000,1.00000,1.00000}%
\definecolor{mycolor2}{rgb}{1.00000,0.00000,1.00000}%

\pgfplotsset{every axis label/.append style={font=\footnotesize},
every tick label/.append style={font=\footnotesize}
}

\begin{tikzpicture}[font=\footnotesize]

\begin{axis}[%
name=ber,
ymode=log,
width  = 0.34\columnwidth,
height = 0.3\columnwidth,
scale only axis,
xmin  = 17,
xmax  = 32,
xlabel= {$d$},
xmajorgrids,
ymin = 0.05 ,
ymax = 32,
xtick       ={20,30},
xticklabels ={$20$,$30$},
ylabel={Magnitude},
ymajorgrids,
]

\addplot+[smooth,color=red,solid, every mark/.append style={solid}, mark=none]
table[row sep=crcr]{
17	31.3117467429126  \\
18	20.473668476159 \\
19	21.476489085728 \\
20	11.4270404865035 \\
21	10.7606811418295 \\
22	7.92876909423696 \\
23	8.22583370072692 \\
24	7.51071277748083 \\
25	5.91247370179514 \\
26	5.79141832625403 \\
27	5.41921686753839 \\
28	4.99196273103189 \\
29	4.99608379445625 \\
30	4.60723239528336 \\
31	4.19463595024129 \\
32	4.19525175417860\\
};

\addplot+[smooth,color=red, loosely dashed, every mark/.append style={solid}, mark=none]
table[row sep=crcr]{
17	7.15591882264547 \\
18	5.00101064525085 \\
19	3.24061586751695 \\
20	3.40104165012516 \\
21	2.85963273882908 \\
22	2.31762919470654 \\
23	2.18043437934302 \\
24	1.98383452401396 \\ 
25	2.02165583884263 \\
26	1.70989622513196 \\
27	1.81143431527117 \\
28	1.60670955959153 \\
29	1.45460664394635 \\
30	1.50720988184987 \\
31	1.44593267702146 \\
32	1.42347544981487 \\
};

\addplot+[smooth,color=red, dotted, every mark/.append style={solid}, mark=|]
table[row sep=crcr]{
17	3.44161106094365  \\
18	1.65617037366492 \\
19	1.56374825409132 \\
20	1.12649486612484 \\
21	1.13171765487545 \\
22	1.03636584071095 \\
23	0.90907151413037 \\
24	0.88275251921608 \\
25	0.881495794879377 \\
26	0.831799801003188 \\
27	0.787012306863727 \\
28	0.782474368788199 \\
29	0.767124104824553 \\
30	0.748202433068175 \\
31	0.748469295532862 \\
32	0.750856063207420 \\
  };

\addplot+[smooth,color=blue,densely dashed, every mark/.append style={solid}, mark=none]
table[row sep=crcr]{
17	0.689622744590171 \\
18	0.526146639450708 \\
19	0.353276668103683 \\
20	0.268892007547359 \\
21	0.206822564016141 \\
22	0.207543900136060 \\
23	0.159531062644482 \\
24	0.152555676225062 \\
25	0.133470759427520 \\
26	0.127841641562393 \\
27	0.115788185798683 \\
28	0.103323075618980 \\
29	0.0975952039848593 \\
30	0.0949684877131030 \\
31	0.0870480312114759 \\
32	0.0816442677236860  \\
};

\addplot+[smooth,color=blue,densely dotted, every mark/.append style={solid}, mark=none]
table[row sep=crcr]{
17	0.417519646405726 \\
18	0.28587600873224 \\
19	0.211655465768610 \\
20	0.183612661153499 \\
21	0.170124945670102 \\
22	0.144814137159100 \\
23	0.126984554014146 \\
24	0.119588129550518 \\
25	0.105195986380609 \\
26	0.09438521288870 \\
27	0.0863755101111228 \\
28	0.08258554202462 \\
29	0.0796402787885433 \\
30	0.0720325923511177 \\
31	0.0673284358722109 \\
32	0.0658352278345123  \\
};

\addplot+[smooth,color=blue,loosely dotted, every mark/.append style={solid}, mark=-]
table[row sep=crcr]{
17	0.344722943296900 \\
18	0.252651532210168 \\
19	0.215234642642353 \\
20	0.173645243483954 \\
21	0.153583836442379 \\
22	0.132924827660386 \\
23	0.120345462988428 \\
24	0.112858538198128 \\
25	0.102366021720951 \\
26	0.0917977827926556 \\
27	0.0799967912260030 \\
28	0.077414623845036 \\
29	0.0736311520818248 \\
30	0.067999141305468 \\
31	0.065559360512153 \\
32	0.0609016695716029 \\
};

\addplot+[smooth,color=cyan,loosely dotted, every mark/.append style={solid}, mark=triangle]
table[row sep=crcr]{
17	0.333333333333333 \\
18	0.250000000000000 \\
19	0.200000000000000 \\
20	0.166666666666667 \\
21	0.142857142857143 \\
22	0.125000000000000 \\
23	0.111111111111111 \\
24	0.100000000000000 \\
25	0.0909090909090908 \\
26	0.083333333333333 \\
27	0.0769230769230769 \\
28	0.0714285714285713 \\
29	0.0666666666666667 \\
30	0.062500000000000 \\
31	0.0588235294117648 \\
32	0.0555555555555556 \\
};

\end{axis}

\begin{axis}[%
name=SumRate,
at={($(ber.east)+(35,0em)$)},
		anchor= west,
ymode=log,
width  = 0.34\columnwidth,
height = 0.3\columnwidth,
scale only axis,
xmin  = 17,
xmax  = 32,
xlabel= {$d$},
xmajorgrids,
ymin = 0.22 ,
ymax = 50.25,
xtick       ={20,30},
xticklabels ={$20$,$30$},
ylabel={},
ymajorgrids,
legend entries = {Th. bound q=0, Th. bound q=1, Th. bound q=2, SVD, Computed q=0, Computed q=1, Computed q=2},
legend style={at={(0.97,1.3)},anchor=north east,draw=black,fill=white,legend cell align=left,font=\tiny, legend columns=4}
]

\addplot+[smooth,color=red,solid, every mark/.append style={solid}, mark=none]
table[row sep=crcr]{
17	50.2418805463157  \\
18	32.8514283910686 \\
19	34.4605239707255 \\
20	18.3354830963168 \\
21	17.2662630725691 \\
22	12.7222627655591 \\
23	13.1989236365208 \\
24	12.0514622605426 \\
25	9.48697624242451 \\
26	9.29273445299759 \\
27	8.69551126447686 \\
28	8.00995221644523 \\
29	8.01656475000940 \\
30	7.39262557127501 \\
31	6.73058581105850 \\
32	6.73157391140666  \\
};
\addplot+[smooth,color=red, loosely dashed, every mark/.append style={solid}, mark=none]
table[row sep=crcr]{
17	11.4821706255613 \\
18	8.0244702255845 \\
19	5.19979407885109 \\
20	5.45720843112371 \\
21	4.5884800886461 \\
22	3.71879761634367 \\
23	3.4986589705612 \\
24	3.1832007967338 \\
25	3.24388773308802 \\
26	2.74364769857885 \\
27	2.90657264293157 \\
28	2.57807750006452 \\
29	2.33401776806246 \\
30	2.41842333051152 \\
31	2.32009978342618 \\
32	2.28406559676840 \\
};

\addplot+[smooth,color=red, dotted, every mark/.append style={solid}, mark=|]
table[row sep=crcr]{
17	5.52230487907702  \\
18	2.65744082440416 \\
19	2.50914309034354 \\
20	1.80753954624690 \\
21	1.81591987490379 \\
22	1.66292124162851 \\
23	1.45866862031052 \\
24	1.41643795814275 \\
25	1.41442145633200 \\
26	1.33468076960321 \\
27	1.26281611289788 \\
28	1.25553467489351 \\
29	1.23090410622073 \\
30	1.20054296476394 \\
31	1.20097116419284 \\
32	1.20480089932009 \\
};
\addplot+[smooth,color=cyan,loosely dotted, every mark/.append style={solid}, mark=triangle]
table[row sep=crcr]{
17	0.627424093054539  \\
18	0.531554213071630  \\
19	0.469094746756132  \\
20	0.424322850473787  \\
21	0.390220583332584  \\
22	0.363130748341580  \\
23	0.340938323441522  \\
24	0.322324776240937  \\
25	0.306420073393977  \\
26	0.292623988368782  \\
27	0.280507315634385  \\
28	0.269753951520816  \\
29	0.260125265102735  \\
30	0.251437286615020  \\
31	0.243545599632438  \\
32	0.236335040753306  \\
};

\addplot+[smooth,color=blue,densely dashed, every mark/.append style={solid}, mark=none]
table[row sep=crcr]{
17	0.959409594377126   \\
18	0.769485551548882  \\
19	0.606172893966728  \\
20	0.514825105470633  \\
21	0.450149608063964  \\
22	0.427437681752296  \\
23	0.385518918103848  \\
24	0.368924986830423  \\
25	0.346834380710657  \\
26	0.33668693383333  \\
27	0.317082087320559  \\
28	0.302212755329866  \\
29	0.292064938001253  \\
30	0.285214461162300  \\
31	0.274875581025456  \\
32	0.264956785505375  \\
};

\addplot+[smooth,color=blue,densely dotted, every mark/.append style={solid}, mark=none]
table[row sep=crcr]{
17	0.686148065520150   \\
18	0.555016865996210  \\
19	0.481205736123134  \\
20	0.435202367108696  \\
21	0.408400590123414  \\
22	0.374787779345583  \\
23	0.349737795042259  \\
24	0.334838747574436  \\
25	0.315203572359814  \\
26	0.301315597796917  \\
27	0.287833390249489  \\
28	0.277892473196315  \\
29	0.268935677393630  \\
30	0.258237261564300  \\
31	0.249839398343478  \\
32	0.243723148665432 \\
};

\addplot+[smooth,color=blue,loosely dotted, every mark/.append style={solid}, mark=-]
table[row sep=crcr]{
17	0.636431805218337  \\
18	0.533385878342474  \\
19	0.476580452714884  \\
20	0.428104782734331  \\
21	0.395461434561671  \\
22	0.367403170734356  \\
23	0.346405737981520  \\
24	0.328645270922117  \\
25	0.312098231056703  \\
26	0.297665337840883  \\
27	0.283152233506987  \\
28	0.273527574371601  \\
29	0.264137436842870  \\
30	0.254928392192765  \\
31	0.247356246110517  \\
32	0.239552870409499  \\
};

\end{axis}

\end{tikzpicture}%
\captionsetup{justification=centering,font=scriptsize}  
\caption{{\bf Q}-based low-rank approximation error for \texttt{LowRankSlowDecay}. Left: Spectral norm. Right: Frobenius norn.} 
\label{figQLRTropSlow}      
\end{center}
\end{figure}
\begin{figure}[t]
\begin{center}       
%
%
%
\usetikzlibrary{positioning,calc}

\definecolor{mycolor1}{rgb}{0.00000,1.00000,1.00000}%
\definecolor{mycolor2}{rgb}{1.00000,0.00000,1.00000}%

\pgfplotsset{every axis label/.append style={font=\footnotesize},
every tick label/.append style={font=\footnotesize}
}

\begin{tikzpicture}[font=\footnotesize]

\begin{axis}[%
name=ber,
ymode=log,
width  = 0.34\columnwidth,
height = 0.3\columnwidth,
scale only axis,
xmin  = 17,
xmax  = 32,
xlabel= {$d$},
xmajorgrids,
ymin = 0.0030 ,
ymax = 7.9,
xtick       ={20,30},
xticklabels ={$20$,$30$},
ylabel={Magnitude},
ymajorgrids,
]

\addplot+[smooth,color=red,solid, every mark/.append style={solid}, mark=none]
table[row sep=crcr]{
17	7.88171885337806 \\
18	5.0630645607837 \\
19	3.1844657383680 \\
20	3.10361215858711 \\
21	2.32438218248270 \\
22	2.02700290943897 \\
23	1.79199825172970 \\
24	1.76831461283612 \\
25	1.75519875981023 \\
26	1.44814463854794 \\
27	1.52931044232375 \\
28	1.38562553979450 \\
29	1.37052397486774 \\
30	1.20443015258838 \\
31	1.21646874791023 \\
32	1.16941505101069 \\
};

\addplot+[smooth,color=red, loosely dashed, every mark/.append style={solid}, mark=none]
table[row sep=crcr]{
17	0.736502272693119  \\
18	0.547631882403940 \\
19	0.491481247130649 \\
20	0.39534233422972 \\
21	0.379485513789138 \\
22	0.362784824794839 \\
23	0.354911715053827 \\
24	0.33889349892831 \\
25	0.33778642709568 \\
26	0.34451552398533 \\
27	0.324840972951782 \\
28	0.327579941364230 \\
29	0.314136344124639 \\
30	0.31668095549338 \\
31	0.308357549583015 \\
32	0.306931758474450 \\
};

\addplot+[smooth,color=red, dotted, every mark/.append style={solid}, mark=|]
table[row sep=crcr]{
17	0.293890848884713  \\
18	0.265724015091821 \\
19	0.265685235430887 \\
20	0.261374260226967 \\
21	0.257845345830659 \\
22	0.258186413972628 \\
23	0.255738835201322 \\
24	0.255352186974539 \\
25	0.255911030968760 \\
26	0.255105294325161 \\
27	0.254807156700262 \\
28	0.254681438598483 \\
29	0.254249030582578 \\
30	0.254057885995541 \\
31	0.253803709791979 \\
32	0.253710183597908 \\
  };

\addplot+[smooth,color=blue,densely dashed, every mark/.append style={solid}, mark=none]
table[row sep=crcr]{
17	0.162046982700268   \\
18	0.0768056896092715  \\
19	0.0680824037747032  \\
20	0.0397176292195313  \\
21	0.0272783666042366  \\
22	0.0227259692726274  \\
23	0.0169231038067254  \\
24	0.0134422838404172  \\
25	0.0128958461206434  \\
26	0.0106985138357546  \\
27	0.00905149421901623  \\
28	0.00807695966801696  \\
29	0.00681042921949265  \\
30	0.0058187731636640  \\
31	0.0049296456568994  \\
32	0.00479210325264493  \\
};

\addplot+[smooth,color=blue,densely dotted, every mark/.append style={solid}, mark=none]
table[row sep=crcr]{
17	0.116178675162133  \\
18	0.0642850418908800 \\
19	0.0423220560893224 \\
20	0.0319661708499032 \\
21	0.0233803453104009 \\
22	0.0175334074247823 \\
23	0.0134762875989659 \\
24	0.0114243153387270 \\
25	0.00920046196648740 \\
26	0.0084083797620255 \\
27	0.00653510707749451 \\
28	0.00599211783869513 \\
29	0.00523828700279815 \\
30	0.00464315352201131 \\
31	0.00397763714007181 \\
32	0.00350013981042489  \\
};

\addplot+[smooth,color=blue,loosely dotted, every mark/.append style={solid}, mark=-]
table[row sep=crcr]{
17	0.115183277219699  \\
18	0.0659839021423997 \\
19	0.0414981813027888 \\
20	0.0306500423583057 \\
21	0.0210861544860526 \\
22	0.0159452375194134 \\
23	0.013039749528868 \\
24	0.0108638236566727 \\
25	0.00895780963812113 \\
26	0.00756534952626258 \\
27	0.00625916627186334 \\
28	0.00554284728944445 \\
29	0.00494053383888694 \\
30	0.00431131283009108 \\
31	0.00364676413938873 \\
32	0.00329639906203906  \\
};

\addplot+[smooth,color=cyan, loosely dotted, every mark/.append style={solid}, mark=triangle]
table[row sep=crcr]{
17	0.111111111111111  \\
18	0.0625000000000001 \\
19	0.0400000000000001 \\
20	0.0277777777777778 \\
21	0.0204081632653061 \\
22	0.0156250000000000 \\
23	0.0123456790123457 \\
24	0.0100000000000000 \\
25	0.00826446280991735 \\
26	0.0069444444444444 \\
27	0.00591715976331360 \\
28	0.00510204081632654 \\
29	0.00444444444444445 \\
30	0.00390625000000000 \\
31	0.00346020761245675 \\
32	0.00308641975308642 \\
};

\end{axis}

\begin{axis}[%
name=SumRate,
at={($(ber.east)+(35,0em)$)},
		anchor= west,
ymode=log,
width  = 0.34\columnwidth,
height = 0.3\columnwidth,
scale only axis,
xmin  = 17,
xmax  = 32,
xlabel= {$d$},
xmajorgrids,
ymin = 0.0078 ,
ymax = 9.04,
xtick       ={20,30},
xticklabels ={$20$,$30$},
ylabel={},
ymajorgrids,
]

\addplot+[smooth,color=red,solid, every mark/.append style={solid}, mark=none]
table[row sep=crcr]{
17	9.04569914915653   \\
18	5.81078310982600  \\
19	3.65475089329373  \\
20	3.56195677421439  \\
21	2.66764931882681  \\
22	2.32635277080346  \\
23	2.05664238505715  \\
24	2.02946112216590  \\
25	2.01440830656009  \\
26	1.66200811884508  \\
27	1.75516057147803  \\
28	1.59025613569008  \\
29	1.57292435622038  \\
30	1.38230162851060  \\
31	1.39611809589362  \\
32	1.34211546094473  \\
};
\addplot+[smooth,color=red, loosely dashed, every mark/.append style={solid}, mark=none]
table[row sep=crcr]{
17	0.845269681066666  \\
18	0.628506718504578 \\
19	0.56406370002536 \\
20	0.453726894208446 \\
21	0.435528322318728 \\
22	0.416361258504810 \\
23	0.407325439870532 \\
24	0.388941637216183 \\
25	0.387671071883834 \\
26	0.395393928679619 \\
27	0.372813819841052 \\
28	0.375957281908013 \\
29	0.360528320488050 \\
30	0.363448722664450 \\
31	0.353896107662301 \\
32	0.352259754265594 \\
};

\addplot+[smooth,color=red, dotted, every mark/.append style={solid}, mark=|]
table[row sep=crcr]{
17	0.337292949819181  \\
18	0.30496640922383  \\
19	0.30492190254294  \\
20	0.299974278114876  \\
21	0.295924210033794  \\
22	0.296315647467557  \\
23	0.293506607761757  \\
24	0.29306285892966  \\
25	0.29370423357612  \\
26	0.292779504921494  \\
27	0.292437338027438  \\
28	0.292293053747899  \\
29	0.291796787274262  \\
30	0.291577414259251  \\
31	0.291285701054167  \\
32	0.291178362816167  \\
};

\addplot+[smooth,color=cyan, loosely dotted, every mark/.append style={solid}, mark=triangle]
table[row sep=crcr]{
17	0.140795003545487   \\
18	0.0864728512946573  \\
19	0.0597603883105515  \\
20	0.0443993694890806  \\
21	0.0346366723684059  \\
22	0.0279858168523434  \\
23	0.0232177802533512  \\
24	0.0196634058498785  \\
25	0.0169307273800340  \\
26	0.0147766093567004  \\
27	0.0130431160172158  \\
28	0.0116236868408425  \\
29	0.0104441024163167  \\
30	0.00945125329587347  \\
31	0.00860624190923420  \\
32	0.00787999765729378  \\
};

\addplot+[smooth,color=blue,densely dashed, every mark/.append style={solid}, mark=none]
table[row sep=crcr]{
17	0.188052955040075   \\
18	0.0990374431831871  \\
19	0.0866266702954142  \\
20	0.054955389883874  \\
21	0.0409044108429085  \\
22	0.0343467831951632  \\
23	0.0281276169445487  \\
24	0.0231006842928204  \\
25	0.0209902543643290  \\
26	0.0181678482633782  \\
27	0.0159905384233323  \\
28	0.0143257355996392  \\
29	0.0128566532544362  \\
30	0.0114382267433518  \\
31	0.0101934032439727  \\
32	0.00947467546676925  \\
};

\addplot+[smooth,color=blue,densely dotted, every mark/.append style={solid}, mark=none]
table[row sep=crcr]{
17	0.145044741337596   \\
18	0.087816288758457  \\
19	0.0614424967117903  \\
20	0.0474058725865523  \\
21	0.0368882791357508  \\
22	0.0297061763695038  \\
23	0.0241881933045769  \\
24	0.0210078778136533  \\
25	0.0175981831047232  \\
26	0.0158222741391481  \\
27	0.0135241558975424  \\
28	0.0122817294390961  \\
29	0.0109985327056985  \\
30	0.00994987426904915  \\
31	0.00895736671173737  \\
32	0.00823206224204698 \\
};

\addplot+[smooth,color=blue,loosely dotted, every mark/.append style={solid}, mark=-]
table[row sep=crcr]{
17	0.144182967970758   \\
18	0.0892072257859350  \\
19	0.0608543569085746  \\
20	0.0463601456703347  \\
21	0.0351296716180557  \\
22	0.0281826370504583  \\
23	0.0236185776313252  \\
24	0.0202092236866281  \\
25	0.0173817716482842  \\
26	0.0152383985051633  \\
27	0.0133417767474944  \\
28	0.0118960965620097  \\
29	0.0107537889033601  \\
30	0.00968115871117220  \\
31	0.0087238013185028  \\
32	0.00801984491237327  \\
};

\end{axis}

\end{tikzpicture}%
\captionsetup{justification=centering,font=scriptsize}  
\caption{{\bf Q}-based low-rank approximation error for \texttt{LowRankFastDecay}. Left: Spectral norm. Right: Frobenius norn.} 
\label{figQLRTropFast}      
\end{center}
\end{figure}
\begin{figure}[t]
\begin{center}       
%
%
%
\usetikzlibrary{positioning,calc}

\definecolor{mycolor1}{rgb}{0.00000,1.00000,1.00000}%
\definecolor{mycolor2}{rgb}{1.00000,0.00000,1.00000}%

\pgfplotsset{every axis label/.append style={font=\footnotesize},
every tick label/.append style={font=\footnotesize}
}

\begin{tikzpicture}[font=\footnotesize]

\begin{axis}[%
name=ber,
ymode=log,
width  = 0.34\columnwidth,
height = 0.3\columnwidth,
scale only axis,
xmin  = 11,
xmax  = 20,
xlabel= {$d$},
xmajorgrids,
ymin = 50,
ymax = 281,
xtick       ={12,18},
xticklabels ={$12$,$18$},
ylabel={Magnitude},
ymajorgrids,
]

\addplot+[smooth,color=red,solid, every mark/.append style={solid}, mark=none]
table[row sep=crcr]{
11	280.932041505892 \\
12	186.533745019071 \\
13	177.844082285281 \\
14	161.879587891545 \\
15	154.811301278483 \\
16	147.004929990261 \\
17	147.544264837342 \\
18	144.213202210987 \\
19	137.544364960477 \\
20	137.603258117630  \\
};

\addplot+[smooth,color=red, loosely dashed, every mark/.append style={solid}, mark=none]
table[row sep=crcr]{
11	114.076284818373 \\
12	114.051825603029 \\
13	114.051928610511 \\
14	114.034592805386 \\
15	114.032682442294 \\
16	114.029468512913 \\
17	114.029441545993 \\
18	114.028955966834 \\
19	114.027036881129 \\
20	114.024590562798  \\
};

\addplot+[smooth,color=red, dotted, every mark/.append style={solid}, mark=|]
table[row sep=crcr]{
11	114.013199648510  \\
12	114.013189695028 \\
13	114.013190962154 \\
14	114.013182922429 \\
15	114.013180048384 \\
16	114.013179180973 \\
17	114.013177211763 \\
18	114.013178092361 \\
19	114.013177298282 \\
20	114.013176015167 \\
  };

\addplot+[smooth,color=blue,densely dashed, every mark/.append style={solid}, mark=none]
table[row sep=crcr]{
11	116.35207073483 \\
12	111.563841482772 \\
13	109.655287513116 \\
14	106.39278001465 \\
15	101.255060711609 \\
16	97.123246021736 \\
17	93.7261321440011 \\
18	89.9517960836531 \\
19	86.8751565791444 \\
20	83.5036463945416  \\
};

\addplot+[smooth,color=blue,densely dotted, every mark/.append style={solid}, mark=none]
table[row sep=crcr]{
11	114.011106895776 \\
12	111.266459371817 \\
13	102.430008705082 \\
14	100.691757078088 \\
15	92.9221058590400 \\
16	87.6093581251586 \\
17	82.4155681958241 \\
18	77.7018428801253 \\
19	70.9481290346923 \\
20	69.1004874548555 \\
};

\addplot+[smooth,color=blue,loosely dotted, every mark/.append style={solid}, mark=-]
table[row sep=crcr]{
11	114.010814674907  \\
12	108.927097969962 \\
13	106.309071930687 \\
14	101.50428878320 \\
15	89.757501529436 \\
16	84.1496639314856 \\
17	80.2867077441277 \\
18	68.9521836020017 \\
19	63.673980986834 \\
20	61.0671876337207  \\
};

\addplot+[smooth,color=cyan, loosely dotted, every mark/.append style={solid}, mark=triangle]
table[row sep=crcr]{
11	114.008777902283  \\
12	100.015010421027 \\
13	100.010004991460 \\
14	86.0174518981040 \\
15	86.0116321768124 \\
16	72.0208436921943 \\
17	72.0138925058439 \\
18	58.0258740096420 \\
19	58.0172442883791 \\
20	50.0413761366775 \\
};

\end{axis}

\begin{axis}[%
name=SumRate,
at={($(ber.east)+(35,0em)$)},
		anchor= west,
ymode=log,
width  = 0.34\columnwidth,
height = 0.3\columnwidth,
scale only axis,
xmin  = 11,
xmax  = 20,
xlabel= {$d$},
xmajorgrids,
ymin = 238 ,
ymax = 905,
xtick       ={12,18},
xticklabels ={$12$,$18$},
ylabel={},
ymajorgrids,
]

\addplot+[smooth,color=red,solid, every mark/.append style={solid}, mark=none]
table[row sep=crcr]{
11	904.373444142305 \\
12	600.487450727909 \\
13	572.513780751115 \\
14	521.121049963049 \\
15	498.366896772924 \\
16	473.236709235967 \\
17	474.972930247236 \\
18	464.249608820808 \\
19	442.781358775923 \\
20	442.970947001907 \\
};
\addplot+[smooth,color=red, loosely dashed, every mark/.append style={solid}, mark=none]
table[row sep=crcr]{
11	367.233164444815 \\
12	367.154425598569 \\
13	367.154757199160 \\
14	367.098949959431 \\
15	367.092800138857 \\
16	367.082453891528 \\
17	367.082367079906 \\
18	367.080803908637 \\
19	367.074626008317 \\
20	367.066750846328 \\
};

\addplot+[smooth,color=red, dotted, every mark/.append style={solid}, mark=|]
table[row sep=crcr]{
11	367.030081336038 \\
12	367.030049293897 \\
13	367.030053373015 \\
14	367.030027491616 \\
15	367.030018239523 \\
16	367.030015447163 \\
17	367.030009107903 \\
18	367.030011942715 \\
19	367.030009386423 \\
20	367.030005255833 \\
};

\addplot+[smooth,color=cyan, loosely dotted, every mark/.append style={solid}, mark=triangle]
table[row sep=crcr]{
11	348.872479692082  \\
12	329.718069944203  \\
13	314.183072965766  \\
14	297.840565134805  \\
15	285.149084180149  \\
16	271.867613626666  \\
17	262.154529644107  \\
18	252.069428330925  \\
19	245.299805634853  \\
20	238.340080577479  \\
};

\addplot+[smooth,color=blue,densely dashed, every mark/.append style={solid}, mark=none]
table[row sep=crcr]{
11	355.34788740266  \\
12	344.96701462046  \\
13	331.406823434341   \\
14	323.262329067268  \\
15	311.639173361325  \\
16	302.084540582579  \\
17	291.756667032056  \\
18	283.535192118263  \\
19	275.180266638631  \\
20	267.495648751716  \\
};

\addplot+[smooth,color=blue,densely dotted, every mark/.append style={solid}, mark=none]
table[row sep=crcr]{
11	352.067981323213  \\
12	336.628058408882  \\
13	319.419602608148  \\
14	307.194263301410  \\
15	293.840101574932  \\
16	281.949545205822  \\
17	272.010112112035  \\
18	265.012360819355  \\
19	255.489453376464  \\
20	248.316692181282 \\
};

\addplot+[smooth,color=blue,loosely dotted, every mark/.append style={solid}, mark=-]
table[row sep=crcr]{
11	350.736079114047  \\
12	333.525426597256  \\
13	318.547470678267  \\
14	304.713161849956  \\
15	288.113689791258  \\
16	276.740650284061  \\
17	267.319656809967  \\
18	256.837119320006  \\
19	249.640390944312  \\
20	243.372321473736  \\
};

\end{axis}

\end{tikzpicture}%
\captionsetup{justification=centering,font=scriptsize}  
\caption{{\bf Q}-based low-rank approximation error for \texttt{impcol\_e}. Left: Spectral norm. Right: Frobenius norn.} 
\label{}      
\end{center}
\end{figure}


\begin{figure}[t]
\begin{center}       
%
%
%
\usetikzlibrary{positioning,calc}

\definecolor{mycolor1}{rgb}{0.00000,1.00000,1.00000}%
\definecolor{mycolor2}{rgb}{1.00000,0.00000,1.00000}%

\pgfplotsset{every axis label/.append style={font=\footnotesize},
every tick label/.append style={font=\footnotesize}
}

\begin{tikzpicture}[font=\footnotesize]

\begin{axis}[%
name=ber,
ymode=log,
width  = 0.34\columnwidth,
height = 0.3\columnwidth,
scale only axis,
xmin  = 17,
xmax  = 32,
xlabel= {$d$},
xmajorgrids,
ymin = 0.0 ,
ymax =  0.73,
xtick       ={20,30},
xticklabels ={$20$,$30$},
ylabel={Magnitude},
ymajorgrids,
]

\addplot+[smooth,color=red,solid, every mark/.append style={solid}, mark=none]
table[row sep=crcr]{
17	0.726585891715203  \\
18	0.344149500998731 \\
19	0.298786733217387 \\
20	0.212432796942439 \\
21	0.157313290496056 \\
22	0.148454478334997 \\
23	0.158690438613206 \\
24	0.127778985304701 \\
25	0.122492007879723 \\
26	0.117673790363809 \\
27	0.0917183910357720 \\
28	0.097626037866027 \\
29	0.0860594526159401 \\
30	0.0903476152395745 \\
31	0.0771465978567405 \\
32	0.0747531943825955 \\
};

\addplot+[smooth,color=red, loosely dashed, every mark/.append style={solid}, mark=none]
table[row sep=crcr]{
17	0.00484596094790707  \\
18	0.00484320715992987 \\
19	0.00484058117208249 \\
20	0.00484097493349677 \\
21	0.00483909575139813 \\
22	0.00483883667593196 \\
23	0.00483819522381632 \\
24	0.00483798849754648 \\
25	0.00483792242146004 \\
26	0.00483795927556527 \\
27	0.00483760424911735 \\
28	0.00483726485531308 \\
29	0.00483736352119273 \\
30	0.00483702633481383 \\
31	0.00483694360547449 \\
32	0.00483693060213907  \\
};

\addplot+[smooth,color=red, dotted, every mark/.append style={solid}, mark=|]
table[row sep=crcr]{
17	0.0048351637052893 \\
18	0.00483516373914916 \\
19	0.0048351636800442 \\
20	0.0048351636183238 \\
21	0.00483516355597924 \\
22	0.00483516355213937 \\
23	0.00483516353212802 \\
24	0.00483516351670264 \\
25	0.00483516352463536 \\
26	0.00483516350783414 \\
27	0.00483516351262422 \\
28	0.00483516350673206 \\
29	0.00483516350525001 \\
30	0.00483516350191618 \\
31	0.00483516349883886 \\
32	0.00483516349652004 \\
  };

\addplot+[smooth,color=blue,densely dashed, every mark/.append style={solid}, mark=none]
table[row sep=crcr]{
17	0.278287582112130  \\
18	0.165788292294573 \\
19	0.145655826615989 \\
20	0.104871080035947 \\
21	0.076727377159349 \\
22	0.0718135766519517 \\
23	0.0760150406905549 \\
24	0.0621617250148673 \\
25	0.0579151133017745 \\
26	0.0563937406999205 \\
27	0.0436565146271251 \\
28	0.0457003305683698 \\
29	0.0401517667427609 \\
30	0.0427196973612829 \\
31	0.0355670638096751 \\
32	0.0347613116615809  \\
};

\addplot+[smooth,color=blue,densely dotted, every mark/.append style={solid}, mark=none]
table[row sep=crcr]{
17	0.00483035586773119 \\
18	0.00481773629474607 \\
19	0.00482249749652058 \\
20	0.00480750379597512 \\
21	0.00481277149024017 \\
22	0.00480513954662048 \\
23	0.00479250903301804 \\
24	0.0047865346738800 \\
25	0.00478387408792891 \\
26	0.00478235799171736 \\
27	0.00478398254807529 \\
28	0.00476784636696265 \\
29	0.00476061063100194 \\
30	0.00476958702230991 \\
31	0.00475555802964814 \\
32	0.00475765855748010 \\
};

\addplot+[smooth,color=blue,loosely dotted, every mark/.append style={solid}, mark=-]
table[row sep=crcr]{
17	0.00482675093924692 \\
18	0.00481795808719408 \\
19	0.00481799402599406 \\
20	0.00481048300348517 \\
21	0.00480508410754269 \\
22	0.00479821523671610 \\
23	0.00478645258033380 \\
24	0.00479294899711563 \\
25	0.00477426341960964 \\
26	0.00477102211562875 \\
27	0.00476435691283832 \\
28	0.00475064900947291 \\
29	0.00474838270875350 \\
30	0.00474643559711136 \\
31	0.00473496580478843 \\
32	0.00473670393205367 \\
};

\addplot+[smooth,color=cyan, loosely dotted, every mark/.append style={solid}, mark=triangle]
table[row sep=crcr]{
17	0.00478625302213134  \\
18	0.00477709128109251 \\
19	0.00472867300418597 \\
20	0.00471389879573033 \\
21	0.00470841216203010 \\
22	0.00467746585214326 \\
23	0.00464465997466443 \\
24	0.00462437768146376 \\
25	0.00461257157911528 \\
26	0.00459677702545871 \\
27	0.00458528212821993 \\
28	0.00457755514204000 \\
29	0.00457151067442002 \\
30	0.00456022026290263 \\
31	0.00454629071614261 \\
32	0.00451858306243967 \\
};

\end{axis}

\begin{axis}[%
name=SumRate,
at={($(ber.east)+(35,0em)$)},
		anchor= west,
ymode=log,
width  = 0.34\columnwidth,
height = 0.3\columnwidth,
scale only axis,
xmin  = 17,
xmax  = 32,
xlabel= {$d$},
xmajorgrids,
ymin = 0.0 ,
ymax = 11,
xtick       ={20,30},
xticklabels ={$20$,$30$},
ylabel={},
ymajorgrids,
]

\addplot+[smooth,color=red,solid, every mark/.append style={solid}, mark=none]
table[row sep=crcr]{
17	10.1468261308984   \\
18	4.80607343121702  \\
19	4.17257899822295  \\
20	2.96663984210729  \\
21	2.19689182647749  \\
22	2.07317785439312  \\
23	2.21612380257357  \\
24	1.78444305326207  \\
25	1.71061001947910  \\
26	1.64332325276331  \\
27	1.28085416666817  \\
28	1.36335489495489  \\
29	1.20182666986949  \\
30	1.26171118050913  \\
31	1.07735798887424  \\
32	1.04393393097530 \\
};
\addplot+[smooth,color=red, loosely dashed, every mark/.append style={solid}, mark=none]
table[row sep=crcr]{
17	0.0676742058113206  \\
18	0.067635748956976  \\
19	0.0675990768409709  \\
20	0.0676045757484688  \\
21	0.0675783328304041  \\
22	0.0675747148222099  \\
23	0.0675657568956070  \\
24	0.0675628699478412  \\
25	0.0675619471903694  \\
26	0.0675624618606957  \\
27	0.067557503889890  \\
28	0.0675527642301180  \\
29	0.0675541421064810  \\
30	0.0675494332735690  \\
31	0.0675482779521792  \\
32	0.0675480963596309 \\
};

\addplot+[smooth,color=red, dotted, every mark/.append style={solid}, mark=|]
table[row sep=crcr]{
17	0.0675234215134359  \\
18	0.0675234219862909  \\
19	0.0675234211608860  \\
20	0.0675234202989564  \\
21	0.0675234194283092  \\
22	0.0675234193746852  \\
23	0.0675234190952252  \\
24	0.0675234188798086  \\
25	0.0675234189905897  \\
26	0.0675234187559594  \\
27	0.0675234188228532  \\
28	0.0675234187405688  \\
29	0.0675234187198718  \\
30	0.0675234186733146  \\
31	0.0675234186303396  \\
32	0.0675234185979571  \\
};

\addplot+[smooth,color=cyan, loosely dotted, every mark/.append style={solid}, mark=triangle]
table[row sep=crcr]{
17	0.0673500792420047\\
18	0.0671797957418185\\
19	0.0670097332840876\\
20	0.0668426810236098\\
21	0.0666762563778716\\
22	0.0665098039350599\\
23	0.0663451229004977\\
24	0.0661823425575273\\
25	0.0660205846494945\\
26	0.0658592573666649\\
27	0.0656986409438340\\
28	0.0655384361247005\\
29	0.0653783802077800\\
30	0.0652183554587717\\
31	0.0650587294673091\\
32	0.0648996881396577  \\
};

\addplot+[smooth,color=blue,densely dashed, every mark/.append style={solid}, mark=none]
table[row sep=crcr]{
17	0.332656815604842 \\
18	0.220767782041101\\
19	0.196364259060752\\
20	0.162023273954572\\
21	0.140254752767000\\
22	0.133298511420245\\
23	0.134076660281047\\
24	0.121728664504140\\
25	0.117684159523836\\
26	0.113868904592474\\
27	0.104194345507563\\
28	0.105200453496978\\
29	0.101237323556455\\
30	0.102014072171892\\
31	0.0969083125453596\\
32	0.0949823636675624 \\
};

\addplot+[smooth,color=blue,densely dotted, every mark/.append style={solid}, mark=none]
table[row sep=crcr]{
17	0.0674005152131451\\
18	0.0672798635135381\\
19	0.0671557252827123\\
20	0.067038484640030\\
21	0.0669112643435437\\
22	0.0667873596669697\\
23	0.0666738256599205\\
24	0.0665469525706769\\
25	0.0664318118964810\\
26	0.0663085022886657\\
27	0.0661861516076614\\
28	0.0660627990166538\\
29	0.0659423384210134\\
30	0.0658291538079711\\
31	0.0656929262580196\\
32	0.0655764524049499 \\
};

\addplot+[smooth,color=blue,loosely dotted, every mark/.append style={solid}, mark=-]
table[row sep=crcr]{
17	0.0673873602329568\\
18	0.0672488720479228\\
19	0.0671121992805588\\
20	0.0669769112108705\\
21	0.0668450928718736\\
22	0.0667007027220277\\
23	0.0665665865248233\\
24	0.0664271188846929\\
25	0.0663002576592798\\
26	0.0661584342998437\\
27	0.0660244018426649\\
28	0.0658841899809893\\
29	0.0657513937042568\\
30	0.0656143962853200\\
31	0.0654781415131021\\
32	0.0653460388399324 \\
};

\end{axis}

\end{tikzpicture}%
\captionsetup{justification=centering,font=scriptsize}  
\caption{{\bf P}-based low-rank approximation error for \texttt{LowRankLargeGap}. Left: Spectral norm. Right: Frobenius norn.} 
\label{}      
\end{center}
\end{figure}
\begin{figure}[t]
\begin{center}       
%
%
%
\usetikzlibrary{positioning,calc}

\definecolor{mycolor1}{rgb}{0.00000,1.00000,1.00000}%
\definecolor{mycolor2}{rgb}{1.00000,0.00000,1.00000}%

\pgfplotsset{every axis label/.append style={font=\footnotesize},
every tick label/.append style={font=\footnotesize}
}

\begin{tikzpicture}[font=\footnotesize]

\begin{axis}[%
name=ber,
ymode=log,
width  = 0.34\columnwidth,
height = 0.3\columnwidth,
scale only axis,
xmin  = 17,
xmax  = 32,
xlabel= {$d$},
xmajorgrids,
ymin = 0.00 ,
ymax =  1.32,
xtick       ={20,30},
xticklabels ={$20$,$30$},
ylabel={Magnitude},
ymajorgrids,
]

\addplot+[smooth,color=red,solid, every mark/.append style={solid}, mark=none]
table[row sep=crcr]{
17	1.31703785353809 \\
18	1.06257330966435 \\
19	0.559328143194871 \\
20	0.461237950177358 \\
21	0.422386286220590 \\
22	0.328509086225262 \\
23	0.268039159930106 \\
24	0.279548587170675 \\
25	0.220441914595202 \\
26	0.208331823834586 \\
27	0.188045373047803 \\
28	0.191907895136685 \\
29	0.183002790176545 \\
30	0.172749887005170 \\
31	0.153144992405919 \\
32	0.154333675131308 \\
};

\addplot+[smooth,color=red, loosely dashed, every mark/.append style={solid}, mark=none]
table[row sep=crcr]{
17	0.00978377044367659  \\
18	0.00967726010789597 \\
19	0.00966738749275604 \\
20	0.00966982270640446 \\
21	0.00965016105529496 \\
22	0.00965192269220248 \\
23	0.00964550731428026 \\
24	0.00964472400135178 \\
25	0.00964276321795158 \\
26	0.00963956559469626 \\
27	0.00963823742337222 \\
28	0.00963719680357879 \\
29	0.00963621653808982 \\
30	0.00963500175456700 \\
31	0.00963481249475699 \\
32	0.00963463675739822  \\
};

\addplot+[smooth,color=red, dotted, every mark/.append style={solid}, mark=|]
table[row sep=crcr]{
17	0.00962095950543359  \\
18	0.00962095790865574 \\
19	0.00962095792703001 \\
20	0.00962095803535516 \\
21	0.00962095597054906 \\
22	0.00962095545978764 \\
23	0.00962095546657317 \\
24	0.00962095505514117 \\
25	0.00962095499524296 \\
26	0.00962095467082961 \\
27	0.00962095468008893 \\
28	0.00962095465972456 \\
29	0.00962095449623256 \\
30	0.00962095435344819 \\
31	0.00962095436896867 \\
32	0.00962095423891173 \\
  };

\addplot+[smooth,color=blue,densely dashed, every mark/.append style={solid}, mark=none]
table[row sep=crcr]{
17	0.463745878040598  \\
18	0.412161857709167 \\
19	0.264533369495758 \\
20	0.219759113256419 \\
21	0.200984687369907 \\
22	0.158218780745991 \\
23	0.128256089588206 \\
24	0.134046854236248 \\
25	0.104587306497122 \\
26	0.0989489259370805 \\
27	0.0897220241931546 \\
28	0.0906164547902027 \\
29	0.085948598507038 \\
30	0.0823777465153792 \\
31	0.0721388669694697 \\
32	0.0712343200378385  \\
};

\addplot+[smooth,color=blue,densely dotted, every mark/.append style={solid}, mark=none]
table[row sep=crcr]{
17	0.00961670549834992 \\
18	0.00961141900386619 \\
19	0.00959210159739829 \\
20	0.00959540863464773 \\
21	0.00958862150012764 \\
22	0.00956783285395827 \\
23	0.00954840324307653 \\
24	0.00956172782116675 \\
25	0.00955390315208731 \\
26	0.00952847415640021 \\
27	0.00951677079205414 \\
28	0.00952409234177322 \\
29	0.00950900753414741 \\
30	0.00951766214898294 \\
31	0.00949325103512509 \\
32	0.00947454828026843 \\
};

\addplot+[smooth,color=blue,loosely dotted, every mark/.append style={solid}, mark=-]
table[row sep=crcr]{
17	0.00961150985897984   \\
18	0.00960456402525665 \\
19	0.00959843159337552 \\
20	0.00958004295269987 \\
21	0.00957170333836353 \\
22	0.00955250574822476 \\
23	0.00955326066309025 \\
24	0.00953699256144073 \\
25	0.00952179623555695 \\
26	0.00953057055069570 \\
27	0.00950886925731551 \\
28	0.00949895094207853 \\
29	0.00948335432960505 \\
30	0.00948150029979015 \\
31	0.00945098053213761 \\
32	0.00944159410072621 \\
};

\addplot+[smooth,color=cyan,loosely dotted, every mark/.append style={solid}, mark=triangle]
table[row sep=crcr]{
17	0.00959605132682474  \\
18	0.00952930555326711 \\
19	0.00943545513400279 \\
20	0.00940354310830665 \\
21	0.00936711788345268 \\
22	0.00931371096921667 \\
23	0.00929876746095547 \\
24	0.00926447065288894 \\
25	0.00925918487863594 \\
26	0.00922224624370206 \\
27	0.00921187951784446 \\
28	0.00914033436668026 \\
29	0.00910696069548622 \\
30	0.00909995623284013 \\
31	0.00907862729681516 \\
32	0.00903410655133817 \\
};

\end{axis}

\begin{axis}[%
name=SumRate,
at={($(ber.east)+(35,0em)$)},
		anchor= west,
ymode=log,
width  = 0.34\columnwidth,
height = 0.3\columnwidth,
scale only axis,
xmin  = 17,
xmax  = 32,
xlabel= {$d$},
xmajorgrids,
ymin = 0.1 ,
ymax =  18.5,
xtick       ={20,30},
xticklabels ={$20$,$30$},
ylabel={},
ymajorgrids,
legend entries = {Th. bound q=0, Th. bound q=1, Th. bound q=2, SVD, Computed q=0, Computed q=1, Computed q=2},
legend style={at={(0.97,1.3)},anchor=north east,draw=black,fill=white,legend cell align=left,font=\tiny, legend columns=4}
]

\addplot+[smooth,color=red,solid, every mark/.append style={solid}, mark=none]
table[row sep=crcr]{
17	18.4699807111559   \\
18	14.901400503385  \\
19	7.84395071733384  \\
20	6.46834562889997  \\
21	5.92349455878828  \\
22	4.60697197861106  \\
23	3.75894899333683  \\
24	3.92035581893313  \\
25	3.09145094012737  \\
26	2.92162048145221  \\
27	2.63712573156850  \\
28	2.69129327754025  \\
29	2.56640915488367  \\
30	2.42262367195324  \\
31	2.14768698420379  \\
32	2.16435692801036 \\
};
\addplot+[smooth,color=red, loosely dashed, every mark/.append style={solid}, mark=none]
table[row sep=crcr]{
17	0.137206421889572   \\
18	0.135712733730100  \\
19	0.135574281360858  \\
20	0.135608432504649  \\
21	0.135332699870415  \\
22	0.135357404854873  \\
23	0.135267436365235  \\
24	0.135256451278785  \\
25	0.135228953487827  \\
26	0.135184110403257  \\
27	0.135165484288093  \\
28	0.135150890761065  \\
29	0.135137143635556  \\
30	0.135120107657300  \\
31	0.135117453500450  \\
32	0.135114988980833 \\
};

\addplot+[smooth,color=red, dotted, every mark/.append style={solid}, mark=|]
table[row sep=crcr]{
17	0.134923180841614   \\
18	0.134923158448593  \\
19	0.134923158706272  \\
20	0.134923160225410  \\
21	0.134923131268818  \\
22	0.134923124105961  \\
23	0.134923124201120  \\
24	0.134923118431247  \\
25	0.134923117591242  \\
26	0.134923113041708  \\
27	0.134923113171560  \\
28	0.134923112885972  \\  
29	0.134923110593180  \\
30	0.134923108590789  \\
31	0.134923108808446  \\
32	0.134923106984544  \\
};

\addplot+[smooth,color=cyan,loosely dotted, every mark/.append style={solid}, mark=triangle]
table[row sep=crcr]{
17	0.134579631403724   \\
18	0.134237077544526  \\
19	0.133898414192878  \\
20	0.133565555102286  \\
21	0.133234120580248  \\
22	0.132904431789724  \\
23	0.132577685819792  \\
24	0.132251184119608  \\
25	0.131926287314395  \\
26	0.131600960406615  \\
27	0.131277427435809  \\
28	0.130953824800625  \\
29	0.130634446137219  \\
30	0.130316621290106  \\
31	0.129998509918496  \\
32	0.129681113148503  \\
};

\addplot+[smooth,color=blue,densely dashed, every mark/.append style={solid}, mark=none]
table[row sep=crcr]{
17	0.569549769587862   \\
18	0.512686357198831  \\
19	0.384221466268029  \\
20	0.333337235947760  \\
21	0.309331948185538  \\
22	0.273325390634501  \\
23	0.248651405510745  \\
24	0.252624794378669  \\
25	0.225325997312638  \\
26	0.220005385169547  \\
27	0.212364231648124  \\
28	0.212144510651327  \\
29	0.206636910212617  \\
30	0.201742865875292  \\
31	0.192258237494843  \\
32	0.190582717253089  \\
};

\addplot+[smooth,color=blue,densely dotted, every mark/.append style={solid}, mark=none]
table[row sep=crcr]{
17	0.134681871481681   \\
18	0.13443861327483  \\
19	0.134189797332209   \\
20	0.133944749972273  \\
21	0.133707193994607  \\
22	0.133464416225649  \\
23	0.133218480416412  \\
24	0.132971833785376  \\
25	0.132734378952566  \\
26	0.132486681891406  \\
27	0.132240650133733  \\
28	0.132007244688905  \\
29	0.131765803383717  \\
30	0.131520332280071  \\
31	0.131264922450512  \\
32	0.131042319751829  \\
};

\addplot+[smooth,color=blue,loosely dotted, every mark/.append style={solid}, mark=-]
table[row sep=crcr]{
17	0.134650155964546   \\
18	0.134376388250014  \\
19	0.134111107055330  \\
20	0.133833227091302  \\
21	0.133557961032621  \\
22	0.133286050041275  \\
23	0.133019568616366  \\
24	0.132748551789638  \\
25	0.132471867635866  \\  
26	0.132200931103072  \\
27	0.131935329694194  \\
28	0.131657081632950  \\
29	0.131391145828111  \\
30	0.131117867096016  \\
31	0.130845647506311   \\
32	0.130587520365342 \\
};

\end{axis}

\end{tikzpicture}%
\captionsetup{justification=centering,font=scriptsize}  
\caption{{\bf P}-based low-rank approximation error for \texttt{LowRankMediumGap}. Left: Spectral norm. Right: Frobenius norn.} 
\label{}      
\end{center}
\end{figure}
 
\begin{figure}[t]
\begin{center}       
%
%
%
\usetikzlibrary{positioning,calc}

\definecolor{mycolor1}{rgb}{0.00000,1.00000,1.00000}%
\definecolor{mycolor2}{rgb}{1.00000,0.00000,1.00000}%

\pgfplotsset{every axis label/.append style={font=\footnotesize},
every tick label/.append style={font=\footnotesize}
}

\begin{tikzpicture}[font=\footnotesize]

\begin{axis}[%
name=ber,
ymode=log,
width  = 0.34\columnwidth,
height = 0.3\columnwidth,
scale only axis,
xmin  = 17,
xmax  = 32,
xlabel= {$d$},
xmajorgrids,
ymin = 0.05 ,
ymax = 75.8,
xtick       ={20,30},
xticklabels ={$20$,$30$},
ylabel={Magnitude},
ymajorgrids,
]

\addplot+[smooth,color=red,solid, every mark/.append style={solid}, mark=none]
table[row sep=crcr]{
17	75.7015169203566  \\
18	36.5984532154262 \\
19	21.9694972379885 \\
20	22.7548999484924 \\
21	15.8860336540673 \\
22	17.3415724404856 \\
23	12.6971806719315 \\
24	12.7023905789882 \\
25	12.3552371125110 \\
26	10.6427961895941 \\
27	9.6668057621410 \\
28	10.619228894209 \\
29	9.44002866270154 \\
30	8.21437533416470 \\
31	8.56015719166760 \\
32	8.02874691046324 \\
};

\addplot+[smooth,color=red, loosely dashed, every mark/.append style={solid}, mark=none]
table[row sep=crcr]{
17	10.8026116644772  \\
18	11.9906746220733 \\
19	7.53051372810077 \\
20	5.19546025295922 \\
21	4.76774083591032 \\
22	5.2842422442982 \\
23	4.38332216769299 \\
24	3.91877994994283 \\
25	3.32802279503889 \\
26	3.14970209838505 \\
27	3.06426136727027 \\
28	2.68786986026122 \\
29	2.62981093261903 \\
30	2.58664663341437 \\
31	2.39941678791722 \\
32	2.39220437812664 \\
};

\addplot+[smooth,color=red, dotted, every mark/.append style={solid}, mark=|]
table[row sep=crcr]{
17	5.17983640514805  \\
18	3.27960683816622 \\
19	2.32938039110366 \\
20	1.88346546568638 \\
21	1.60492758964131 \\
22	1.61383058169618 \\
23	1.4806389868147 \\
24	1.18270707716168 \\
25	1.21118595742013 \\
26	1.22437078893274 \\
27	1.10419111377375 \\
28	1.09540833525035 \\
29	1.07420303596249 \\
30	0.993285667706181 \\
31	1.00194380168089 \\
32	0.973100009243637 \\
  };

\addplot+[smooth,color=blue,densely dashed, every mark/.append style={solid}, mark=none]
table[row sep=crcr]{
17	0.945230659678021   \\
18	0.859275535692054  \\
19	0.712553499605930  \\
20	0.724510447320748  \\
21	0.617556092973322  \\
22	0.542297553714895  \\
23	0.457781557928060  \\
24	0.427316719932490  \\
25	0.401282066606413  \\
26	0.378491882811200  \\
27	0.330306229062682  \\
28	0.320079118734565  \\
29	0.316846842542562  \\
30	0.259680693342337  \\
31	0.266310538887220  \\
32	0.239059676459357  \\
};

\addplot+[smooth,color=blue,densely dotted, every mark/.append style={solid}, mark=none]
table[row sep=crcr]{
17	0.501579429327556    \\
18	0.324069583577585  \\
19	0.254645163888528  \\
20	0.206531464253582  \\
21	0.178548546970126  \\
22	0.149812183884663  \\
23	0.146539849331434  \\
24	0.122202972622741  \\
25	0.113363663938561  \\
26	0.102037648444218  \\
27	0.0956329136580108  \\
28	0.0872866683798577  \\
29	0.0865093905767792  \\
30	0.0808767401953916  \\
31	0.0766209238439626  \\
32	0.0694737222426055  \\
};

\addplot+[smooth,color=blue,loosely dotted, every mark/.append style={solid}, mark=-]
table[row sep=crcr]{
17	0.434643127880937   \\
18	0.263179034976335  \\
19	0.219231209695323  \\
20	0.176517806859702  \\
21	0.156884645268550  \\
22	0.138049733360362  \\
23	0.125202062007892  \\
24	0.117647569853330  \\
25	0.0987167817900514  \\
26	0.092868202353837  \\
27	0.086897303093059  \\
28	0.0777733221382892  \\
29	0.0728995219351086  \\
30	0.0706864304875462  \\
31	0.0676742472627202  \\
32	0.0646645712166681 \\
};

\addplot+[smooth,color=cyan,loosely dotted, every mark/.append style={solid}, mark=triangle]
table[row sep=crcr]{
17	0.333333333333334   \\
18	0.250000000000000  \\
19	0.200000000000000  \\
20	0.166666666666667  \\
21	0.142857142857143  \\
22	0.125000000000000  \\
23	0.111111111111111  \\
24	0.100000000000000  \\
25	0.0909090909090910  \\
26	0.0833333333333335  \\
27	0.0769230769230770  \\
28	0.0714285714285715  \\
29	0.0666666666666667  \\
30	0.0625000000000001  \\
31	0.0588235294117647  \\
32	0.0555555555555556  \\
};

\end{axis}

\begin{axis}[%
name=SumRate,
at={($(ber.east)+(35,0em)$)},
		anchor= west,
ymode=log,
width  = 0.34\columnwidth,
height = 0.3\columnwidth,
scale only axis,
xmin  = 17,
xmax  = 32,
xlabel= {$d$},
xmajorgrids,
ymin = 0.2 ,
ymax = 121.5,
xtick       ={20,30},
xticklabels ={$20$,$30$},
ylabel={},
ymajorgrids,
]

\addplot+[smooth,color=red,solid, every mark/.append style={solid}, mark=none]
table[row sep=crcr]{
17	121.468361427277   \\
18	58.7247696440179  \\
19	35.2515899210725  \\
20	36.5118233243980  \\
21	25.490248492223  \\
22	27.8257619478662  \\
23	20.3735115716139  \\
24	20.3818712306947  \\
25	19.8248392919405  \\
26	17.0771084483622  \\
27	15.5110637663758  \\
28	17.0392930799245  \\
29	15.1471841005677  \\
30	13.180537888552  \\
31	13.7353689850854  \\
32	12.8826841416441  \\
};
\addplot+[smooth,color=red, loosely dashed, every mark/.append style={solid}, mark=none]
table[row sep=crcr]{
17	17.3335435193423  \\
18	19.2398733605722  \\
19	12.0832342662349  \\
20	8.33647817720015  \\
21	7.65017255410122  \\
22	8.47893507173692  \\
23	7.03334600121270  \\
24	6.28795562728836  \\
25	5.34004458763090  \\
26	5.05391659822877  \\  
27	4.91682098865757  \\
28	4.31287457554102  \\
29	4.21971497856303  \\
30	4.15045485129158  \\
31	3.85003151146936  \\
32	3.83845869714748 \\
};

\addplot+[smooth,color=red, dotted, every mark/.append style={solid}, mark=|]
table[row sep=crcr]{
17	8.31140862417112   \\
18	5.26235780951199  \\
19	3.7376532301970  \\
20	3.02215164542187  \\
21	2.57521820504942  \\
22	2.58950367647322  \\
23	2.3757884770386  \\
24	1.89773595768781   \\
25	1.94343230646662  \\
26	1.96458828780867  \\
27	1.77175162069428  \\
28	1.75765904026238  \\
29	1.72363365922833  \\
30	1.59379610070963  \\
31	1.60768868027355  \\
32	1.56140680446401  \\
};

\addplot+[smooth,color=blue,densely dashed, every mark/.append style={solid}, mark=none]
table[row sep=crcr]{
17	1.57465620385962   \\
18	1.37231595173353  \\
19	1.19640508017073  \\
20	1.16379374116945  \\
21	1.03851735679864  \\
22	0.943552765986809  \\
23	0.849151630896603  \\
24	0.785715682682583  \\
25	0.750372919358458  \\
26	0.706330187951081  \\
27	0.656747391201188  \\
28	0.616480692417111  \\
29	0.605203264620506  \\  
30	0.554897410593399  \\
31	0.553205570558690  \\
32	0.519619720008773  \\
};

\addplot+[smooth,color=blue,densely dotted, every mark/.append style={solid}, mark=none]
table[row sep=crcr]{
17	0.749731790389086   \\
18	0.584117089325138  \\
19	0.503309071731319  \\
20	0.452223587701031  \\
21	0.414989040811839  \\
22	0.383195822796104  \\
23	0.365943264409190  \\
24	0.339089030693477  \\
25	0.322742268664946  \\
26	0.307108677100304  \\
27	0.295618995516873  \\
28	0.282547673955071  \\
29	0.274476375290682  \\
30	0.266255031889762  \\
31	0.257325545063777  \\
32	0.248430297721527 \\
};

\addplot+[smooth,color=blue,loosely dotted, every mark/.append style={solid}, mark=-]
table[row sep=crcr]{
17	0.694174387348973   \\
18	0.539370599015114  \\
19	0.479033620842475  \\
20	0.429850301869957  \\
21	0.398070663375412  \\
22	0.370473064721922  \\
23	0.348521336691303  \\
24	0.330890748192548  \\
25	0.311706610569134  \\
26	0.298552069405415  \\
27	0.286755908192277  \\
28	0.274781691534938  \\
29	0.264291172383017  \\
30	0.256908985101578  \\
31	0.248435502709075  \\
32	0.241372473078274  \\
};

\addplot+[smooth,color=cyan,loosely dotted, every mark/.append style={solid}, mark=triangle]
table[row sep=crcr]{
17	0.627424093054539  \\
18	0.531554213071630  \\
19	0.469094746756132  \\
20	0.424322850473787  \\
21	0.390220583332584  \\
22	0.363130748341580  \\
23	0.340938323441522  \\
24	0.322324776240937  \\
25	0.306420073393977  \\
26	0.292623988368782  \\
27	0.280507315634385  \\
28	0.269753951520816  \\
29	0.260125265102735  \\
30	0.251437286615020  \\
31	0.243545599632438  \\
32	0.236335040753307  \\
};

\end{axis}

\end{tikzpicture}%
\captionsetup{justification=centering,font=scriptsize}  
\caption{{\bf P}-based low-rank approximation error for \texttt{LowRankSlowDecay}. Left: Spectral norm. Right: Frobenius norn.} 
\label{}      
\end{center}
\end{figure}
\begin{figure}[t]
\begin{center}       
%
%
%
\usetikzlibrary{positioning,calc}

\definecolor{mycolor1}{rgb}{0.00000,1.00000,1.00000}%
\definecolor{mycolor2}{rgb}{1.00000,0.00000,1.00000}%

\pgfplotsset{every axis label/.append style={font=\footnotesize},
every tick label/.append style={font=\footnotesize}
}

\begin{tikzpicture}[font=\footnotesize]

\begin{axis}[%
name=ber,
ymode=log,
width  = 0.35\columnwidth,
height = 0.3\columnwidth,
scale only axis,
xmin  = 17,
xmax  = 32,
xlabel= {$d$},
xmajorgrids,
ymin = 0.003 ,
ymax =30.6,
xtick       ={20,30},
xticklabels ={$20$,$30$},
ylabel={Magnitude},
ymajorgrids,
]

\addplot+[smooth,color=red,solid, every mark/.append style={solid}, mark=none]
table[row sep=crcr]{
17	30.5965679948811  \\
18	20.0414833243133 \\
19	14.3619584874612 \\
20	9.59033922363657 \\
21	11.0828948937010 \\
22	7.00931046540275 \\
23	7.07268751358384 \\
24	5.88192011384219 \\
25	5.87157624144315 \\
26	5.49954678663893 \\
27	4.88860505327535 \\
28	4.96431299294603 \\
29	4.34608375540750 \\
30	4.31433511344509 \\
31	4.22089033258075 \\
32	4.08457340203634  \\
};

\addplot+[smooth,color=red, loosely dashed, every mark/.append style={solid}, mark=none]
table[row sep=crcr]{
17	1.67261313294223 \\
18	1.51278039000494 \\
19	1.46635059860087 \\
20	0.86920921379983 \\
21	0.814974801490419 \\
22	0.803644490220768 \\
23	0.661096289558895 \\
24	0.664512990563888 \\
25	0.612258928372031 \\
26	0.599927631435108 \\
27	0.569798592503786 \\
28	0.570335062314899 \\
29	0.503690207660317 \\
30	0.491183551652858 \\
31	0.474975620688575 \\
32	0.485505954128998  \\
};

\addplot+[smooth,color=red, dotted, every mark/.append style={solid}, mark=|]
table[row sep=crcr]{
17	0.353257609243787  \\
18	0.315399731096371 \\
19	0.312594390133840 \\
20	0.291085134904953 \\
21	0.287515540600979 \\
22	0.275786558797310 \\
23	0.277426468137945 \\
24	0.273789936447863 \\
25	0.270354253490037 \\
26	0.271753726311520 \\
27	0.267297173265130 \\
28	0.268297619455700 \\
29	0.265942522447107 \\
30	0.265814008579200 \\
31	0.265093151737111 \\
32	0.263605889010095 \\
  };

\addplot+[smooth,color=blue,densely dashed, every mark/.append style={solid}, mark=none]
table[row sep=crcr]{
17	0.575804310631972   \\
18	0.329634382205102  \\
19	0.248129776050919  \\
20	0.120147155667957  \\
21	0.116083154443244  \\
22	0.081384037821754  \\
23	0.0585815653016246  \\
24	0.048889504324760  \\
25	0.0411727413428059  \\
26	0.0371557311894364  \\
27	0.0298925119221863  \\
28	0.0243165178727163  \\
29	0.022273300147498  \\
30	0.0185420715906860  \\
31	0.0160061763524710  \\
32	0.0152618285445664  \\
};

\addplot+[smooth,color=blue,densely dotted, every mark/.append style={solid}, mark=none]
table[row sep=crcr]{
17	0.120318779073906   \\
18	0.0669430997012926  \\
19	0.0452354460272852  \\
20	0.0308768312647593  \\
21	0.0243650619383768  \\
22	0.0184388211608594  \\
23	0.0136859125556790  \\
24	0.0119578540121250  \\
25	0.0101082437277676  \\
26	0.0085276985836819  \\
27	0.00716886555650309  \\
28	0.00602647586904796  \\
29	0.00531931277652040  \\
30	0.00482644645566057  \\
31	0.00404062932688560  \\
32	0.00378765111940000  \\
};

\addplot+[smooth,color=blue,loosely dotted, every mark/.append style={solid}, mark=-]
table[row sep=crcr]{
17	0.113507858818005  \\
18	0.0632945636831570  \\
19	0.0426315173204687  \\
20	0.0299871151330119  \\
21	0.0214574531008389  \\
22	0.0169383778207443  \\
23	0.0128923503215717  \\
24	0.0108875823127112  \\
25	0.00920040447244673  \\
26	0.00791689654863835  \\
27	0.00689201401219045  \\
28	0.00573452380549204  \\
29	0.00507123329331521  \\
30	0.00424227997488127  \\
31	0.00397872921654173  \\
32	0.00345911053078922   \\
};

\addplot+[smooth,color=cyan, loosely dotted, every mark/.append style={solid}, mark=triangle]
table[row sep=crcr]{
17	0.111111111111111   \\
18	0.0624999999999999  \\
19	0.04000000000000  \\
20	0.0277777777777778  \\
21	0.0204081632653061  \\
22	0.0156250000000000  \\
23	0.0123456790123457  \\
24	0.0100000000000000  \\
25	0.00826446280991737  \\
26	0.00694444444444445  \\
27	0.00591715976331360  \\
28	0.00510204081632653  \\
29	0.00444444444444444  \\
30	0.00390625000000000  \\
31	0.00346020761245674  \\
32	0.00308641975308643  \\
};

\end{axis}

\begin{axis}[%
name=SumRate,
at={($(ber.east)+(35,0em)$)},
		anchor= west,
ymode=log,
width  = 0.35\columnwidth,
height = 0.3\columnwidth,
scale only axis,
xmin  = 17,
xmax  = 32,
xlabel= {$d$},
xmajorgrids,
ymin = 0.0077 ,
ymax = 35.2,
xtick       ={20,30},
xticklabels ={$20$,$30$},
ylabel={},
ymajorgrids,
]

\addplot+[smooth,color=red,solid, every mark/.append style={solid}, mark=none]
table[row sep=crcr]{
17	35.1151004275908 \\
18	23.0012300650481\\
19	16.4829471955307\\
20	11.0066503219904\\
21	12.7196281388766\\
22	8.04445259879993\\
23	8.11718923423093\\
24	6.75056808787148\\
25	6.73869662182460\\
26	6.31172547690215\\
27	5.61055924394237\\
28	5.69744780952093\\
29	4.98791784632145\\
30	4.95148053246533\\
31	4.84423573085769\\
32	4.68778737668774   \\
};
\addplot+[smooth,color=red, loosely dashed, every mark/.append style={solid}, mark=none]
table[row sep=crcr]{
17	1.91962634990957 \\
18	1.73618934413823\\
19	1.68290275368604\\
20	0.99757491886916\\
21	0.935331113119607\\
22	0.922327529901508\\
23	0.758727665274486\\
24	0.762648948175936\\
25	0.702677952673315\\
26	0.688525557201545\\
27	0.653947031007457\\
28	0.65456272722857\\
29	0.578075692324222\\
30	0.563722040575154\\
31	0.545120505800747\\
32	0.557205969646181 \\
};

\addplot+[smooth,color=red, dotted, every mark/.append style={solid}, mark=|]
table[row sep=crcr]{
17	0.405427054023886 \\
18	0.361978285739012\\
19	0.358758649155936\\
20	0.334072885131309\\
21	0.329976129492221\\
22	0.316514999668179\\
23	0.318397092496309\\
24	0.314223514089517\\
25	0.310280445961039\\
26	0.311886594359130\\
27	0.30677189300402\\
28	0.307920086110529\\
29	0.305217185968679\\
30	0.305069693041412\\
31	0.304242379324131\\
32	0.302535476117532  \\
};

\addplot+[smooth,color=cyan, loosely dotted, every mark/.append style={solid}, mark=triangle]
table[row sep=crcr]{
17	0.140795003545487   \\
18	0.0864728512946572  \\
19	0.0597603883105515  \\
20	0.0443993694890806  \\
21	0.0346366723684060  \\
22	0.0279858168523434  \\
23	0.0232177802533512  \\
24	0.0196634058498785  \\
25	0.0169307273800340  \\
26	0.0147766093567004  \\
27	0.0130431160172158  \\
28	0.0116236868408425  \\
29	0.0104441024163167  \\
30	0.00945125329587347  \\
31	0.00860624190923420  \\
32	0.00787999765729379  \\
};

\addplot+[smooth,color=blue,densely dashed, every mark/.append style={solid}, mark=none]
table[row sep=crcr]{
17	0.657654609426848    \\
18	0.388554039577762  \\
19	0.283427947756222  \\
20	0.154341257996740  \\
21	0.141787055704898  \\
22	0.105589556669859  \\
23	0.0822475662090453  \\
24	0.0687329298781591  \\
25	0.0585842105203315  \\
26	0.0528929867586976  \\
27	0.0434323081059210  \\
28	0.0367526765919879  \\
29	0.0339761146725007  \\
30	0.0285856809998696  \\
31	0.0263711201023966  \\
32	0.0249435795787485  \\
};

\addplot+[smooth,color=blue,densely dotted, every mark/.append style={solid}, mark=none]
table[row sep=crcr]{
17	0.148474660131650   \\
18	0.0900319801471364  \\
19	0.0636807291228696  \\
20	0.0471292162375147  \\
21	0.0379528836674086  \\
22	0.0302868381296949  \\
23	0.0244222018142985  \\
24	0.0214379871037096  \\
25	0.0184386723694590  \\
26	0.0159850707565702  \\
27	0.0141227747555979  \\
28	0.0123763478322872  \\
29	0.0112219014583106  \\
30	0.0103814520818195  \\
31	0.00915661910587406  \\
32	0.00849875964919054 \\
};

\addplot+[smooth,color=blue,loosely dotted, every mark/.append style={solid}, mark=-]
table[row sep=crcr]{
17	0.142725007112984   \\
18	0.0870561612461422  \\
19	0.0617906995660773  \\
20	0.0459604897141749  \\
21	0.0356666360578312  \\
22	0.028878561373040  \\
23	0.0236842365396951  \\
24	0.020295882727854  \\
25	0.0175555844582813  \\
26	0.0154205062524997  \\
27	0.0136718933174149  \\
28	0.0120137610125321  \\
29	0.0108574396074060  \\
30	0.00965897290876398  \\
31	0.00891086938865661  \\
32	0.00811413255357704  \\
};

\end{axis}

\end{tikzpicture}%
\captionsetup{justification=centering,font=scriptsize}  
\caption{{\bf P}-based low-rank approximation error for \texttt{LowRankFastDecay}. Left: Spectral norm. Right: Frobenius norn.} 
\label{}      
\end{center}
\end{figure}

\begin{figure}[t]
\begin{center}       
%
%
%
\usetikzlibrary{positioning,calc}

\definecolor{mycolor1}{rgb}{0.00000,1.00000,1.00000}%
\definecolor{mycolor2}{rgb}{1.00000,0.00000,1.00000}%

\pgfplotsset{every axis label/.append style={font=\footnotesize},
every tick label/.append style={font=\footnotesize}
}

\begin{tikzpicture}[font=\footnotesize]

\begin{axis}[%
name=ber,
ymode=log,
width  = 0.35\columnwidth,
height = 0.3\columnwidth,
scale only axis,
xmin  = 11,
xmax  = 20,
xlabel= {$d$},
xmajorgrids,
ymin = 50,
ymax = 5.64e+03,
xtick       ={12,18},
xticklabels ={$12$,$18$},
ylabel={Magnitude},
ymajorgrids,
]

\addplot+[smooth,color=red,solid, every mark/.append style={solid}, mark=none]
table[row sep=crcr]{
11	5636.35687787020  \\
12	4925.86215170909 \\
13	2326.4980947720 \\
14	1901.70983647178 \\
15	2012.24786207971 \\
16	1659.59451971991 \\
17	1521.39860440510 \\
18	1265.26559533083 \\
19	1321.95862261687 \\
20	1033.68273204428  \\
};

\addplot+[smooth,color=red, loosely dashed, every mark/.append style={solid}, mark=none]
table[row sep=crcr]{
11	116.539256422154  \\
12	115.600960021492 \\
13	115.241444228000 \\
14	115.073334059353 \\
15	115.051182603433 \\
16	114.847050178712 \\
17	114.726961920963 \\
18	114.654921250964 \\
19	114.573933704008 \\
20	114.553936895188  \\
};

\addplot+[smooth,color=red, dotted, every mark/.append style={solid}, mark=|]
table[row sep=crcr]{
11	114.014758349723  \\
12	114.014064889646 \\
13	114.013846975214 \\
14	114.013804642516 \\
15	114.013650434060 \\
16	114.013588876112 \\
17	114.013520858650 \\
18	114.013563338970 \\
19	114.013481771221 \\
20	114.013410139427 \\
  };

\addplot+[smooth,color=blue,densely dashed, every mark/.append style={solid}, mark=none]
table[row sep=crcr]{
11	1016.82764845172    \\
12	997.305005740881  \\
13	495.495192912889  \\
14	374.074324170057  \\
15	370.315635025897  \\
16	328.896774769830  \\
17	298.268237947811  \\
18	224.223025565706  \\
19	262.008420625404  \\
20	179.382296138917  \\
};

\addplot+[smooth,color=blue,densely dotted, every mark/.append style={solid}, mark=none]
table[row sep=crcr]{
11	114.010451138475   \\
12	111.796033963906  \\
13	106.023388590607  \\
14	104.378534129815  \\
15	100.868709600386  \\
16	92.3170817465667  \\
17	85.7547429896962  \\
18	83.693568430520  \\
19	77.1809701712118  \\
20	73.1063584475662  \\
};

\addplot+[smooth,color=blue,loosely dotted, every mark/.append style={solid}, mark=-]
table[row sep=crcr]{
11	114.010616819829   \\
12	109.976335182664  \\
13	105.57913220923  \\
14	96.705666802201  \\
15	95.7106902280970  \\
16	82.1149450385976  \\
17	77.3888018010797  \\
18	72.549006078358  \\
19	67.8242990495210  \\
20	61.9132254982089  \\
};

\addplot+[smooth,color=cyan, loosely dotted, every mark/.append style={solid}, mark=triangle]
table[row sep=crcr]{
11	114.008777902283  \\
12	100.015010421027  \\
13	100.010004991460  \\
14	86.0174518981040  \\
15	86.0116321768124  \\
16	72.0208436921943  \\
17	72.0138925058439  \\
18	58.0258740096420  \\
19	58.0172442883791  \\
20	50.0413761366775  \\
};

\end{axis}

\begin{axis}[%
name=SumRate,
at={($(ber.east)+(35,0em)$)},
		anchor= west,
ymode=log,
width  = 0.34\columnwidth,
height = 0.3\columnwidth,
scale only axis,
xmin  = 11,
xmax  = 20,
xlabel= {$d$},
xmajorgrids,
ymin = 238 ,
ymax = 1.82e+04,
xtick       ={12,18},
xticklabels ={$12$,$18$},
ylabel={},
ymajorgrids,
]

\addplot+[smooth,color=red,solid, every mark/.append style={solid}, mark=none]
table[row sep=crcr]{
11	18144.5001956024  \\
12	15857.2831195474 \\
13	7489.43795617341 \\
14	6121.96410687179 \\
15	6477.80694484738 \\
16	5342.54905077218 \\
17	4897.67022801591 \\
18	4073.12956567864 \\
19	4255.63515696220 \\
20	3327.62047190587  \\
};
\addplot+[smooth,color=red, loosely dashed, every mark/.append style={solid}, mark=none]
table[row sep=crcr]{
11	375.161936471659  \\
12	372.141382664608 \\
13	370.984033240736 \\
14	370.442854771381 \\
15	370.371545039764 \\
16	369.714404106259 \\
17	369.327817262412 \\
18	369.095904702719 \\
19	368.835190451754 \\
20	368.770816937194 \\
};

\addplot+[smooth,color=red, dotted, every mark/.append style={solid}, mark=|]
table[row sep=crcr]{
11	367.035099090430  \\
12	367.032866711114 \\
13	367.032165203289 \\
14	367.032028926319 \\
15	367.031532500094 \\
16	367.031334333403 \\
17	367.031115372311 \\
18	367.031252124507 \\
19	367.030989542476 \\
20	367.030758946164 \\
};

\addplot+[smooth,color=cyan, loosely dotted, every mark/.append style={solid}, mark=triangle]
table[row sep=crcr]{
11	348.872479692082   \\
12	329.718069944203  \\
13	314.183072965766  \\
14	297.840565134805  \\
15	285.149084180149  \\
16	271.867613626666  \\
17	262.154529644107  \\
18	252.069428330925  \\
19	245.299805634853  \\
20	238.340080577479  \\
};

\addplot+[smooth,color=blue,densely dashed, every mark/.append style={solid}, mark=none]
table[row sep=crcr]{
11	1268.97533084640   \\
12	1214.99693709223  \\
13	722.848291893808  \\
14	613.176856514876  \\
15	604.717522609172  \\
16	546.86470986630  \\
17	512.30090205267  \\
18	454.191469807033  \\
19	465.534633374322  \\
20	402.848256082848  \\
};

\addplot+[smooth,color=blue,densely dotted, every mark/.append style={solid}, mark=none]
table[row sep=crcr]{
11	353.582301692196   \\
12	337.864180124980  \\
13	323.875903662293  \\
14	314.673183958138  \\
15	302.272554577604  \\
16	289.639587527687  \\
17	278.920244044250  \\
18	271.553999628027  \\
19	261.505442130138  \\
20	254.475924281927  \\
};

\addplot+[smooth,color=blue,loosely dotted, every mark/.append style={solid}, mark=-]
table[row sep=crcr]{
11	350.426205959107   \\
12	334.632438091783  \\
13	318.998599869687  \\
14	303.969787120491  \\
15	291.698607881013  \\
16	278.322731752821  \\
17	268.770997578563  \\
18	259.558127589681  \\
19	252.163479515586  \\
20	244.863797736909  \\
};

\end{axis}

\end{tikzpicture}%
\captionsetup{justification=centering,font=scriptsize}  
\caption{{\bf P}-based low-rank approximation error for \texttt{impcol\_e}. Left: Spectral norm. Right: Frobenius norn.} 
\label{figPLRimpcole}       
\end{center}
\end{figure}

\appendices


\section{Proof of Theorem 1}
\label{Proo_R22bound}

\begin{remark}\label{ReProof}
Before proving our results, we should add that we use the SVD partitioning of the matrix A  as in \eqref{eqSVD}, a $2\times 2$ block $\bf \Sigma$ and $2\times 1$ blocks $\bf U$ and $\bf V$. It is however possible to use a $3\times 3$ block for $\bf \Sigma$ and $3\times 1$ blocks for the orthogonal matrices to derive more bounds on the principal angles between subspaces based on the CS decomposition \cite{PAIGE1994}. This is particularly useful when the SVD is used to construct the low-rank factorization; see \cite{STEWART2001, Nakatsukasa17}.
\end{remark}

Exploiting the Cauchy interlacing theorem \cite[Theorem 4.3.17]{HornJohnson852012}, we have
\begin{equation}\notag
\begin{aligned}
\lambda_i({\bf A}^T{\bf A})  & \ge \lambda_i({\bf R}^T{\bf R}) \ge \lambda_i({\bf R}_{11}^T{\bf R}_{11})\\ &  = \lambda_i(\bar{\bf P}_1^T{\bf A}^T{\bf Q}_1{\bf Q}_1^T{\bf A}\bar{\bf P}_1).
\end{aligned}
\end{equation}
The last relation is due to the first equality in \eqref{eqAPbarPart}. In order to compute ${\bf Q}_1{\bf Q}_1^T$, we have for $\mathcal{R}({\bf Q})$: 
\begin{equation}\notag
\mathcal{R}({\bf Q}) = \mathcal{R}({\bf A}\bar{\bf P})= \mathcal{R}({\bf A}({\bf A}^T{\bf A})^q{\bf A}^T {\bf \Phi}).
\end{equation}

We define a non-singular ${\bf Y} \triangleq [\widehat{\bf \Phi}_1^\dagger{\bf \Sigma}_k^{-(2q+2)} \quad \bar{\bf Y}] \in \mathbb R^{d \times d}$ such that $\widehat{\bf \Phi}_1\bar{\bf Y}={\bf 0}$. Assuming the rank of $\widehat{\bf \Phi}_1$ is $k$, we get $\widehat{\bf \Phi}_1\widehat{\bf \Phi}_1^\dagger = {\bf I}$. Then, the following product is formed and its QR factorization is computed:
\begin{equation} \notag
\begin{aligned}
{\bf A}({\bf A}^T{\bf A})^q{\bf A}^T{\bf \Phi}{\bf Y} & = {\bf U}\begin{bmatrix}
{\bf I} & {\bf 0} \\
{\bf S} & {\bf \Sigma}_\perp^{2q+2}\widehat{\bf \Phi}_2\bar{\bf Y}
\end{bmatrix}=\dot{\bf Q}\dot{\bf R} \\& = [\dot{\bf Q}_1\quad \dot{\bf Q}_2]
\begin{bmatrix}
\dot{\bf R}_{11} & \dot{\bf R}_{12}  \\
{\bf 0} & \dot{\bf R}_{22}
\end{bmatrix},
\end{aligned}
\end{equation}
where ${\bf S} \triangleq {\bf \Sigma}_\perp^{2q+2}\widehat{\bf \Phi}_2\widehat{\bf \Phi}_{1}^\dagger {\bf \Sigma}_k^{-(2q+2)}$. Non-singularity of $\bf Y$ implies $\mathcal{R}({\bf Q}) =\mathcal{R}(\dot{\bf Q})$. This means that there is an orthogonal matrix $\bf W$ such that ${\bf Q}= \dot{\bf Q}{\bf W}$, thus giving ${\bf Q}{\bf Q}^T = \dot{\bf Q}\dot{\bf Q}^T$. From the above equation, we obtain $ \dot{\bf Q}_1 = {\bf U}[{\bf I}^T  \quad {\bf S}^T]^T\dot{\bf R}_{11}^{-1}$, and hence we get
\begin{equation}\label{eqQ1erQ1dot}
{\bf Q}_1{\bf Q}_1^T =  \dot{\bf Q}_1\dot{\bf Q}_1^T = {\bf U}
\begin{bmatrix}
\bar{\bf S} & \bar{\bf S}{\bf S}^T  \\
{\bf S}\bar{\bf S} & {\bf S}\bar{\bf S}{\bf S}^T 
\end{bmatrix}{\bf U}^T,
\end{equation}
where $\bar{\bf S}$ is defined as:
\begin{equation}
\begin{aligned}
\bar{\bf S} \triangleq \dot{\bf R}_{11}^{-1} \dot{\bf R}_{11}^{-T} = 
(\dot{\bf R}_{11}^T \dot{\bf R}_{11})^{-1} = ({\bf I}+{\bf S}^T{\bf S}) ^{-1}.
\end{aligned}
\end{equation}

Now, the product $\bar{\bf P}_1^T{\bf A}^T{\bf Q}_1{\bf Q}_1^T{\bf A}\bar{\bf P}_1$ is formed:
\begin{equation}\notag
\begin{aligned}
& \bar{\bf P}_1^T{\bf A}^T{\bf Q}_1{\bf Q}_1^T{\bf A}\bar{\bf P}_1 = \\
& \bar{\bf P}_1^T{\bf V}
\begin{bmatrix}
{\bf \Sigma}_k\bar{\bf S}{\bf \Sigma}_k & {\bf \Sigma}_k\bar{\bf S}{\bf S}^T{\bf \Sigma}_\perp  \\
{\bf \Sigma}_\perp{\bf S}\bar{\bf S}{\bf \Sigma}_k & {\bf \Sigma}_\perp{\bf S}\bar{\bf S}{\bf S}^T{\bf \Sigma}_\perp 
\end{bmatrix}{\bf V}^T\bar{\bf P}_1.
\end{aligned}
\end{equation}

As ${\bf \Sigma}_k\bar{\bf S}{\bf \Sigma}_k$ is a submatrix,  we have for $i=1,..., k$
\begin{equation}
\lambda_i(\bar{\bf P}_1^T{\bf A}^T{\bf Q}_1{\bf Q}_1^T{\bf A}\bar{\bf P}_1) \ge \lambda_i({\bf \Sigma}_k\bar{\bf S}{\bf \Sigma}_k).
\end{equation}

Due to matrix partial orderings \cite[Section 7.7]{HornJohnson852012}, it follows 
\begin{equation}\label{equpbStS}
\begin{aligned}
{\bf S}^T{\bf S} \preceq \sigma_{k+1}^{4q+4}\|\widehat{\bf \Phi}_{2}\widehat{\bf \Phi}_{1}^\dagger\|_2^2 {\bf \Sigma}_k^{-(4q+4)}  = {\bf \Delta}^{(4q+4)}\|\widehat{\bf \Phi}_{2}\widehat{\bf \Phi}_{1}^\dagger\|_2^2, 
\end{aligned}
\end{equation}
where ${\bf \Delta} = \text{diag}(\delta_1, ..., \delta_k) \in \mathbb R^{k \times k}$ has entries $\delta_i= \frac{\sigma_{k+1}}{\sigma_i}$. Accordingly, we have
\begin{equation}\notag
 {\bf \Sigma}_k({\bf I} + {\bf S}^T{\bf S})^{-1}{\bf \Sigma}_k \succeq {\bf \Sigma}_k({\bf I} + {\bf \Delta}^{(4q+4)}\|\widehat{\bf \Phi}_{2}\widehat{\bf \Phi}_{1}^\dagger\|_2^2)^{-1}{\bf \Sigma}_k,
\end{equation}
which results in 
\begin{equation}\notag
\begin{aligned}
\lambda_i({\bf A}^T{\bf A}) & \ge \lambda_i({\bf R}_{11}^T{\bf R}_{11})  \ge \lambda_i({\bf \Sigma}_k({\bf I} + {\bf S}^T{\bf S})^{-1}{\bf \Sigma}_k) \\ & \ge \frac{\sigma_i^2}{1 + \delta_i^{4q+4}\|\widehat{\bf \Phi}_{2}\widehat{\bf \Phi}_{1}^\dagger\|_2^2}.
\end{aligned}
\end{equation}

The bounds for the first $k$ singular values of ${\bf R}_{11}$ are obtained through taking the square root of the last identity. 

To prove \eqref{eqThR22}, we get from the second equality in  \eqref{eqAPbarPartSep}
\begin{equation}\label{eqR22_decomposed}
\begin{aligned}
\|{\bf R}_{22}\|_{2,F} & = \|{\bf Q}_2^T{\bf A}\bar{\bf P}_2\|_{2,F} 
= \|{\bf Q}_2{\bf Q}_2^T{\bf A}\bar{\bf P}_2\|_{2,F} \\
& \le \|({\bf I} - {\bf Q}_1{\bf Q}_1^T){\bf A}\|_{2,F} \\
& \le  \|({\bf I} - {\bf Q}_1{\bf Q}_1^T){\bf A}_k\|_{2,F} +  \|{\bf A}_\perp\|_{2,F}.
\end{aligned}
\end{equation}

The second equality results from the unitary invariance of the 2- and Frobenius norms, the first inequality follows because $\bar{\bf P}_2$ has orthonormal columns,  and the second inequality is due to the triangle inequality after writing ${\bf A}={\bf A}_k+{\bf A}_\perp$.  To bound $\|({\bf I} - {\bf Q}_1{\bf Q}_1^T){\bf A}_k\|_{2,F}$ from above, we write ${\bf A}_k = {\bf U}[{\bf \Sigma}_k \quad {\bf 0}; {\bf 0} \quad {\bf 0}]{\bf V}^T$, which, together with \eqref{eqQ1erQ1dot}, gives
\begin{equation}\notag
\begin{aligned}
({\bf I} - {\bf Q}_1{\bf Q}_1^T){\bf A}_k = {\bf U}
	\begin{bmatrix}
	({\bf I} - \bar{\bf S}){\bf \Sigma}_k  \\
	- {\bf S} \bar{\bf S}{\bf \Sigma}_k
	\end{bmatrix} {\bf V}^T.
\end{aligned}
\end{equation}

It follows that
\begin{equation}
\begin{aligned}\notag
& \|({\bf I} - {\bf Q}_1{\bf Q}_1^T){\bf A}_k \|_{2,F}^2 = \|{\bf \Sigma}_k ({\bf I} - \bar{\bf S}){\bf \Sigma}_k\|_{2,F}^2 
\\&= \|{\bf \Sigma}_k {\bf S}^T({\bf I}+{\bf S}{\bf S}^T)^{-1}{\bf S}{\bf \Sigma}_k\|_{2,F}^2=
\|{\bf \Sigma}_k ^{-(2q+1)}(\widehat{\bf \Phi}_2\widehat{\bf \Phi}_{1}^\dagger)^T
\\&\times {\bf \Sigma}_\perp^{2q+2}({\bf I}+{\bf S}{\bf S}^T)^{-1}{\bf \Sigma}_\perp^{2q+2}\widehat{\bf \Phi}_2\widehat{\bf \Phi}_{1}^\dagger {\bf \Sigma}_k^{-(2q+1)}\|_{2,F}^2 \\
& \le \delta_k^{4q+2}\|\widehat{\bf \Phi}_{2}\widehat{\bf \Phi}_{1}^\dagger \|_2^2 \|({\bf I}+{\bf S}{\bf S}^T)^{-1}\|_2^2\|{\bf \Sigma}_\perp\|_{2,F}^2.
\end{aligned} 
\end{equation} 

The second equality results from the Sherman–Morrison –Woodbury formula \cite[Section 0.7.4]{HornJohnson852012}, and the last relation follows from the strong submultiplicativity property for the Frobenius norm \cite[equation 9.3.13]{Bernstein09} holding for any matrices $\bf A$ and $\bf B$ with appropriate dimensions:   
\begin{equation}\label{eqStrongSub}
\begin{aligned}
 & \|{\bf A}{\bf B}\|_F \le \|{\bf A}\|_F\sigma_1({\bf B}),  \\
 & \|{\bf A}{\bf B}\|_F \le \sigma_1({\bf A})\|{\bf B}\|_F.    
\end{aligned}
\end{equation}
The non-zero eigenvalues of $({\bf I}+{\bf S}{\bf S}^T)^{-1}$ are:
\begin{equation}
{1}/{1+\sigma_i^2({\bf S})}.
\end{equation}
$\sigma_1({\bf S})$ is bounded from below by:  
\begin{equation}
\sigma_1({\bf S}) \ge \gamma^{2q+2}\|\widehat{\bf \Phi}_{2}\widehat{\bf \Phi}_{1}^\dagger \|_2,
\end{equation}
 where $\gamma=\sigma_n/\sigma_1$, and we have used the following relations \cite[p. 80]{StewartSun90} that hold for any unitary invariant norm: 
\begin{equation}\notag
\begin{aligned}
 & \|{\bf A}{\bf B}\| \ge \|{\bf A}\| \sigma_\text{min}({\bf B}),  \\
 & \|{\bf A}{\bf B}\| \ge \sigma_\text{min}({\bf A})\|{\bf B}\|.
\end{aligned}
\end{equation}
Thus
\begin{equation}\label{eq2nIminSSt}
\|({\bf I}+{\bf S}{\bf S}^T)^{-1}\|_2  \le \dfrac{1}{1 + \gamma^{4q+4}\|\widehat{\bf \Phi}_{2}\widehat{\bf \Phi}_{1}^\dagger \|_2^2}.
\end{equation}

Upon substitution and taking the square root, we obtain 
\begin{equation}
\begin{aligned}\label{eqIminQ1QitUB}
\|({\bf I} - {\bf Q}_1{\bf Q}_1^T){\bf A}_k \|_{2,F}  \le \dfrac{\delta_k^{2q+1}\|\widehat{\bf \Phi}_{2}\widehat{\bf \Phi}_{1}^\dagger\|_2}{1 + \gamma^{4q+4}\|\widehat{\bf \Phi}_{2}\widehat{\bf \Phi}_{1}^\dagger \|_2^2}\|{\bf \Sigma}_\perp\|_{2,F},
\end{aligned} 
\end{equation}

Inserting this bound into \eqref{eqR22_decomposed} gives the desired result.   \QEDB

\section{Proof of Theorem 2}
\label{Proo_CanABk}

We use a result from \cite[equation 13]{BjorckGolub73}, which states in our notation that the SVD of the matrix product $({\bf I}-{\bf P}_Q){\bf U}_k$ gives the sines of canonical angles  between the two subspaces spanned by $\bf Q$ and ${\bf U}_k$. In particular, let
\begin{equation}\notag
({\bf I}-{\bf P}_Q){\bf U}_k = {\bf Y}_Q{\bf S}_Q{\bf Z}_Q^T.
\end{equation}
Then, ${\bf S}_Q \in \mathbb R^{k \times k}$ contains on the diagonal the sines of canonical angles $\text{sin}\theta_i$. Further 
\begin{equation}\notag
\begin{aligned}
{\bf U}_k^T({\bf I}-{\bf P}_{{Q}}){\bf U}_k = {\bf Z}_Q{\bf S}_Q^2{\bf Z}_Q^T & = {\bf U}_k^T({\bf I}-{\bf P}_{\dot {Q}}){\bf U}_k \\ & \preceq {\bf U}_k^T({\bf I}-{\bf P}_{\dot {Q}_1}){\bf U}_k.
\end{aligned}
\end{equation}

The last relation holds because $\dot{\bf Q}_1$ contains $k$ columns of $\dot{\bf Q}$. Upon substitution, we get
\begin{equation}\notag
\begin{aligned}
{\bf U}_k^T({\bf I}-{\bf P}_{\dot {Q}_1}){\bf U}_k=
 [{\bf I} \quad {\bf 0}]
\begin{bmatrix}
{\bf I} - \bar{\bf S} & -\bar{\bf S}{\bf S}^T  \\
-{\bf S}\bar{\bf S} & {\bf I} -{\bf S}\bar{\bf S}{\bf S}^T 
\end{bmatrix}
\begin{bmatrix}
{\bf I} \\
{\bf 0}
\end{bmatrix} = {\bf I} - \bar{\bf S}.
\end{aligned}
\end{equation}

The matrix ${\bf I}-\bar{\bf S} = {\bf S}^T{\bf S}({\bf I}+{\bf S}^T{\bf S})^{-1}$ is positive semidefinite, whose eigenvalues are given by \cite[p. 148]{StewartSun90}:   
\begin{equation} \notag
\begin{aligned}
\lambda_i({\bf I}-\bar{\bf S}) = \frac{{\sigma}_i^2({\bf S})}{1+{\sigma}_i^2({\bf S})}, \quad i=1,...k.
\end{aligned}
\end{equation}

We obtain for the singular values of ${\bf S}$ 
\begin{equation}\notag
\begin{aligned}
{\sigma}_i({\bf S}) & \le {\sigma}_1({\bf \Sigma}_\perp^{2q+2}\widehat{\bf \Phi}_2\widehat{\bf \Phi}_1^\dagger){\sigma}_i({\bf \Sigma}_k^{-(2q+2)})  \\
& \le \frac{{\sigma}_{k+1}^{2q+2}}{{\sigma}_{k-i+1}^{2q+2}}\|\widehat{\bf \Phi}_2\widehat{\bf \Phi}_1^\dagger\|_2,
\end{aligned}
\end{equation}
where we have made use of \cite[equation 9.6.2]{Bernstein09}:
\begin{equation}\label{eqSigm1xSigmai}
\sigma_i({\bf A}{\bf B}) \le \sigma_1({\bf A})\sigma_i({\bf B}).
\end{equation}

Since ${\bf Z}_Q{\bf S}_Q^2{\bf Z}_Q^T \preceq {\bf I}-\bar{\bf S}$,  the Weyl's inequality \cite[Theorem 8.4.9]{Bernstein09} implies that
\begin{equation} 
\text{sin}^2\theta_i  \le \lambda_i({\bf I}-\bar{\bf S}).
\label{eq_sineTtilde}
\end{equation}

Upon substitution, taking the square root, and renaming $i \leftarrow k-i+1$, the result in \eqref{eqThCanABkSinT} follows.  

The proof of \eqref{eqThCanABkSinP} is similar to that of  \eqref{eqThCanABkSinT}. We first have
\begin{equation}\label{eqCanAPVk}
({\bf I}-{\bf P}_P){\bf V}_k = {\bf Y}_P{\bf S}_P{\bf Z}_P^T,
\end{equation}
where the matrix ${\bf S}_P \in \mathbb R^{k \times k}$ contains on the diagonal the sines of canonical angles $\text{sin}\phi_i$ between the two subspaces spanned by ${\bf P}$ and ${\bf V}_k$. To compute ${\bf P}_P$, as the columns of $\widetilde{\bf P}$ are orthonormal, we have for $\mathcal{R}({\bf P})$: 
\begin{equation}\notag
\mathcal{R}({\bf P}) =\mathcal{R}(\bar{\bf P}) = \mathcal{R}(({\bf A}^T{\bf A})^q{\bf A}^T{\bf \Phi}).
\end{equation}

We define a non-singular matrix ${\bf X} \triangleq [\widehat{\bf \Phi}_1^\dagger{\bf \Sigma}_k^{-(2q+1)} \quad \bar{\bf X}] \in \mathbb R^{d \times d}$ such that $\widehat{\bf \Phi}_1\bar{\bf X}={\bf 0}$. Then the following matrix product is formed and its QR factorization is computed:
\begin{equation} \notag
\begin{aligned}
({\bf A}^T{\bf A})^q{\bf A}^T{\bf \Phi}{\bf X} & = {\bf V}\begin{bmatrix}
{\bf I} & {\bf 0} \\
{\bf J} & {\bf \Sigma}_\perp^{2q+1}\widehat{\bf \Phi}_2\bar{\bf X}
\end{bmatrix}=\bar{\bf Q}\bar{\bf R} \\
& = [\bar{\bf Q}_1\quad \bar{\bf Q}_2]
\begin{bmatrix}
\bar{\bf R}_{11} & \bar{\bf R}_{12}  \\
{\bf 0} & \bar{\bf R}_{22}
\end{bmatrix},
\end{aligned}
\end{equation}
where ${\bf J} \triangleq {\bf \Sigma}_\perp^{2q+1}\widehat{\bf \Phi}_2\widehat{\bf \Phi}_{1}^\dagger {\bf \Sigma}_k^{-(2q+1)}$. Non-singularity $\bf X$ implies $\mathcal{R}(\bar{\bf P}) =\mathcal{R}(\bar{\bf Q})$. This means that for an orthogonal $\bf W$, we have $\bar{\bf P}= \bar{\bf Q}{\bf W}$, which gives $\bar{\bf P}\bar{\bf P}^T = \bar{\bf Q}\bar{\bf Q}^T$. From the above equation, we get $\bar{\bf Q}_1 = {\bf V}[{\bf I}^T  \quad {\bf J}^T]^T\bar{\bf R}_{11}^{-1}$, and hence obtain for $\bar{\bf P}_1\bar{\bf P}_1^T$
\begin{equation} \notag
\begin{aligned}
\bar{\bf P}_1\bar{\bf P}_1^T= \bar{\bf Q}_1\bar{\bf Q}_1^T = {\bf V}\begin{bmatrix}
\bar{\bf J} & \bar{\bf J}{\bf J}^T \\
{\bf J}\bar{\bf J} & {\bf J}\bar{\bf J}{\bf J}^T
\end{bmatrix}{\bf V}^T,
\end{aligned}
\end{equation}

where $\bar{\bf J}$ is defined as:
\begin{equation}
\begin{aligned}
\bar{\bf J} \triangleq \bar{\bf R}_{11}^{-1} \bar{\bf R}_{11}^{-T} = 
(\bar{\bf R}_{11}^T \bar{\bf R}_{11})^{-1} = ({\bf I}+{\bf J}^T{\bf J}) ^{-1}.
\end{aligned}
\end{equation}

Having obtained ${\bf P}_P $, it follows from \eqref{eqCanAPVk} that 
\begin{equation}\notag
\begin{aligned}
{\bf V}_k^T({\bf I}-{\bf P}_{{P}}){\bf V}_k & = {\bf Z}_P{\bf S}_P^2{\bf Z}_P^T = {\bf V}_k^T({\bf I}-{\bf P}_{\bar {Q}}){\bf V}_k \\
& \preceq {\bf V}_k^T({\bf I}-{\bf P}_{\bar {Q}_1}){\bf V}_k \\
& = [{\bf I} \quad {\bf 0}]
\begin{bmatrix}
{\bf I} - \bar{\bf J} & -\bar{\bf J}{\bf J}^T  \\
-{\bf J}\bar{\bf J} & {\bf I} -{\bf J}\bar{\bf J}{\bf J}^T 
\end{bmatrix}
\begin{bmatrix}
{\bf I} \\
{\bf 0}
\end{bmatrix} = {\bf I} - \bar{\bf J}.
\end{aligned}
\end{equation}

The rest of the proof is analogous to that of \eqref{eqThCanABkSinT}, we thus omit it. \QEDB

\section{Proof of Theorem 3}
\label{ProoThCanASMax}

From the proof of Theorem \ref{ThCanABk}, we have
\begin{equation}
\begin{aligned}\notag
{\bf Z}_Q{\bf S}_Q^2{\bf Z}_Q^T \preceq {\bf I}-\bar{\bf S}= {\bf S}^T{\bf S}({\bf I}+{\bf S}^T{\bf S})^{-1}.
\end{aligned} 
\end{equation}
It follows that
\begin{equation}
\begin{aligned}\notag
 \|{\bf Z}_U{\bf S}_U^2{\bf Z}_U^T\|_{2,F} \le \|{\bf S}^T{\bf S}\|_{2,F}\|({\bf I}+{\bf S}^T{\bf S})^{-1}\|_2, 
\end{aligned} 
\end{equation}
where we have made use of \eqref{eqStrongSub}. Using \eqref{eq2nIminSSt} and \eqref{eqStrongSub}, the result follows. The proof for the other bound is analogous, we therefore omit it. \QEDB

\section{Proof of Theorem 4}
\label{Proo_ThCanASIn}

The proof of this theorem follows that of Theorem \ref{ThCanABk}. For $i=1,..., k$, let ${\bf u}_i$ be the $i$th left singular vector of $\bf A$. We then have
\begin{equation}\notag
\begin{aligned}   
&\text{sin}^2\angle(\mathcal{R}({\bf Q}), \mathcal{R}({\bf u}_i)) = 
{\bf u}_i^T({\bf I}-{\bf P}_Q){\bf u}_i \preceq {\bf u}_i^T({\bf I}-{\bf P}_{\dot {Q}_1}){\bf u}_i \\
& = [{\bf e}_i^T \quad {\bf 0}]
\begin{bmatrix}
{\bf I} - \bar{\bf S} & -\bar{\bf S}{\bf S}^T  \\
-{\bf S}\bar{\bf S} & {\bf I} -{\bf S}\bar{\bf S}{\bf S}^T 
\end{bmatrix}
\begin{bmatrix}
{\bf e}_i \\
{\bf 0}
\end{bmatrix} \preceq {\bf e}_i^T({\bf I} - \bar{\bf S}){\bf e}_i .
\end{aligned}
\end{equation}
Writing ${\bf I} - \bar{\bf S} = {\bf S}^T({\bf I}+{\bf S}{\bf S}^T)^{-1}{\bf S}$, it follows
\begin{equation}
\begin{aligned}\notag
\text{sin}^2\angle(\mathcal{R}({\bf Q}), \mathcal{R}({\bf u}_i)) & \le \delta_i^{4q+4}\|\widehat{\bf \Phi}_{2}\widehat{\bf \Phi}_{1}^\dagger \|_2^2 \|({\bf I}+{\bf S}{\bf S}^T)^{-1}\|_2^2 \\
& \le \dfrac{\delta_i^{4q+4}\|\widehat{\bf \Phi}_{2}\widehat{\bf \Phi}_{1}^\dagger \|_2^2}{1 + \gamma^{4q+4}\|\widehat{\bf \Phi}_{2}\widehat{\bf \Phi}_{1}^\dagger \|_2^2}.
\end{aligned} 
\end{equation}

The result follows after taking the square root. The proof for the other bound follows likewise, we therefore omit it. \QEDB

\section{Proof of Theorem 5}
\label{Pro_ThLRAEB}

To prove the first bound, we have
\begin{equation}
\begin{aligned}
 \|{\bf A} - & {\bf Q}{\bf Q}^T{\bf A}\|_{2,F}   \le \|{\bf A} - {\bf Q}[{\bf Q}^T{\bf A}]_k\|_{2,F} 
\\& \le 
\|{\bf A} - {\bf Q}_1{\bf Q}_1^T{\bf A}\|_{2,F}
\\& \le  \|{\bf A}_\perp -{\bf Q}_1{\bf Q}_1^T{\bf A}_\perp\|_{2,F} +  \|{\bf A}_k - {\bf Q}_1{\bf Q}_1^T{\bf A}_k\|_{2,F} \\
& \le \|{\bf A}_\perp\|_{2,F} + \|{\bf A}_k - {\bf Q}_1{\bf Q}_1^T{\bf A}_k\|_{2,F}.
\end{aligned}
\end{equation}

The first line is due the truncation, see \cite[Lemma 2.2]{Gu2015}, the second line is due to the optimality of the SVD, the third line is due to the triangle inequality, and the first term in the fourth line is due to  \eqref{eqStrongSub}. Inserting \eqref{eqIminQ1QitUB} into the last relation, the result follows. The second bound follows similarly. \QEDB

\section{Proof of Theorem 6}
\label{Proo_ThSigmR22Probab}

Owing to the statistical independence of $\widehat{\bf \Phi}_1$ and $\widehat{\bf \Phi}_2$, we take expectations of the identity \eqref{eqThSigR11Prob} in turn:
\begin{equation}\notag
\begin{aligned}
  & \mathbb{E}(\sigma_i({\bf R}_{11})) \ge \mathbb{E}_{\widehat{\bf \Phi}_1}  \left(\mathbb{E}_{\widehat{\bf \Phi}_2}\left[\frac{\sigma_i}{\sqrt{1 + \delta_i^{4q+4}{\|{\widehat{\bf \Phi}_2}}{\widehat{\bf \Phi}_1}^\dagger\|_2^2}}\right]\right) \\
& \ge \mathbb{E}_{\widehat{\bf \Phi}_1} \left(\frac{\sigma_i}{\sqrt{1 +\delta_i^{4q+4} \eta_1^2\|{\widehat{\bf \Phi}_1}^\dagger\|_2^2 }}\right)\ge \frac{\sigma_i}{\sqrt{1 + \delta_i^{4q+4}\eta^2}}.
\end{aligned}
\label{eqThr3_proof}
\end{equation}             
The second and third inequalities follow from Proposition 5.4 and Proposition 5.5 of \cite{Gu2015}, respectively.  

To prove \eqref{eqThSigR22Prob}, the identity is first simplified to:
\begin{equation}\notag
\begin{aligned}
\|{\bf R}_{22}\|_{2,F} \le \|{\bf \Sigma}_\perp\|_{2,F} + {\delta_k^{2q+1}\|\widehat{\bf \Phi}_{2}\widehat{\bf \Phi}_{1}^\dagger\|_2}\|{\bf \Sigma}_\perp\|_{2,F}.
\end{aligned}
\end{equation}
 
We then take expectations in turn:
\begin{equation}\notag
\begin{aligned}
& \mathbb{E} \|{\bf R}_{22}\|_{2,F} 
\\& \le  \mathbb{E}_{\widehat{\bf \Phi}_1} \big(\mathbb{E}_{\widehat{\bf \Phi}_2} \big[\|{\bf \Sigma}_\perp\|_{2,F} + {\delta_k^{2q+1}\|\widehat{\bf \Phi}_{2}\widehat{\bf \Phi}_{1}^\dagger\|_2}\|{\bf \Sigma}_\perp\|_{2,F} \big]\big)\\
& \le  \mathbb{E}_{\widehat{\bf \Phi}_1} \big(\|{\bf \Sigma}_\perp\|_{2,F} + {\delta_k^{2q+1}\mathbb{E}_{\widehat{\bf \Phi}_2}\|\widehat{\bf \Phi}_{2}\widehat{\bf \Phi}_{1}^\dagger\|_2}\|{\bf \Sigma}_\perp\|_{2,F} \big)\\
& 
 \le \mathbb{E}_{\widehat{\bf \Phi}_1} \big(\|{\bf \Sigma}_\perp\|_{2,F} + \delta_k^{2q+1}\big[\|\widehat{\bf \Phi}_{1}^\dagger\|_F +
 \sqrt{\bar{m}}\|\widehat{\bf \Phi}_{1}^\dagger\|_2\|\big]{\bf \Sigma}_\perp\|_{2,F} \big)\\
 & \le \|{\bf \Sigma}_\perp\|_{2,F} + \delta_k^{2q+1}\Bigg[\sqrt{\dfrac{k}{p-1}}+
 {\dfrac{e\sqrt{\bar{m}(p+k)}}{p}}\Bigg]\|{\bf \Sigma}_\perp\|_{2,F},
\end{aligned}
\end{equation}
where $\bar{m} = m-k$. In the above relations, the fourth line results from \cite[Proposition 10.1]{HMT2011}, and the last line from \cite[Proposition 10.2]{HMT2011} together with the H$\ddot{\text{o}}$lder's inequality: 
\begin{equation}\notag
\mathbb{E}_{\widehat{\bf \Phi}_1}\|\widehat{\bf \Phi}_1^\dagger\|_F \le \Big(\mathbb{E}_{\widehat{\bf \Phi}_1}\|\widehat{\bf \Phi}_1^\dagger\|_F^2\Big)^{1/2}.
\end{equation}
 \QEDB

\section{Proof of Theorem 7}
\label{Proo_ThSigmCanAngiProbab}

To prove the first bound, through leveraging the statistical independence of $\widehat{\bf \Phi}_1$ and $\widehat{\bf \Phi}_2$, we take expectations of the first identity over $\widehat{\bf \Phi}_2$ and $\widehat{\bf \Phi}_1$:
\begin{equation}\notag
\begin{aligned}
&  \mathbb{E}\text{sin}\theta_i  \le \mathbb{E}_{\widehat{\bf \Phi}_1} \Bigg( \mathbb{E}_{\widehat{\bf \Phi}_2} \left[
\frac{\delta_i^{2q+2}{\|{\widehat{\bf \Phi}_2}\widehat{\bf \Phi}_1}^\dagger\|_2} {\sqrt{
{1 + \delta_i^{4q+4}{\|{\widehat{\bf \Phi}_2}\widehat{\bf \Phi}_1}^\dagger\|_2^2}}}, \right]\Bigg)\\ 
& \le \mathbb{E}_{\widehat{\bf \Phi}_1} \Bigg(
\sqrt{\frac{\delta_i^{4q+4}{\omega_1^2}\|\widehat{\bf \Phi}_1^\dagger\|_2^2}{1 + \delta_i^{4q+4}{{\omega_1^2}}\|{\widehat{\bf \Phi}_1}^\dagger\|_2^2}} \Bigg) 
 \le 
\sqrt{\frac{\delta_i^{4q+4}{\omega_1^2}{\omega_2^2}}{1 + \delta_i^{4q+4}{{\omega_1^2}}{\omega_2^2}}}.
\end{aligned}
\end{equation}
The second and third inequalities follow from Proposition 3 and Proposition 5  of  \cite{MFKDeTSP18}, respectively. The second bound follows analogously, we therefore omit it. \QEDB

\bibliographystyle{IEEEtran}
\bibliography{Mybib}




\begin{figure}[!htb]
	\centering
	\includegraphics[width=0.46\textwidth,height=5.7cm]{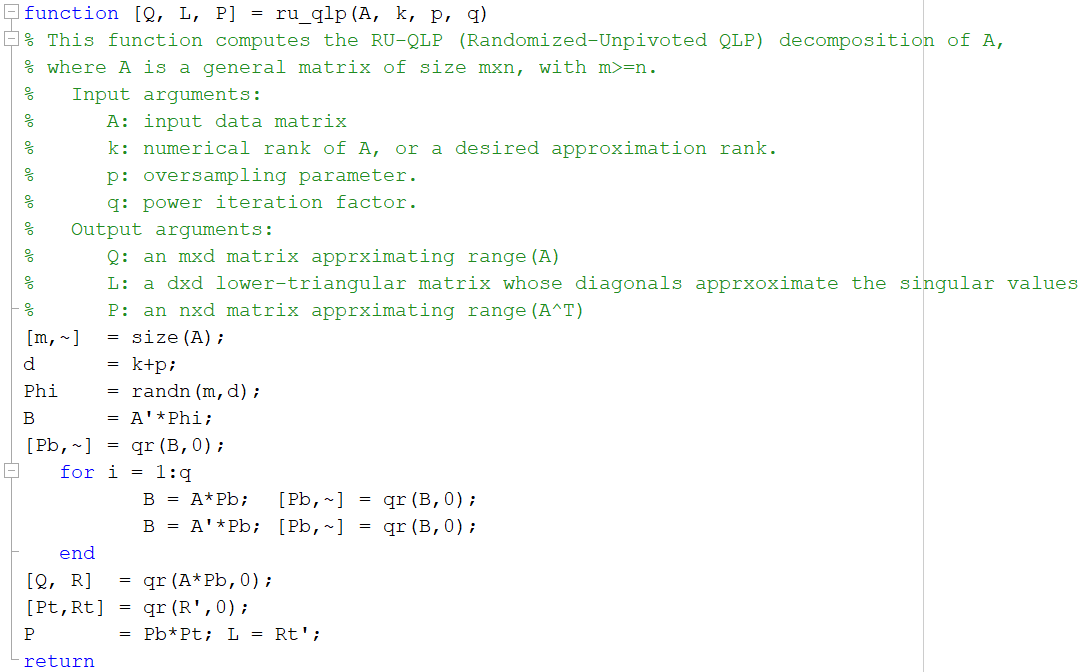}
	\label{fig_Butterfly}
\end{figure}

\begin{figure}[!htb]
	\centering
	\includegraphics[width=0.46\textwidth,height=5cm]{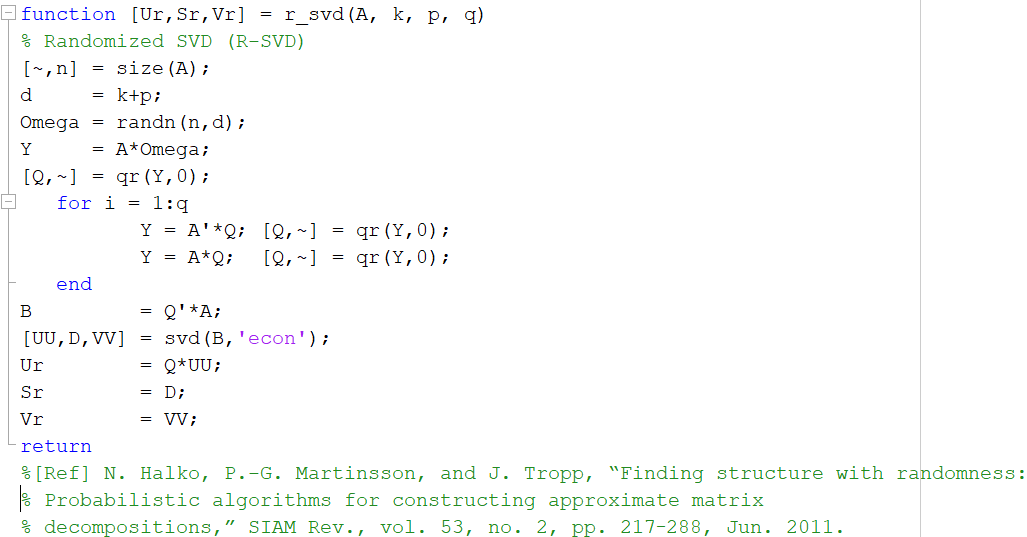}
	\label{fig_Butterfly}
\end{figure}

\begin{figure}[!htb]
	\centering
	\includegraphics[width=0.46\textwidth,height=5.7cm]{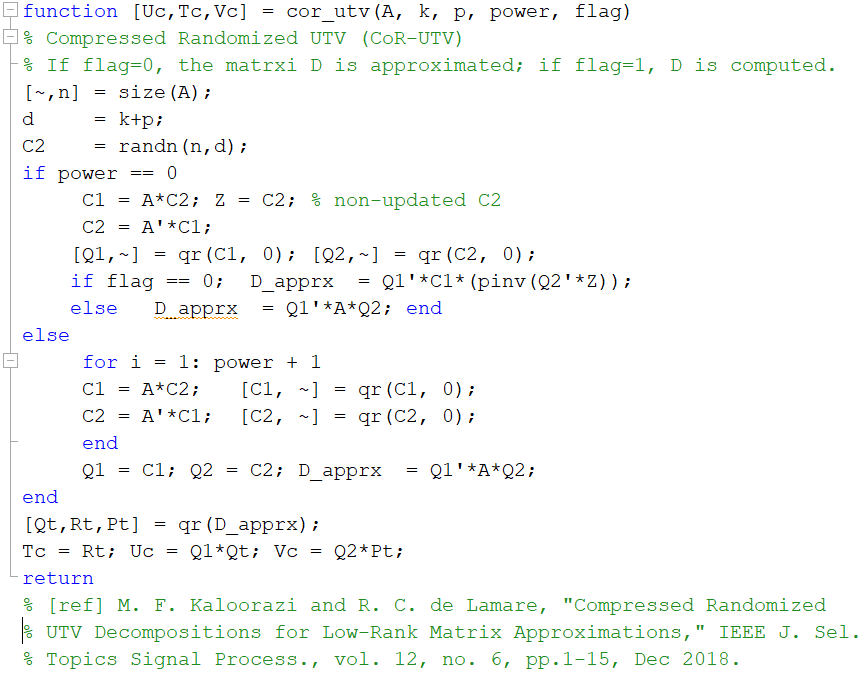}
	\label{fig_Butterfly}
\end{figure}

\begin{figure}[!htb]
	\centering
	\includegraphics[width=0.45\textwidth,height=5.5cm]{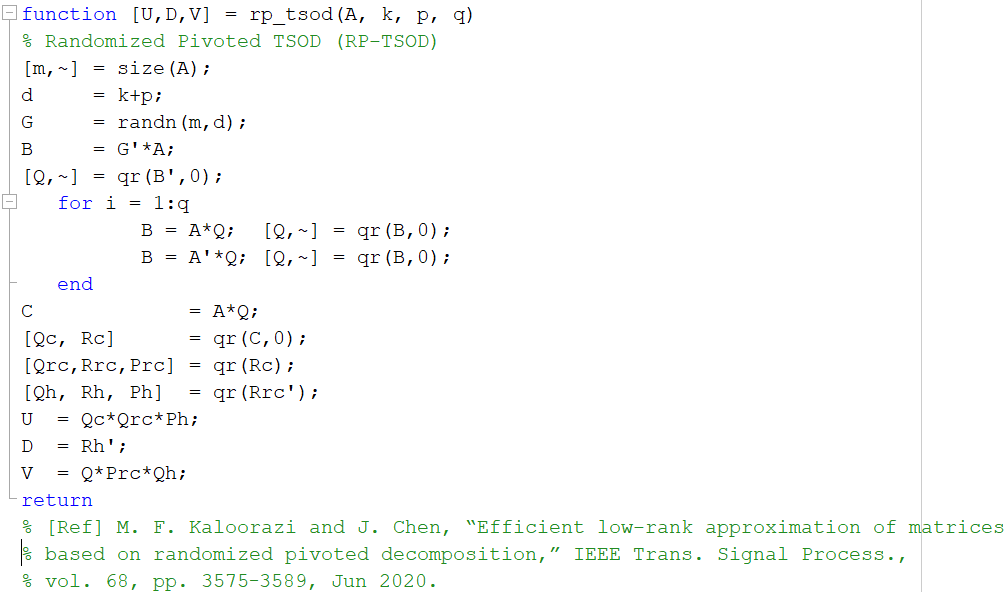}
	\label{fig_Butterfly}
\end{figure}

\begin{figure}[!htb]
	\centering
	\includegraphics[width=0.3\textwidth,height=3cm]{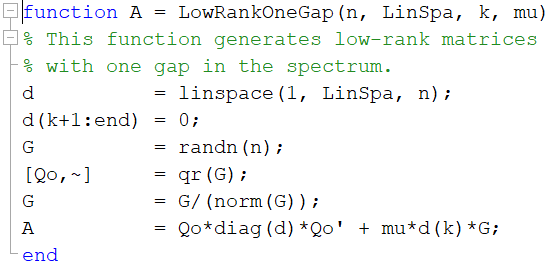}
	\label{fig_Butterfly}
\end{figure}

\begin{figure}[!htb]
	\centering
	\includegraphics[width=0.4\textwidth,height=6.2cm]{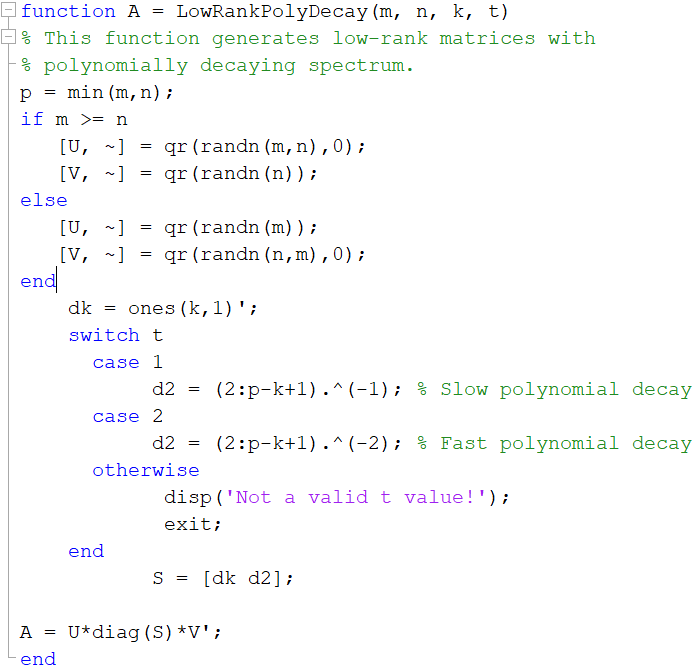}
	\label{fig_Butterfly}
\end{figure}

\EOD

\end{document}